\documentclass[a4 paper,12pt]{article}
\usepackage[T1]{fontenc}
\usepackage[latin1]{inputenc}
\usepackage[german,english]{babel}
\usepackage{graphicx}
\usepackage{mathrsfs}
\usepackage{amsfonts}
\usepackage{amsxtra}
\usepackage{amssymb}
\usepackage{pifont}
\usepackage{eufrak}
\usepackage{pst-all}
\usepackage{multido}
\usepackage{amsmath}
\usepackage{color}
\usepackage{makeidx}
\usepackage{multicol}
\usepackage[absolut]{overpic}

\usepackage[rm,bf,tiny,center]{titlesec}
\titlelabel{\thetitle.\enspace}
\newcommand{\ger}[1]{\mathfrak{#1}}

\newcommand{\OO}{\mathcal{O}}

\newcommand{\sur}{\twoheadrightarrow}
\newcommand{\A}{\mathcal{A}}
\newcommand{\B}{\mathcal{B}}
\newcommand{\M}{\mathcal{M}}
\newcommand{\N}{\mathcal{N}}
\newcommand{\LL}{\mathcal{L}}
\newcommand{\ma}{\leqslant}
\newcommand{\De}{\Delta}

\newcommand{\CC}{\mathbb{C}}

\newcommand{\NN}{\mathbb{N}}
\newcommand{\DD}{\mathcal{D}}
\newcommand{\II}{\mathcal{I}}

\DeclareMathOperator{\BGG}{BGG}

\DeclareMathOperator{\Id}{Id}

\DeclareMathOperator{\Ann}{Ann}
\DeclareMathOperator{\add}{add}
\DeclareMathOperator{\dimm}{dim}
\DeclareMathOperator{\domdim}{dom.dim}
\DeclareMathOperator{\Hom}{Hom}
\DeclareMathOperator{\Kern}{ker}

\DeclareMathOperator{\End}{End}
\DeclareMathOperator{\topp}{top}
\DeclareMathOperator{\rad}{rad}
\DeclareMathOperator{\modd}{mod}
\DeclareMathOperator{\soc}{soc}
\DeclareMathOperator{\id}{id}

\DeclareMathOperator{\im}{im}
\DeclareMathOperator{\spann}{span}

\newtheorem{num}{}[section]
\newtheorem{subnum}{}[subsection]

\usepackage{geometry}

\usepackage{fancyhdr}

\begin{document}

\begin{center}
\begin{large}\textbf{Quasi-hereditary algebras via generator-cogenerators of  local self-injective algebras and   transfer of Ringel duality} \end{large}
\end{center}
\begin{center}
Daiva Pu\v{c}inskait\.{e}
\end{center}

\begin{abstract}
 The dominant dimension of algebras in the class $\textbf{A}$ of   1-quasi-hereditary algebras  introduced in \cite{P} is at least two. By the Morita-Tachikawa Theorem this implies that $\textbf{A}$  is related to  a certain class   $\textbf{B}$ of algebras  via   bimodules  satisfying the   double centralizer condition.  In this paper we specify  the class  $\textbf{B}$ and   the modules over algebras in $\textbf{B}$  connected with $\textbf{A}$. 
 
 The class $\textbf{A}$ is  not closed under taking the Ringel-dual. However the dominant dimension of the  Ringel-dual $R(A)$ of $A\in \textbf{A}$  is at least two. This fact induces  a corresponding concept of modules over $B\in \textbf{B}$ which yield the  algebras $A$ and $R(A)$ for $A\in \textbf{A}$.
\end{abstract}

\begin{center}
\textbf{Introduction}
\end{center}
Let $\A$, $\B$ be algebras.   An $\A$-$\B$-bimodule ${}_{\A}\M_{\B}$  satisfying the   double centralizer condition $\A\cong \End_{\B}(\M_{\B})$ and $\B\cong \End_{\A}({}_{\A}\M)$   provides a relationship  between the representation theories  of the algebras $\A$ and $\B$  that may differ in terms of their  homological properties.  Soergel's 'Struktursatz' relating an algebra $\A_{\Theta}(\ger{g})$ corresponding to  a block $\Theta$ of the Bernstein-Gelfand-Gelfand category $\OO(\ger{g})$  of a complex semisimple Lie algebra $\ger{g}$ with 
a subalgebra of  the corresponding coinvariant algebra as well as the Schur-Weyl duality between the Schur algebra $S(n,r)$ for $n\geq r$ and the group algebra $K\Sigma_r$ of the symmetric group  are prominent examples for this connection (see for instance \cite{KSX}). 

 In this paper we present a further example for  this situation  which relates a 1-quasi-hereditary algebra $A$ defined in \cite{P} with a local self-injective algebra $B$ via an $A$-$B$-bimodule $L$  whose structure has a precise description: The $A$-module $L$ is a   projective-injective  indecomposable  and  the  $ \End_A(L)$-module $L$ decomposes into a  direct sum of  local ideals of $B:=\End_A(L)$ generated by the endomorphisms   corresponding   to 
  certain    paths  in the quiver of $A$ (see Proposition ~\ref{End(P(1))}).  The $B$-module $L$ is a generator-cogenerator of $\modd B$ (i.e., any projective resp. injective indecomposable $B$-module is a direct summand of $L$).   Thus,  any  1-quasi-hereditary algebra can be defined  as an endomorphism algebra of a generator-cogenerator of a local self-injective algebra. 

 The algebras $\A_{\Theta}(\ger{g})$ and  $S(n,r)$ belong  to the class $ \mathscr{A}$    of quasi-hereditary algebras   with a duality (induced  by an anti-automorphism) and with   dominant dimension at least two (see \cite{KSX} and  \cite{FK}). 
The class of 1-quasi-hereditary algebras has a non-empty intersection with  $\mathscr{A}$.  Many   factor algebras of $\A_{\Theta}(\ger{g})$  related to certain  saturated subsets of weights  are 1-quasi-hereditary. Note that  a 1-quasi-hereditary algebra does not have a duality in general. 
  The results in this paper  clarify  the connection between  $\A_{\Theta}(\ger{g})$ and the   coinvariant algebras.

Our first result presents a one-to-one correspondence (up to isomorphism) between the class of 1-quasi-hereditary algebras (over an  algebraically closed field $K$) and  
 the class of pairs $(B,L)$   yielding 1-quasi-hereditary algebras via  the double centralizer condition, where $B$ is a local  self-injective  algebra and $L\in \modd B$ satisfies certain properties.
\\[-3mm]

\text{  }
\\[-2pt]
 \hspace*{0.2cm}\begin{minipage}{6cm}
 \hrule
 \vspace{0.5cm} \textcolor{black}{ }
\end{minipage}
\\[-4mm]
\hspace*{0.3cm}\begin{footnotesize}\rm{Partly supported by the D.F.G. priority program SPP 1388 ``Darstellungstheorie''.}\end{footnotesize}

Any  algebra $\A$ in this paper is  basic, thus given by a quiver and relations $(Q(\A),\II(\A))$. 
\\

\textbf{Theorem A.}
\textit{Let   $A$,  $B$ be   finite dimensional  basic  $K$-algebras and   $L$ be  a $B$-module.  
Let    $n\in \NN$  and    $(\Lambda :=\left\{1,\ldots , n\right\},\ma)$ be a    partially ordered set.
 The following statements are
  equivalent:
\begin{itemize}
	\item[(i)]  $A$  with   $\left(\Lambda , \ma \right)$ is 1-quasi-hereditary 
	(here  we  identify $\Lambda$  with the vertices in $Q(A)$), i.e.,       $A\cong \End_B(L)^{op}$,  where  $L$ is a multiplicity-free  generator-cogenerator of $\modd B$.
	\item[(ii)] $B$ is   local, self-injective  with  $\dimm_KB=n$  and   $\displaystyle L\cong \bigoplus_{i\in \Lambda} L(i)$  where   $L(i)$ are    local  submodules    of $B$ and  $L(1)=B$,  moreover,  for all     $i, j\in \Lambda $   the following properties~  hold:
\begin{itemize}
	\item[(a)] $L(i)\sur L(j)$ if and only if  \ $i\ma j$,
	\item[(b)] $\displaystyle \rad (L(i))=\sum_{i<j}L(j)$.
\end{itemize}
\end{itemize}}

An algebra of the form  $A_{\Theta}(\ger{g})$ is 1-quasi-hereditary if $\text{rank} (\ger{g})\leq 2$, hence Theorem A is applicable for these  algebras.

Dlab,  Heath  and  Marko have shown in \cite{DHM} that a pair $(B,L)$ with the properties in $(ii)$ yields a quasi-hereditary BGG-algebra  (defined by Irving in \cite{Ir}) if $B$ is commutative. The next theorem strengthens the main theorem in  \cite{DHM} by  determining  the properties of a 1-quasi-hereditary algebra  $A\cong \End_B(L)^{op}$ for which  $B$ is commutative. 

In the quiver $Q$ of a 1-quasi-hereditary algebra  $A$ between  two vertices $i$ and $j$  either there are no arrows or  $i \leftrightarrows j$ (see \cite[Theorem 2.7]{P}). Thus for any path $p$ in $Q$ the opposite path $p^{op}$ also belongs to $Q$.
\\

\textbf{Theorem B.}
\textit{Let     $A\cong \End_B(L)^{op}$   be a   1-quasi-hereditary  algebra where the  $B$-module $L$ satisfies  the   conditions $(ii)$ in Theorem A.   
  The following statements are equivalent:
\begin{itemize}
	\item[(i)] $B$ is commutative.
		\item[(ii)] If  $\rho$ is a relation of $A$, then   $\rho^{op}$ is also a  relation of $A$.
		\item[(iii)] $A$ has a duality induced by the  anti-automorphism $p\mapsto p^{op}$.
\end{itemize}}

The coinvariant algebra related to the algebra $\A_{\Theta}(\ger{g})$  is commutative.  Note that Theorem B is also true for  all algebras $\A_{\Theta}(\ger{g})$, where the quiver and relations are known, also for non 1-quasi-hereditary algebras (see \cite{Str}).

The concept of Ringel duality  introduced in \cite{Rin1} is essential  in the theory of quasi-hereditary algebras  (for a (basic) quasi-hereditary algebra $\A$ there exists another quasi-hereditary algebra $R(\A)$ such that $R(R(\A))\cong \A$). In this paper we show how how  the Ringel duality $R(-)$ induces a corresponding concept $\widetilde{R}(-)$ on certain   generator-cogenerators of a local self-injective algebra. 
 This  is based on the fact that the class of 1-quasi-hereditary algebras is  not closed under Ringel duality, however  for any 1-quasi-hereditary algebra $A$ there exists an  $R(A)$-$ \widetilde{R}(B)$-bimodule $ \widetilde{R}(L)$ having the double centralizer condition (see Lemma ~\ref{33}).

Our next result explicitly determines the aforementioned correspondence for those   1-quasi-hereditary algebras, whose Ringel duals are also 1-quasi-hereditary.
\\

\textbf{Theorem C.}
\textit{Let $(A,\ma)$ and $(R(A),\geqslant )$ be 1-quasi-hereditary algebras as well as $(B,L)$  and $( \widetilde{R}(B), \widetilde{R}(L))$  the corresponding pairs (w.r.t. Theorem A \textit{(ii)}),  where  $\displaystyle L\cong \bigoplus_{i\in \Lambda}L(i)$ and  $ \displaystyle \widetilde{R}(L)\cong \bigoplus_{i\in \Lambda}\widetilde{R}(L(i))$. Then $B\cong \widetilde{R}(B)$ and    for every  $i\in \Lambda$  we have 
\begin{center}
$\displaystyle \widetilde{R}(L(i))\cong B/\left(\sum_{j\not\ma i}L(j)\right) \cong \bigcap_{j\not\ma i} \ker\left(B\sur L(j)\right).$
\end{center}}

In particular, if a 1-quasi-hereditary algebra is Ringel self-dual, then $L \cong \widetilde{R}(L)$. There exists a permutation $\sigma\in \text{Sym}(\dimm_KB)$ with $L(\sigma(i))\cong \widetilde{R}(L(i))$. The algebras of the form $\A_{\Theta}(\ger{g})$ are Ringel self-dual. In case of $\text{rank}(\ger{g})\leq 2$ we have this situation.
\\

The paper is organized   as follows: In Section 1, we introduce the Morita-Tachikawa  Theorem  which shows  that a minimal faithful module over  an algebra of dominant dimension at least two    has the  double centralizer property. The results of this paper rest on this theorem. 
We also recall the  relevant definitions and give some examples 
which show  the diversity of modules over  a local self-injective  algebra satisfying  the double centralizer condition.

 Section 2 is devoted to the proof of Theorem A. The paths in the quiver  of a 1-quasi-hereditary algebra of  the  form $p(j,i,k)$  defined in \cite[Section 3]{P} play  an  important role. The  other part  of the proof   is based      on the structure of  the  $B$-module $L$ which  will be analyzed    in Lemma ~\ref{L5}.
 We also determine   an easier transition  from  the $B$-maps of   $L$     to the  relations of the algebra $A=\End_B(L)^{op}$ (see Remark~\ref{streles}).
 This is used in the proof of Theorem B in Section 3.
  
 In   Section 4 we describe the transfer of  Ringel duality. Consequently we obtain a decomposition of the  class  of 1-quasi-hereditary algebras with  their  Ringel duals into  subclasses  which are  closed under Ringel duality. Moreover, the algebras in a fixed subclass  arise  from the same local self-injective algebra.  Subsequently we present the proof of Theorem C.

\section{Preliminaries}
\begin{small}
Unless otherwise specified,  any algebra $\A$   is an  associative,   finite dimensional,  basic   $K$-algebra  over an algebraically closed field $K$. Thus $\A$ is determined  by a quiver $Q$ and relations $\II$, i.e.,  $\A=KQ/\II$ (Theorem of Gabriel). Furthermore,   $\modd \A$   is the category  of finite dimensional left  $\A$-modules.  For  $\M\in \modd \A$ we denote by $\soc \M$, $\rad \M$ and  $\topp \M$ the socle, the  radical and  the top of $\M$, respectively, and  $\add (\M)$  is the full subcategory of $\modd \A$ whose objects are the
direct sums of direct summands of  $\M$ (for references see \cite{ARS} or \cite{ASS}).  To distinguish between an  arbitrary  algebra $\A$ and a 1-quasi-hereditary algebra, we denote the latter by $A$.
 \end{small}
 \\
 
 We repeat  some notations and  facts  about  bound quiver algebras  $\A=KQ/\II$. Throughout, we identify  the set of   vertices $Q_0:=Q_0(\A)$    with the set  $\Lambda=\left\{1, \ldots , \left|Q_0\right|\right\}$. The product of arrows $(k\to i)$  and
$(i\to j)$ is given by  $(k\to i\to j)=(i \to j )\cdot(k \to i)$. We denote 
  by   $P(i)$, $I(i)$,  $S(i)$ and $e_i$  the    projective indecomposable,  injective indecomposable,    simple $\A$-module   and  the  primitive idempotent, respectively, corresponding to  $i\in \Lambda$.    Let $\M\in \modd \A$, then for each vertex $i$ there exists a $K$-subspace of $\M$ corresponding to $i$, denoted by  $\M_i$. We have    $\M_i\cong \Hom_{\A}(P(i),\M)$ via $m\mapsto f_{(m)}:P(i)\to \M$, where  $f_{(m)}(a\cdot e_i)=a\cdot m$ for all $a\in \A$. Recall that     $\End_{\A}(\M)^{op}$ is a $K$-algebra with    product $F\circ G=\left(\M\stackrel{G}{\rightarrow}\M \stackrel{F}{\rightarrow}\M\right)$. 
 The  (left) $\End_{\A}(\M)^{op}$-module  $\M$  (written  ${}_{\End_{\A}(\M)^{op}}\M$) is isomorphic to   $\bigoplus_{i\in \Lambda}\Hom_{\A}(P(i),\M)$.  The   Jordan-Hölder multiplicity  of $S(i)$ in $\M$ is  denoted  by  $[ \M:S(i) ]$. In particular, we have $\dimm_K\Hom_{\A}(P(i),\M)=[ \M:S(i) ]$.

\subsection{Theorem of  Morita-Tachikawa}
  Based on various works  by Morita and Tachikawa,   Ringel has
described  in  \cite{Rin2}   a relationship between the  algebras $\A$  and $\End_{\A}(\M)^{op}$  via an  $\A$-module $\M$ having  the   double centralizer condition. 
  We recall some notations and terminology needed in the statement of the  theorem below. 
  Let $\B$ be an artin (basic) algebra,  then   the   \textit{dominant
dimension} of $\B$  is at least $2$  (written $\domdim \B\geq 2$),  if  there exists   an exact sequence
$0\to {}_{\B}\B \rightarrow \M_1 \to  \M_2$
such that  $\M_1,\M_2$    are projective and injective  $\B$-modules  (all finite dimensional algebras are artin algebras).
 A  $\B$-module  
   $\N$ is  called      \textit{minimal faithful} if $\N$ is  faithful (i.e.,   ${}_{\B}\B$ can be
embedded into  a direct sum of copies of $\N$)  and  $\N$  is  a direct summand of any faithful $\B$-module.   A minimal faithful $\B$-module is   unique (up to isomorphism) and will be  denoted   by  $\mathscr{F}(\B)$.  A  $\B$-module  $\M$  is a \textit{generator-cogenerator} of $\modd \B$  if 
    every projective indecomposable  as well as  every  injective indecomposable $\B$-module is  a  direct summand of $\M$.  
 We denoted by  $[\B]$  resp.   $[\B,\M]$   the isomorphism class  of  $\B$  and  a $\B$-module $\M$.

\begin{num}\begin{normalfont}\textbf{Theorem  (Morita-Tachikawa)}.\end{normalfont} \textit{There is are bijections  $\Psi$  and $\Phi$   between 
\begin{itemize}
	\item[] $ \textbf{X} := \left\{ \ \ \ [\A] \ \ \  \    \left| 
\begin{array}{ll}
\A \text{ is a basic artin  algebra, }  
\domdim \A \geq 2
\end{array}
\right.\right\} $  \   and 
	\item[] $\textbf{Y}:= \left\{ \ [\B,\M]  \  \left|   
\begin{array}{l}
\B \text{  is a basic artin algebra,  }\\
\M \text{  is a multiplicity-free,  
 generator-cogenerator of }\modd \B
\end{array}
\right.\right\} $ \end{itemize}
defined as follows: 
\\[7pt]
$
\begin{array}{ccl}
\textbf{X} & \stackrel{\Phi}{\longrightarrow} & \textbf{Y} \\
\left[ \A\right] &  \mapsto &  \bigl[ \B(\A):=\End_{\A}\left(\mathscr{F}(\A)\right)^{op},{}_{\B(\A)}\mathscr{F}(\A)\bigr]
\end{array}
$ \ \  and  \ \  $
\begin{array}{ccl}
\textbf{Y} & \stackrel{\Psi}{\longrightarrow} & \textbf{X} \\
\left[ \B,\M\right] &  \mapsto &  \left[\End_{\B}(\M)^{op}\right]
\end{array}
$,
\\[7pt]
such that    $\Psi\circ \Phi = \Id_{\textbf{X}}$  and $\Phi\circ \Psi = \Id_{\textbf{Y}}$.}
\label{Morita}
\end{num}

This theorem  also  provides    correspondences between the subsets of $\textbf{X}$ and their  image under  $\Phi$ in  $\textbf{Y}$. By the Theorem of König, Slung\r{a}rd and Xi \cite[Theorem 1.3]{KSX}  the set of  isomorphism classes of algebras which corresponds to  the blocks of the BGG-category $\OO$ is a subset of   $\textbf{X}$. Theorem ~\ref{Morita} restricted to this subset is known as Soergel's 'Struktursatz'. Moreover, the dominant dimension of  Schur algebras $S(n,r)$ (with $n>r$) is at least two. The Schur-Weyl duality  is a special case of  Theorem ~\ref{Morita}. The module $ \mathscr{F}(\A)$ has  the   double centralizer condition.

Note that  the Theorem of Morita-Tachikawa provides some  connections   between a finite dimensional algebra $\A=KQ/\II$ with $\domdim \A\geq 2$ and a pair $(\B,\M)$ with $\Phi[A]=[\B,\M]$: The $\B$-module $\M$ has $\left|Q_0\right|$  pairwise non-isomorphic,  indecomposable direct  summands (they correspond to the  vertices of $Q$). Because $\A\hookrightarrow \mathscr{F}(\A)^m$ for some $m\in \NN$, we have  $[\mathscr{F}(\A):S(i)]\neq 0$ for all $i\in Q_0$, thus $\Hom_{\A}(P(i),\mathscr{F}(\A))\neq 0$.  Since the  $\B$-module  $\M\cong \mathscr{F}(\A)$ is isomorphic to $\bigoplus_{i\in Q_0}\Hom_{\A}(P(i),\mathscr{F}(\A))$, we obtain that the  $\End_{\A}\left(\mathscr{F}(\A)\right)^{op}$-module  $\Hom_{\A}(P(i),\mathscr{F}(\A))$ is indecomposable for any $i\in Q_0$.

\subsection{Quasi-hereditary algebras and local self-injective algebras}
In this paper we consider a subclass of quasi-hereditary algebras  with dominant dimension at least 2 which  are related to   local self-injective algebras. 
We recall some necessary definitions.

 Quasi-hereditary algebras  were  defined by Cline, Parshall and Scott in  \cite{CPS}. We use the equivalent definition and terminology  given by Dlab and Ringel in \cite{DRin}: Let   $\A$  be a bound quiver algebra  and  $(\Lambda,\ma)$  a  poset (with $\Lambda=Q_0(\A)$).   
For  every   $i\in \Lambda$ the  \textit{standard}  module    $\De(i)$  is   the largest  factor module    of     $P(i)$  such that      $[\De(i):S(k)]=0$ for all $k\in \Lambda$ with  $k\not\ma i$. 
The modules in  the full subcategory  $\ger{F}(\De)$  of $\modd \A$  consisting of the modules having a filtration  such that each subquotient is isomorphic to a standard module  are called    \textit{$\De$-good}  and  these   filtrations are   \textit{$\De$-good filtrations}.  For   $M\in \ger{F}(\De)$  we   denote  by $\left(M:\De(i)\right)$  the (well-defined) number of  subquotients  isomorphic to  $\De(i)$ in a $\De$-good filtration  of $M$. 
 \\[8pt]
\hspace*{2mm}\parbox{16cm}{
An algebra $\A$  with $(\Lambda,\ma)$  is  \textit{quasi-hereditary}    if  for all $i,k\in \Lambda$ the following    conditions are satisfied:
\begin{itemize}
	\item   $[\De(i):S(i)]=1$,
	\item $P(i)$ is a  $\De$-good module   with $\left(P(i):\De(k)\right)=0$ for all $k\not\geqslant i$ and $\left(P(i):\De(i)\right)=1$.
\end{itemize}
     }
     \\ 
      Throughout, $(\A,\ma)$ denotes  an algebra $\A$  with a partial order  $\ma$ on the vertices $\Lambda = Q_0(\A)$.
  
 We can identify the vertices of the quivers of $\A$ and $\A^{op}$.
     An algebra $(\A,\ma)$ is quasi-hereditary if and only if $\left(\A^{op},\ma\right)$ is quasi-hereditary (see \cite{DRin}).  The standard duality $\DD:=\Hom_K(-,K)$  provides   the \textit{costandard} $\A$-module $\nabla(i)\cong \DD(\De_{\A^{op}}(i))$ corresponding to   $i\in \Lambda$  and also   the  subcategory $\ger{F}(\nabla)$  of $\modd \A$  of all $\nabla$-\textit{good}  modules.

\begin{subnum}  \normalfont{
\textbf{Remark.}} Let $(\A,\ma)$  be a quasi-hereditary algebra  and let   $M$  be  a projective-injective $\A$-module.  Moreover, let    $\soc (\De(i))\in \add \left(\soc M\right)$ and  $P(i)\hookrightarrow M$ with  $M/P(i)\in \ger{F}(\De)$ for all $i\in \Lambda$. In this case we have $\soc\left(M/P(i)\right)\in \add\left(\bigoplus_{j\in \Lambda}\soc \De(j)\right)\subseteq \add \left(\soc M\right)$ and consequently $M/P(i)$ can  be embedded into some copies of $M$ for any $i\in \Lambda$. In other words, there exist $m,r\in \NN$ such that the sequence $0\to {}_{\A}\A\to M^m\to M^r$ is exact, i.e.,
 we have $\domdim \A\geq 2$. 
\label{dd2r}    
 \end{subnum}

    We  recall the definition of a 1-quasi-hereditary algebra from  \cite{P}  and show that we have the situation described in the foregoing  remark.

\begin{subnum}  \normalfont{
\textbf{Definition (1-quasi-hereditary).}}
A quasi-hereditary algebra $A$ with $(\Lambda,\ma)$
is called 
\textit{1-quasi-hereditary}
if for all $i,j\in \Lambda=\left\{1, \ldots , n\right\}$ the
following conditions are satisfied:
\begin{itemize}
    \item[(1)] There is a smallest and a largest element with respect to  $\ma$, 
    \\ without loss of generality we will assume  them to be $1$ resp. $n$,
    \item[(2)] $[\Delta(i):S(j)]=\big(P(j):\Delta(i)\big)=1$ for $j\leqslant i$,  
\item[(3)] $\soc P(j) \cong \topp I(j)\cong  S(1)$,
\item[(4)] $\Delta(i) \hookrightarrow \Delta(n)$ and $\nabla(n)\twoheadrightarrow
\nabla(i). $
\end{itemize} \label{def1qh} 
\end{subnum} 

\begin{subnum} \begin{normalfont}\textbf{Lemma.}\end{normalfont} \textit{Let $(A,\ma)$ be  a 1-quasi-hereditary algebra  with $\left\{1\right\}=\min (\Lambda,\ma)$. Then    $P(1)$ is a minimal faithful $A$-module and $\domdim A\geq 2$.  
}
\label{mini}
 \end{subnum} 
    
\textit{Proof.} According
to  \cite[2.6]{P} we have  $P(1)\cong I(1)$  and $\De(j)\hookrightarrow P(i)\hookrightarrow P(1)$  for every $i,j\in \Lambda$. We obtain ${}_AA\hookrightarrow P(1)^{\left|\Lambda\right|}$, thus   $P(1)$ is a minimal faithful $A$-module because  $P(1)$    is indecomposable (we have  $\mathscr{F}(A)\cong P(1)$). 

Since    $\soc \De(i)\cong \soc P(1)\cong S(1)$  and $P(1)/P(i)\in \ger{F}(\De)$ (see \cite[4.3]{P}),  we have $\domdim A\geq 2$ according to   Remark~\ref{dd2r}.
\hfill $\Box$ 
\\

 An  (finite dimensional, basic)  algebra $\B$ is local and  self-injective   if and only if  the socle and the top of $\B$-module $\B$ are simple. An ideal $I$ of $\B$  is  a two-sided,  local ideal  if   $\B$$ \cdot$$ I \subseteq I$ as well as  $I$$\cdot$$ \B\subseteq I$ and  $\rad (I)$ is  the (uniquely  determined)  maximal submodule of $I$.

\begin{subnum}  \normalfont{
\textbf{Definition (\fbox{$\ma$}).}} Let $B$ be an algebra,  $L\in \modd B$ and $\left(\Lambda=\left\{1,\ldots,n\right\},\ma\right)$  be a poset.  We say  that the  pair $(B,L)$   \textit{satisfies  the
    condition  \fbox{$\ma$}} if  
\begin{itemize}
	\item[(1)]  $B$ is a  local, self-injective algebra,  $\dimm_K B=n$,
	\item[(2)]  $\displaystyle L\cong \bigoplus_{i\in \Lambda}^{} L(i)$, where   $L(1),\ldots , L(n)$ are  two-sided  local  ideals  of $B$  with   $L(1)=B$  and for all     $i, j\in \Lambda $   the following is satisfied:
	\begin{itemize}
	\item[(a)] $L(i)\sur L(j)$ if and only if  \ $i\ma j$,
	\item[(b)] $\displaystyle \rad \left(L(i)\right)=\sum_{i<j}L(j)$.
\end{itemize}
\end{itemize}
	   \label{property_ma} 
	\end{subnum} 
	
\begin{subnum}  \normalfont{
\textbf{Remark}.}	 If $(B,L)$ satisfies  the   condition  \fbox{$\ma$}, then $[B,L]\in \textbf{Y}$:  Since $B$ is local and self-injective, any projective resp.  injective,  indecomposable $B$-module is isomorphic to $_{B}B$.  Hence $L$ is a generator-cogenerator of $\modd B$,  because $B (=L(1))$ is a direct summand of $L$. The condition $(a)$ implies $L(i)\cong L(j)$ if and only if  $i=j$, therefore $L$ is multiplicity-free. Moreover, the quiver of $B$ consists of one vertex and  finitely many loops ($B$ is isomorphic to  $K\left\langle x_1,\ldots, x_r\right\rangle /J$, where $J$ is an  ideal  with $J\subseteq \left\langle  x_1,\ldots, x_r\right\rangle^2$  and    $r\in \NN$). 
\label{BLinY} 
	\end{subnum}  
It should be noted that for a fixed local, self-injective  algebra $B$ a poset $(\Lambda,\ma)$ and a $B$-module $L$ with   the   condition  \fbox{$\ma$}  are not uniquely determined.

\begin{subnum} \normalfont{
\textbf{Example}.} The algebra $B=\CC[x,y]/\left\langle xy, \ x^3-y^3\right\rangle$ is local and self-injective with $\dimm_{\CC}B=6$. The next diagrams present  the partial orders  $\ma_{(i)}$ on $\Lambda=\left\{1,\ldots , 6\right\}$ and  the  generators-cogenerators $L_i=\bigoplus_{k=1}^6L_i(k)$ of $\modd B$ for $i=1,2,3.$ We  write $  L_{i}(k)\to  L_{i}(k')$ if  $ L_{i}(k) \supset  L_{i}(k')$ and 
there does not exist a submodule $U$ with $ L_{i}(k) \supset U \supset L_{i}(k')$. In other words, 
$k <_{(i)}k'$ and $k,k'$ are neighbours  for all $1\leq i\leq 3$ and $1\leq k,k'\leq 6$. In the diagram on the right-hand side,   $\texttt{i}\in \CC$ is such that  $\texttt{i}^2=-1$ and $\omega = \frac{1}{2}+\texttt{i}\frac{\sqrt{3}}{2}$ is a  6th  root of unity.
\label{exampleB}
\end{subnum} 

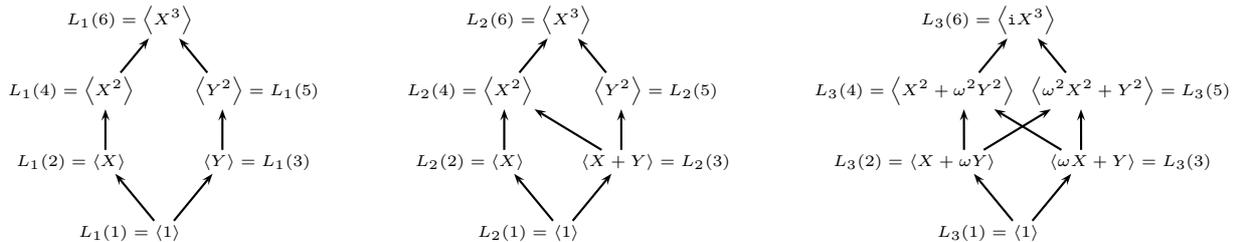
\text{ }\\[-7mm]
   
    \hspace{5mm}  \psset{xunit=0.77mm,yunit=0.94mm,runit=1.2mm}
\begin{pspicture}(-10,-17)(0,0)
\begin{tiny}
              \rput(0,0){\rnode{1}{$L_1(6)=\left\langle X^3\right\rangle \hspace{9mm}$}}
              \rput(-10,-10){\rnode{2}{$L_1(4)=\left\langle X^2\right\rangle \hspace{9mm}$}}
              \rput(10,-10){\rnode{3}{$\hspace{9mm}\left\langle Y^2\right\rangle = L_1(5)$}}
              \rput(-10,-20){\rnode{4}{$L_1(2)= \left\langle X\right\rangle \hspace{9mm}$}}
              \rput(10,-20){\rnode{5}{$\hspace{9mm} \left\langle Y\right\rangle = L_1(3)$}}
              \rput(0,-30){\rnode{6}{$L_1(1)=\left\langle 1\right\rangle \hspace{9mm} $}} \end{tiny}
              \psset{nodesep=1.5pt}
              \ncline{<-}{1}{2}
              \ncline{<-}{1}{3}
              \ncline{<-}{2}{4}
              \ncline{<-}{3}{5}
              \ncline{<-}{4}{6}
              \ncline{<-}{5}{6}
\end{pspicture} \hspace{42mm}  \psset{xunit=0.77mm,yunit=0.94mm,runit=1.2mm}
\begin{pspicture}(-10,-17)(0,0)
\begin{tiny}
              \rput(0,0){\rnode{1}{$L_2(6)=\left\langle X^3\right\rangle \hspace{9mm}$}}
              \rput(-10,-10){\rnode{2}{$L_2(4)=\left\langle X^2\right\rangle \hspace{9mm}$}}
              \rput(10,-10){\rnode{3}{$\hspace{9mm}\left\langle Y^2\right\rangle = L_2(5)$}}
              \rput(-10,-20){\rnode{4}{$L_2(2)= \left\langle X\right\rangle \hspace{9mm}$}}
              \rput(10,-20){\rnode{5}{$\hspace{9mm} \left\langle X+Y\right\rangle = L_2(3)$}}
              \rput(0,-30){\rnode{6}{$L_2(1)=\left\langle 1\right\rangle \hspace{9mm} $}} \end{tiny}
              \psset{nodesep=1.5pt}
              \ncline{<-}{2}{5}
              \ncline{<-}{1}{2}
              \ncline{<-}{1}{3}
              \ncline{<-}{2}{4}
              \ncline{<-}{3}{5}
              \ncline{<-}{4}{6}
              \ncline{<-}{5}{6}
\end{pspicture} \hspace{50mm}  \psset{xunit=0.77mm,yunit=0.94mm,runit=1.2mm}
\begin{pspicture}(-10,-17)(0,0)
\begin{tiny}
              \rput(0,0){\rnode{1}{$L_3(6)=\left\langle \texttt{i}X^3\right\rangle \hspace{9mm}$}}
              \rput(-10,-10){\rnode{2}{$L_3(4)=\left\langle  X^2+\omega^2Y^2\right\rangle \hspace{13mm}$}}
              \rput(10,-10){\rnode{3}{$\hspace{13mm}\left\langle \omega^2X^2+ Y^2\right\rangle = L_3(5)$}}
              \rput(-10,-20){\rnode{4}{$L_3(2)= \left\langle X +\omega Y\right\rangle \hspace{13mm}$}}
              \rput(10,-20){\rnode{5}{$\hspace{13mm} \left\langle \omega X+Y\right\rangle = L_3(3)$}}
              \rput(0,-30){\rnode{6}{$L_3(1)=\left\langle 1\right\rangle \hspace{9mm} $}} \end{tiny}
              \psset{nodesep=1.5pt}
              \ncline{<-}{1}{2}
              \ncline{<-}{2}{5}
              \ncline{<-}{3}{4}
              \ncline{<-}{1}{3}
              \ncline{<-}{2}{4}
              \ncline{<-}{3}{5}
              \ncline{<-}{4}{6}
              \ncline{<-}{5}{6}
\end{pspicture} 
\\[1.6cm]
\label{examp} 
It is easy to check that $(B,L_i)$ satisfies the condition  \fbox{$\ma_{(i)}$} and therefore  $[B,L_{i}]\in \textbf{Y}$.
In view of Theorem A,  $L_{i}$ is an $A_i$-$B$-bimodule, where $A_i=\End_B(L_{i})^{op}$ is a 1-quasi-hereditary algebra for   $i=1,2,3$ (the quiver and relations of $A_i$ are given in \ref{Examp_relat}) .  
Note that  the  algebra $A_3$  is  associated to a regular  block  of the $\BGG$-category $\OO(\ger{sl}_3)$. 
\\

In the next  section we prove Theorem A  which can be rewritten   as follows:
\\

\begin{normalfont}\textbf{Theorem A.}\textit{
Let   $A$,  $B$ be    finite dimensional  basic  $K$-algebras  and  $n\in \NN$. Moreover, let     $\left(\Lambda :=\left\{1,\ldots , n\right\},\ma\right)$ be    partially ordered.
 The following statements are equivalent:
\begin{itemize}
	\item[(i)] The algebra  $(A,\ma)$  is 1-quasi-hereditary, i.e.,     $A\cong \End_B(L)^{op}$  and    $L$ is a multiplicity-free  generator-cogenerator of $\modd B$,
	\item[(ii)]  The  pair $(B,L)$   satisfies  the
   condition  \fbox{$\ma$}.
\end{itemize}}
\label{theorem}
\end{normalfont} 

For   any   1-quasi-hereditary algebra $A$ there exists an (up to isomorphism)   uniquely determined  algebra $B$ and  a multiplicity-free generator-cogenerator  $L$  of $\modd B$ with $A\cong \End_B(L)^{op}$ (see   Theorem ~\ref{Morita}  and Lemma ~\ref{mini}). According to Remark ~\ref{BLinY}, 
Theorem A  provides  bijections   between  the isomorphism classes of 1-quasi-hereditary algebras and the pairs   defined in ~\ref{property_ma}.

\subsection{BGG-algebras}
 Dlab, Heath and  Marko  have shown in \cite{DHM} that if  for a commutative algebra $B$ and a $B$-module $L$ the pair  $(B,L)$  satisfies  the
   condition  \fbox{$\ma$}, then  the algebra $\End_B(L)^{op}$ is a $\BGG$-algebra as defined  by Irving in  \cite{Ir}. Our Theorem B  elaborates  on  the Theorem in \cite{DHM}.

We refer to the definition of $\BGG$-algebras given by Xi  in \cite{Chan} (these  algebras are also $\BGG$-algebras in the sense of \cite{Ir}):
 A quasi-hereditary algebra $\A$ is called a BGG-\textit{algebra} if there is a duality $\delta$ of $\modd \A$ such that $\delta$ induces a $K$-linear map on  $\Hom_{\A}(M,N)$ for all $M,N\in \modd \A$ 
and  $\delta(S(i))\cong S(i)$ for all $i\in Q_0(\A)$.

If there is an  anti-automorphism  $\epsilon$ of $\A$ (i.e.,   a $K$-map  $\epsilon:\A\to \A$  with $\epsilon (a\cdot a')=\epsilon(a')\cdot \epsilon(a)$ and $\epsilon^2(a)=a$ for all $a,a'\in \A$) such that $\A\cdot \epsilon(e_i)\cong \A\cdot e_i$ for all $i\in Q_0(\A)$, then $\A$ is a $\BGG$-algebra (see \cite{Chan}).

According to  \cite[Theorem 2.7]{P}, the quiver $Q(A)$ of  a 1-quasi-hereditary  algebra $(A,\ma)$  is the  double    of the  quiver of the  incidence algebra of  $(\Lambda,\ma)$: Let $i,j\in \Lambda$, we  write 
\\[3pt]
 \parbox{13.7cm}{
\begin{center}
$i \triangleleft j$ \ \  and  \ \   $i\triangleright j$   
\end{center}
   if $i$ is a small neighbour    of $j$ and  $i$ is a large  neighbour  of $j$  w.r.t.  $\ma$,  respectively. We have $\left|\left\{\alpha\in Q_1(A)\mid i\stackrel{\alpha}{\rightarrow}j\right\}\right|$$=$$\begin{footnotesize} \left\{
\begin{array}{cl}
1 & \text{if } i\triangleleft j, \\
1 & \text{if } i\triangleright j, \\
0 & \text{else,} 
\end{array}
\right. \end{footnotesize}$ thus  for any  path~ $p=$ 
}
\psset{xunit=0.48mm,yunit=0.53mm,runit=1mm}
\begin{pspicture}(0,0)(0,0)
\rput(30,21){\begin{tiny}
\rput(0,10){\rnode{A}{$n$}}
\rput(0,-44){\rnode{B}{$1$}}
\rput(20,-17){\rnode{aa}{}}
\rput(-20,-17){\rnode{bb}{}}
\rput(-11,3){\rnode{A1}{$t_1$}}
\rput(-11,-37){\rnode{B1}{$l_1$}}
\rput(11,3){\rnode{A2}{$t_{v}$}}
\rput(11,-37){\rnode{B2}{$l_w$}}
              \rput(-10,-7){\rnode{3}{$i_1$}}
              \rput(0,-7){\rnode{33}{$\cdots $}}
              \rput(10,-7){\rnode{4}{$i_r$}}
              \rput(31,-7){\rnode{nx}{$i\triangleright j$}}
              \rput(0,-17){\rnode{5}{$j$}}
              \rput(-10,-27){\rnode{6}{$k_1$}}
              \rput(0,-27){\rnode{66}{$\cdots $}}
              \rput(10,-27){\rnode{7}{$k_{m}$}}
              \rput(31,-27){\rnode{ny}{$k\triangleleft j$}}
              \end{tiny}
                      }
                      
                       \psset{nodesep=1.5pt,arrows=<-}
                      \ncarc[arcangle=15]{4}{nx}
                      \ncarc[arcangle=-15]{7}{ny}
                      
                      \psset{nodesep=1pt,offset=1.3pt,arrows=<-}    
              \ncarc[arcangle=35,linestyle=dotted]{B1}{1}
              \ncarc[arcangle=-25,linestyle=dotted]{A}{3}
              \ncarc[arcangle=25,linestyle=dotted]{3}{A}
              \ncarc[arcangle=-25,linestyle=dotted]{6}{B}
              \ncarc[arcangle=25,linestyle=dotted]{B}{6}
              
              \psset{nodesep=1pt,offset=1.3pt,arrows=->}
              \ncarc[arcangle=35,linestyle=dotted]{2}{B2}
              \ncarc[arcangle=-35,linestyle=dotted]{B2}{2}
              \ncarc[arcangle=-25,linestyle=dotted]{4}{A}
              \ncarc[arcangle=25,linestyle=dotted]{A}{4}
              \ncarc[arcangle=25,linestyle=dotted]{7}{B}
              \ncarc[arcangle=-25,linestyle=dotted]{B}{7}
                      
             \psset{nodesep=1pt,offset=1.3pt}
               \ncline{3}{5}
               \ncline{5}{3}
               \ncline{5}{6}
               \ncline{6}{5}
               \ncline{5}{7}
               \ncline{7}{5}
               \ncline{5}{4}
               \ncline{4}{5}
              \ncline[linestyle=dotted]{21}{22}
              \ncline[linewidth=0.7pt]{->}{A}{A1}
              \ncline[linewidth=0.7pt]{->}{A1}{A}
              \ncline[linewidth=0.7pt]{->}{B1}{B}
              \ncline[linewidth=0.7pt]{->}{B}{B1}
              \ncline[linewidth=0.7pt]{->}{A}{A2}
              \ncline[linewidth=0.7pt]{->}{A2}{A}
              \ncline[linewidth=0.7pt]{->}{B2}{B}
              \ncline[linewidth=0.7pt]{->}{B}{B2}
              \ncarc[arcangle=25,linestyle=dotted]{-}{A2}{aa}
              \ncarc[arcangle=-25,linestyle=dotted]{->}{aa}{A2}
              \ncarc[arcangle=25,linestyle=dotted]{->}{aa}{B2}
              \ncarc[arcangle=-25,linestyle=dotted]{-}{B2}{aa}
              \ncarc[arcangle=-25,linestyle=dotted]{-}{A1}{bb}
              \ncarc[arcangle=-25,linestyle=dotted]{->}{bb}{B1}
              \ncarc[arcangle=25,linestyle=dotted]{-}{B1}{bb}
              \ncarc[arcangle=25,linestyle=dotted]{->}{bb}{A1}
\end{pspicture} 
\\
$(i_1\to i_2 \to \cdots \to i_m)$   there exists a   uniquely determined  path     $\ger{o}(p):=(i_m \to \cdots \to i_2 \to  i_1)$
\\[5pt]
in $Q(A)$ running through the same  vertices  in the  opposite direction. 
 Obviously, $\ger{o}(\ger{o}(p))=p$.

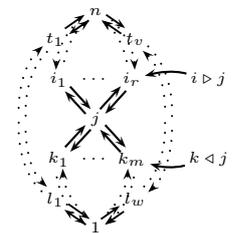
\begin{subnum}\begin{normalfont}\textbf{Definition.}  A 1-quasi-hereditary algebra  $A=KQ/\II$ is a  \textit{$ \BGG_{(\leftrightarrows)}$ -algebra}, if   the $K$-linear map  $\epsilon : KQ \rightarrow KQ$  given by    $\epsilon (p)=\ger{o}(p)$ for all paths $p$ in $Q$  induces 
 an  anti-automorphism  of  $A$. 
 \end{normalfont} 
\end{subnum}

\textbf{Theorem B.} \textit{Let     $A\cong \End_B(L)^{op} \cong KQ/\II$  with $(\Lambda,\ma)$  be a   1-quasi-hereditary  algebra such that $(B,L)$ satisfies  the   condition  \fbox{$\ma$}.  
  The following statements are equivalent:}\textit{
\begin{itemize}
	\item[(i)]   $B$ is commutative. 
	\item[(ii)] $A$  is a $\BGG_{(\leftrightarrows)}$-algebra.
	\item[(iii)]  $\displaystyle \sum_{t=1}^{r}c_t\cdot p_t\in \II$ \    if and only if \   $\displaystyle\sum_{t=1}^{r}c_t \cdot \ger{o}(p_t) \in \II$.
\end{itemize}}
\label{commut}
\label{BGG-algebra}

\text{  }
\\
The Theorems A and B are special cases  of Morita-Tachikawa Theorem ~\ref{Morita}.  For the subsets 
\\[5pt]
$\textbf{X}(1):=\left\{[A]\mid A \text{ is a 1-quasi-hereditary algebra}\right\}$, \ \   $\textbf{X}':=\left\{[A]\mid A\text{ is a }\BGG_{(\leftrightarrows)}\text{-algebra}\right\}$\\[3pt]
\hspace*{7cm}  and   \\[3pt]
 $\textbf{Y}(1):=\left\{[B,L]\mid (B,L) \text{ has the property  \fbox{$\ma$}}\right\}$,\  \   $\textbf{Y}':= \left\{[B,L]\in \textbf{Y}(1) \mid B\text{ is commutative}\right\}$\\[5pt]
\\
\parbox{9.5cm}{of $\textbf{X}$ and $\textbf{Y}$ respectively (defined in ~\ref{Morita}) we have 
  $\Phi\left(\textbf{X}(1)\right)=\textbf{Y}(1)$ and $\Psi\left(\textbf{Y}(1)\right)=\textbf{X}(1)$ as well as $\Phi\left(\textbf{X}'\right)=\textbf{Y}'$ and $\Psi\left(\textbf{Y}'\right)=\textbf{X}'$. 
  The function $\Phi$ restricted to $\textbf{X}(1)$ maps $[A,\ma]$ to $\left[\End_A(P(1))^{op},P(1)\right]$, where $\left\{1\right\}=\min(\Lambda,\ma)$. \\
  The picture to the right  visualises this situation. 
}
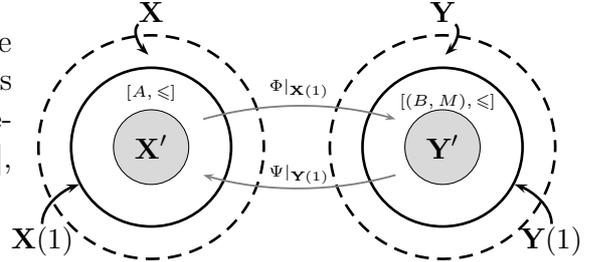
\psset{xunit=1.3cm,yunit=1.2cm,runit=1.2cm}
\begin{pspicture}(-1.3,1.3)(0,0)
\pscircle[linewidth=1pt,linestyle=dashed](0,0){1.5cm}
\pscircle[linewidth=1pt](0,0){1.07cm}
\pscircle*[linecolor=gray!30,linewidth=1pt](0,0){0.5cm}
\pscircle[linewidth=0.3pt](0,0){0.5cm}
\rput(-3,0){
\pscircle[linewidth=1pt,linestyle=dashed](6,0){1.5cm}
\pscircle[linewidth=1pt](6,0){1.07cm}
\pscircle*[linecolor=gray!30,linewidth=1pt](6,0){0.5cm}
\pscircle[linewidth=0.3pt](6,0){0.5cm}
\begin{small}
\rput(6,1.5){\rnode{xx}{$\textbf{Y}$}}
\rput(6,0){\rnode{yb}{$\textbf{Y}'$}}
\rput(6,0.99){\rnode{yy}{}}
\rput[b](7.1,-1.2){\rnode{xx1}{$\textbf{Y}(1)$}}
\rput(6.7,-0.4){\rnode{yy1}{}}
\end{small}
}
\begin{small}
\rput(0,1.5){\rnode{x}{$\textbf{X}$}}
\rput(0,0){\rnode{xb}{$\textbf{X}'$}}
\rput(0,0.99){\rnode{y}{}}
\rput[b](-1.1,-1.2){\rnode{x1}{$\textbf{X}(1)$}}
\rput(-0.7,-0.4){\rnode{y1}{}}
\end{small}
\begin{tiny} 
\rput(0,0.6){\rnode{A}{$[A,\ma]$}}
\rput(3,0.5){\rnode{A}{$\left[(B,M),\ma \right]$}}
\rput(1.5,0.65){\rnode{Psi}{$\Phi|_{\textbf{X}(1)}$}}
\rput(1.5,-0.3){\rnode{Phi}{$\Psi|_{\textbf{Y}(1)}$}}
\end{tiny}
\rput(0.5,0.3){\rnode{l}{}}
\rput(2.5,0.3){\rnode{r}{}}
\rput(0.5,-0.3){\rnode{ll}{}}
\rput(2.5,-0.3){\rnode{rr}{}}
\psset{nodesep=1pt}
\ncarc[linewidth=1pt,arcangle=-45]{->}{x}{y}
\ncarc[linewidth=1pt,arcangle=45]{->}{xx}{yy}\ncarc[linewidth=0.7pt,arcangle=16,linecolor=gray]{->}{l}{r}
\ncarc[linewidth=0.7pt,arcangle=16,linecolor=gray]{->}{rr}{ll}
\ncline[linewidth=0.7pt]{->}{RT1}{RT2}
\ncline[linewidth=0.7pt]{<-}{b2}{b1}
\ncline[linewidth=0.7pt]{->}{a1}{a2}
\ncarc[linewidth=1pt,arcangle=35]{->}{x1}{y1}
\ncarc[linewidth=1pt,arcangle=35]{<-}{yy1}{xx1}

\psset{nodesep=2pt,offset=2pt}
\ncline[linestyle=dashed,linecolor=black]{->}{PsiBM}{RPsiBM}
\ncline[linestyle=dashed,linecolor=black]{->}{RPsiBM}{PsiBM}
\ncline[linestyle=dashed,linecolor=black]{->}{BM}{BMv}
\ncline[linestyle=dashed,linecolor=black]{->}{BMv}{BM}
\end{pspicture}
\\[0.2cm]

\begin{subnum}\begin{normalfont}\textbf{Example.} Let $n\geq 3$ and  $C=\left(c_{ij}\right)_{2\leq i,j\leq n-1}\in \text{GL}_{n-2}(K)$. We define  $B:= B_n(C)= K\left\langle x_2,\ldots , x_{n-1}\right\rangle/\texttt{I}$ \    with \   
$\texttt{I}:= \left\langle \left\{ c_{mk}\cdot x_i\cdot x_j-c_{ij}\cdot x_m\cdot x_k, \ \  x_i^3 \mid  2\leq j,i,k,m\leq n-1\right\}\right\rangle$. 

 Let $X_k=x_k+\texttt{I}$ for any $k\in\Gamma:= \left\{2,\ldots,n-1\right\}$. 
Since $\text{det} \ C\neq 0$, for every $i\in \Gamma$   there exist  $l(i),r(i)\in \Gamma$ such that $c_{l(i)i}\neq 0$ and $c_{ir(i)}\neq 0$. 
Therefore   $c_{l(i)i}  X_{j} X_{r(j)} =c_{jr(j)} X_{l(i)}  X_i$  implies $\left\langle X_{j} X_{r(j)} \right\rangle = \left\langle X_{l(i)}  X_i \right\rangle$ for all $i,j\in \Gamma$. Moreover, for  any two $l_1(i), l_2(i)\in \Gamma$ with $c_{l_1(i)i}\neq 0$ and $c_{l_2(i)i}\neq 0$ we have $c_{l_1(i)i}  X_{l_2(i)} X_{i} =c_{l_2(i)i} X_{l_1(i)}  X_i$, hence $\left\langle  X_{l_2(i)} X_{i}\right\rangle  = \left\langle  X_{l_1(i)}  X_i\right\rangle$. A similar situation holds for any two $r_1(i), r_2(i)\in \Gamma$ with $c_{ir_1(i)}\neq 0$ and $c_{ir_2(i)}\neq 0$.
Since $X_m^3=0$ for all $m\in \Gamma$, we obtain  $X_i X_j X_k=0$ for all $i,j,k\in \Gamma$.   

Furthermore,  we have $c_{ij}=0$ iff $X_iX_j=0$. 
Thus,    $\left\langle X_i\right\rangle = B\cdot X_i=\spann_K\left\{X_i,X_{l(i)} X_i\right\} =\spann_K\left\{X_i,X_i X_{r(i)}\right\}=X_i\cdot B$  is a two-sided local ideal  of $B$ and  $\soc B= \left\langle X_{j} X_{r(j)}\right\rangle $ for all $j\in \Gamma$. 
The algebra $B$ is self-injective and $\dimm_KB=n$. Let    $\left(\Lambda=\left\{1,\ldots, n\right\},\ma\right)$  be  the  poset given by  $1\triangleleft i\triangleleft n$ for all $2\leq i\leq n-1$ and  let the  $B$-module   $L:= \bigoplus_{i=1}^{n} L(i)$ be  given by  $L(1)=B$,  $L(i)=\left\langle X_i\right\rangle $  for all $i\in \Gamma$    and $L(n)=\soc B$. The pair $(B,L)$
\\[3pt]
\parbox{12.7cm}{
  satisfies the property \fbox{$\ma$}. 
The quiver and relations of the 1-quasi-hereditary algebra $A_{n}(C):=\End_B(L)^{op}$ can by found  in \cite[Example 3]{P1}. 
 The algebra $B$ is commutative if and only if  $C=C^{t}$ and only  in this case $A_{n}(C)$ is a $\BGG_{(\leftrightarrows)}$-algebra. }
\psset{xunit=1.1mm,yunit=1mm,runit=1mm}
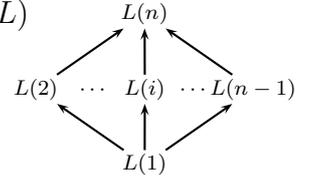
\begin{pspicture}(-13,5)(0,0)
\rput(5,15){\begin{scriptsize}
              \rput(0,0){\rnode{0}{$L(n)$}}
              \rput(-13,-10){\rnode{1}{$L(2)$}}
              \rput(-6,-10){\rnode{a}{$\cdots$}}
              \rput(0,-10){\rnode{3}{$L(i)$}}
              \rput(6,-10){\rnode{a}{$\cdots$}}
              \rput(13,-10){\rnode{6}{$L(n-1)$}}
              \rput(0,-20){\rnode{7}{$L(1)$}} 
\end{scriptsize}
              \psset{nodesep=1.5pt,arrows=<-}
              \ncline{0}{1}
              \ncline{0}{2}
              \ncline{0}{3}
              \ncline{0}{6}
              \ncline{1}{7}
              \ncline{2}{7}
              \ncline{3}{7}
              \ncline{6}{7}}
\end{pspicture}
\end{normalfont}
\label{Example_C}
\end{subnum}

\section{Proof of  the  Theorem  A}
In this section let $(\Lambda,\ma)$ be a poset. 
Until the end of this paper  for any $j\in \Lambda$  we denote by $\Lambda_{(j)}$ and $\Lambda^{(j)}$ the following subsets of $\Lambda$:
\begin{center}
$\Lambda_{(j)}:=\left\{i\in \Lambda \mid i\ma j\right\}$ \  and \  $\Lambda^{(j)}:=\left\{i\in \Lambda \mid i\geqslant j\right\}$.
\end{center}
We also  adopt  all notation of the previous section.

\subsection{Proof of  the  Theorem  A    $(i)\Rightarrow (ii)$}
In this subsection  $A$ with $(\Lambda,\ma)$ denotes   a 1-quasi-hereditary algebra and  $\left\{1\right\}=\min(\Lambda,\ma)$.  Furthermore,   $B:=\End_A(P(1))^{op}$ and  $L(i):= \Hom_A(P(i),P(1))$ for any $i\in \Lambda$. Since $P(1)$ is  a  minimal faithful $A$-module,     Theorem~\ref{Morita}  and   Lemma ~\ref{mini} provide  isomorphisms 
\begin{center}
 $\displaystyle A\cong \End_B\left({}_{B}P(1)\right)$  \  and \     $\displaystyle {}_{B}P(1) \cong \bigoplus_{i\in \Lambda}L(i)$.
\end{center}
In particular, for any  $[\B,\LL]\in \textbf{Y}$ with $A\cong \End_{\B}\left(\LL\right)^{op}$ we obtain   $\B\cong B$ and $\LL\cong {}_{B}P(1)$. For the implication   $(i)\Rightarrow (ii)$  we  have to show that $\left(B,\bigoplus_{i\in \Lambda}L(i)\right)$ satisfies the property~ \fbox{$\ma$}.
We recall some  notations and  properties of 1-quasi-hereditary algebras from \cite{P} and \cite{P1} which  we will use in the  proof:  
Let $Q$ be the quiver of $A$ and $\II$ be the corresponding ideal of $KQ$ generated by the  relations of $A$. The structure of $Q$ (see~\ref{BGG-algebra}) shows that  for all $j,i,k \in  \Lambda$ with $i\in \Lambda^{(j)}\cap \Lambda^{(k)}$ 
  there exists a path $ (j\to \lambda_1\to \cdots \to \lambda_m\to i) $  with $ j\ma \lambda_1\ma \cdots \ma \lambda_m\ma i$  and a path 
  $ (i\to \mu_1\to \cdots \to \mu_r\to k) $  with $i\geqslant \mu_1\geqslant  \cdots \geqslant \mu_m\geqslant  k$. We write $p(j,i,i)$  resp.  $p(i,i,k)$ for the residue class  $A$ of such a path.   The concatenation   of  these two paths  is denoted   by   $p(j,i,k)=p(i,i,k)\cdot p(j,i,i)$. 
  
In the next subsection (Lemma ~\ref{polynom}) we will show, that for  any  two paths $p$ and $q$ in $Q$ which  yield  paths of the form $p(j,i,k)$ in $A$  we have $p-q\in \II$, thus  we can talk  about   \textit{the}  path $p(j,i,k)$ in $A$.    The $A$-map  from $P(k)$ to $P(j)$   corresponding to $p(j,i,k)$  (i.e., $e_k\mapsto p(j,i,k)$)  is denoted by $\ger{f}_{(j,i,k)} $.  For  all   $j,i,k$ with $i\in \Lambda^{(j)}\cap \Lambda^{(k)}$  we have
\begin{center}
$\ger{f}_{(j,i,k)}=\ger{f}_{(j,i,i)}\circ \ger{f}_{(i,i,k)}:\left(P(k)  \stackrel{\ger{f}_{(i,i,k)}}{\longrightarrow} P(i)  \stackrel{\ger{f}_{(j,i,i)}}{\longrightarrow}  P(j)\right)$; \ \  $e_k\mapsto p(j,i,k)=p(i,i,k)\cdot p(j,i,i)$.
\end{center}
For any  $i\in \Lambda^{(j)}$  the map  $\ger{f}_{(j,i,i)}$ is an  inclusion and  $\im (\iota)=\im \left(\ger{f}_{(j,i,i)}\right)$  for any inclusion $P(i)\stackrel{\iota}{\hookrightarrow} P(j)$  (see \cite[3.1(a)]{P}).  The  contravariant    functor  $\Hom_A(-, P(1)): \modd A \rightarrow \modd B$ is exact  since    $P(1)\cong I(1)$ (see \cite[2.1(3)]{P}). The inclusion  $\ger{f}_{(j,i,i)}: P(i)\hookrightarrow P(j)$  induces the surjective $B$-map  $ \underbrace{\Hom_A(P(j), P(1))}_{L(j)} \sur \underbrace{\Hom_A(P(i), P(1))}_{L(i)}$ with $g\mapsto g\circ \ger{f}_{(j,i,i)}$. Since $P(i)\hookrightarrow P(j)$ if and only if $i\in \Lambda^{(j)}$ (see \cite[2.2]{P}), we obtain $L(j)\sur L(i)$ if and only if $j\ma i$, thus  part (2)(a) of the Definition \fbox{$\ma$}  for $(B,\bigoplus_{i\in \Lambda}L(i))$  is satisfied. 
In particular,  the surjection 
  $ {}_{B}B\cong L(1) \sur L(i)$   with $F\mapsto F\circ \ger{f}_{(1,i,i)}$  provides  
$L(i)= B\circ \ger{f}_{(1,i,i)}$ for any $i\in \Lambda$.

It is enough to  show the following two  statements:
\begin{itemize}
	\item[\ding{182}] The algebra $B$  is local, self-injective and  $\dimm_KB=\left|\Lambda\right|$, 
	\item[\ding{183}]    $L(i) \cong B \circ \ger{f}_{(1,i,1)}=\ger{f}_{(1,i,1)}\circ B$  and  $\displaystyle \rad \left(B \circ \ger{f}_{(1,i,1)}\right)=\sum_{i<i'}\left(B \circ \ger{f}_{(1,i',1)}\right)$  for any $i\in \Lambda$ .
\end{itemize}
The second  statement implies that  $L(i)$ is isomorphic to a two-sided ideal of $B$ for any $i\in \Lambda$ and   $L(1)= {}_{B}B$  because  $\ger{f}_{(1,1,1)}=\id_{P(1)}$. Moreover, we obtain the following  explicit  expression of the  $B$-module   $P(1)$:

\begin{subnum}\begin{normalfont}\textbf{Proposition.} \textit{For a 1-quasi-hereditary algebra $(A,\ma)$ with $\left\{1\right\}=\min \left(Q_0(A),\ma\right)$  and $B=\End_A(P(1))^{op}$  we have  \  
${}_{B}P(1)\cong \bigoplus_{i\in \Lambda}B\circ \ger{f}_{(1,i,1)}$,   \  
where  $\ger{f}_{(1,i,1)}$ is the endomorphism of $P(1)$ as described above.}
\end{normalfont} 
\label{End(P(1))}
\end{subnum}

\begin{subnum} \begin{normalfont}\textbf{Remark.} The set   
$\left\{\ger{f}_{(j,i,k)} \mid i\in \Lambda^{(j)}\cap \Lambda^{(k)}\right\}$ is a $K$-basis of $\Hom_A(P(k),P(j))$ for all  $j,k\in \Lambda$ because   the set $\left\{p(j,i,k)\mid i\in \Lambda^{(j)}\cap \Lambda^{(k)}\right\}$ is a $K$-basis of $P(j)_k$ (see \cite[3.2]{P}).
  In particular, since  $\Lambda^{(1)}=\Lambda$ we obtain that $\left\{\ger{f}_{(1,i,1)} \mid i\in \Lambda\right\}$     is a $K$-basis of    $B$    and \begin{center}
    $\left\{\ger{f}_{(1,t,1)} \mid   t\in \Lambda^{(i)}\right\}$   is a $K$-basis of $L(i)=\Hom_A(P(i),P(1))$ for any $i\in \Lambda$.
\end{center}
 This situation gives rise to some implications: Let $f\in \End_A(P(j))^{op}$ and $i\in \Lambda^{(j)}$,  then     $f\circ \ger{f}_{(j,i,i)}\in  \Hom_A(P(i),P(j)) = \spann_K\left\{\ger{f}_{(j,t,i)}\mid t\in \Lambda^{(i)}\right\}$,  since $\Lambda^{(j)}\cap \Lambda^{(i)}=\Lambda^{(i)}$. 
\\
\parbox{13cm}{Thus   $f\circ \ger{f}_{(j,i,i)} = \sum_{t\in \Lambda^{(i)}} c_t\cdot \ger{f}_{(j,t,i)}$ for  $c_t\in K$. Let $k\in \Lambda_{(i)}$, then   
\begin{center}
$p(j,t,k)=p(i,i,k)\cdot p(j,t,i)$   \   implies   \    $\ger{f}_{(j,t,k)}= \ger{f}_{(j,t,i)}\circ \ger{f}_{(i,i,k)}.$
\end{center} 
 We obtain $f$$\circ $$\ger{f}_{(j,i,k)}$$ =$$ \left(f\circ \ger{f}_{(j,i,i)}\right)$$\circ$$ \ger{f}_{(i,i,k)}$$ =$$  \sum_{t\in \Lambda^{(i)}} c_t$$\cdot$$ \ger{f}_{(j,t,k)}$  for every  $t\in \Lambda^{(i)}$. In other  words, for any  $i\in \Lambda^{(j)}\cap \Lambda^{(k)}$ and  $f\in \End_A(P(j))^{op}$ we have
}
\psset{xunit=0.4mm,yunit=0.43mm,runit=1mm}
\begin{pspicture}(0,0)(0,0)
\rput(32,21){
\begin{tiny}
\rput[l](25,2){\rnode{p1}{$p(j,t,i)$}}
\rput[l](25,-6){\rnode{p2}{$p(j,t,k)$}}
\rput[l](25,-14){\rnode{p3}{$p(j,i,i)$}}
\rput[l](25,-22){\rnode{p4}{$p(i,i,k)$}}
\rput(24,2){\rnode{p1r}{}}
\rput(2.5,-9){\rnode{p1l}{}}
\rput(24,-6){\rnode{p2r}{}}
\rput(9,-14){\rnode{p2l}{}}
\rput(24,-14){\rnode{p3r}{}}
\rput(-6.9,-23){\rnode{p3l}{}}
\rput(24,-22){\rnode{p4r}{}}
\rput(7,-21){\rnode{p4l}{}}
\rput(0,10){\rnode{A}{$n$}}
\rput(0,-44){\rnode{B}{$1$}}
\rput(20,-17){\rnode{aa}{}}
\rput(-20,-17){\rnode{bb}{}}
\rput(-11,3){\rnode{A1}{}}
\rput(-11,-37){\rnode{B1}{}}
\rput(11,3){\rnode{A2}{}}
\rput(11,-37){\rnode{B2}{}}
              \rput(0,1){\rnode{01i}{$t$}}              
              \rput(-7,-33){\rnode{1}{$j$}}              
              \rput(2,-15){\rnode{2}{$i$}}
              \rput(8,-32){\rnode{2k}{$k$}}
              \rput(6,-10){\rnode{1r}{}}
              \rput(-7,-10){\rnode{2l}{}}
              \end{tiny}
                      }
                      
             \psset{nodesep=1pt,offset=1pt}
              \ncline{1}{01}
              \ncline{->}{01}{2}
              \ncline[linestyle=dotted]{21}{22}
              \ncline[linewidth=0.4pt]{->}{A}{A1}
              \ncline[linewidth=0.4pt]{->}{A1}{A}
              \ncline[linewidth=0.4pt]{->}{B1}{B}
              \ncline[linewidth=0.4pt]{->}{B}{B1}
              \ncline[linewidth=0.4pt]{->}{A}{A2}
              \ncline[linewidth=0.4pt]{->}{A2}{A}
              \ncline[linewidth=0.4pt]{->}{B2}{B}
              \ncline[linewidth=0.4pt]{->}{B}{B2}
              \ncarc[linewidth=0.6pt,arcangle=25,linestyle=dotted]{A2}{aa}
              \ncarc[linewidth=0.6pt,arcangle=-25,linestyle=dotted]{->}{aa}{A2}
              \ncarc[linewidth=0.6pt,arcangle=25,linestyle=dotted]{->}{aa}{B2}
              \ncarc[linewidth=0.6pt,arcangle=-25,linestyle=dotted]{B2}{aa}
              \ncarc[linewidth=0.6pt,arcangle=-25,linestyle=dotted]{A1}{bb}
              \ncarc[linewidth=0.6pt,arcangle=-25,linestyle=dotted]{->}{bb}{B1}
              \ncarc[linewidth=0.6pt,arcangle=25,linestyle=dotted]{B1}{bb}
              \ncarc[linewidth=0.6pt,arcangle=25,linestyle=dotted]{->}{bb}{A1}
              \ncarc[linewidth=0.3pt, arcangle=-20]{->}{p1r}{p1l}
              \ncarc[linewidth=0.3pt, arcangle=-20]{->}{p2r}{p2l}
              \ncarc[linewidth=0.3pt, arcangle=-20]{->}{p3r}{p3l}
              \ncarc[linewidth=0.3pt, arcangle=-20]{->}{p4r}{p4l}
              \psset{nodesep=0pt,offset=1pt}
              \ncarc[arcangle=25]{1}{01}
              \ncarc[arcangle=25]{->}{1}{2}
              \ncarc[arcangle=25]{->}{01}{2k}
              \ncarc[arcangle=25]{1}{01i}
              \ncarc[arcangle=15]{->}{01i}{2}
              \ncarc[arcangle=25]{->}{2}{2k}
              \ncarc[arcangle=30]{->}{01i}{2k}
\end{pspicture} 
\\
$f\circ \ger{f}_{(j,i,k)}\in \spann_K\left\{\ger{f}_{(j,t,k)} \mid t\in \Lambda^{(i)}\right\}$. 
 In particular,  $B\circ \ger{f}_{(1,i,k)}\subseteq \spann_K\left\{\ger{f}_{(1,t,k)} \mid t\in \Lambda^{(i)}\right\}$. 
\end{normalfont} 
\label{Proofhin}
\end{subnum}

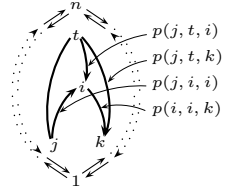
\textit{Proof of} \ding{182}. The algebra $A^{op}$ is also 1-quasi-hereditary with $P_{A^{op}}(1)$ as a minimal faithful $A^{op}$-module (see \cite[1.3]{P}).  The  algebras $B^{op}\cong \End_A(P(1))$ and $B \cong \End_{A^{op}}(P_{A^{op}}(1))\cong \End_A(P(1))^{op}$ are local, since $P(1)$ and $P_{A^{op}}(1)$ are  indecomposable. Therefore  $\soc B$ and  $\topp B$ are simple, thus $B$ is   local, self-injective  and 
$\dimm_KB= \left|\left\{\ger{f}_{(1,i,1)} \mid i\in \Lambda\right\}\right|=\left|\Lambda\right|$. 
\hfill $\Box$ 
\\

 The proof of \ding{183} is based on  the following  properties of $B$-modules generated by the maps~ $\ger{f}_{(1,i,k)}$.

\begin{subnum}\begin{normalfont}\textbf{Lemma.}\end{normalfont} 
\textit{For all  $i,k\in \Lambda$ with $i\in \Lambda^{(k)}$   the following is satisfied:
\begin{itemize}
		\item[(1)] $B\circ \ger{f}_{(1,i,k)}=\spann_K\left\{ \ger{f}_{(1,t,k)}\mid t\in \Lambda^{(i)}\right\}$  for all $k\in \Lambda_{(i)}$,
	\item[(2)] $B \circ \ger{f}_{(1,i,k)} \cong B \circ \ger{f}_{(1,i,k')}$  for all    $k,k'\in \Lambda_{(i)}$,
	\item[(3)] $B\circ \ger{f}_{(1,i,1)} = \ger{f}_{(1,i,1)} \circ B = \spann_{K}\left\{\ger{f}_{(1,t,1)} \mid t\in \Lambda^{(i)}\right\}$.  
\end{itemize}}
\label{hin}
\end{subnum}
\begin{small}
We will clarify    the   statements  of this  lemma   using an example:  Let $A$ be the algebra corresponding to a regular block of $\OO(\ger{sl}_3(\CC))$, then  $P(1)$ is the minimal faithful $A$-module  (the presentation of  $P(1)$ via  the  $\CC$-basis $p(1,i,k)$ can by found  in ~\cite[Sec.2]{P1}).     
\end{small}
\\[2pt]
\parbox{10cm}{\begin{small} The picture  presents  $P(1)$  as a $\End_A(P(1))^{op}$-module (in particular, $\End_A(P(1))^{op}\cong \CC[x,y]/\left\langle xy,\   x^3-y^3\right\rangle$,  see Example ~\ref{examp}): The   circle $C(k)$   represents the  spaces  $L(k)=\Hom_A(P(j),P(1))$ for every $k\in \Lambda=\left\{1,\ldots , 6\right\}$.  
Inside  of  $C(k)$ we find  the  $\CC$-basis $\left\{\ger{f}_{(i,k)}:=\ger{f}_{(1,i,k)} \mid i\in \Lambda^{(k)}\right\}$ of $L(k)$.  There is an arrow $\ger{f}_{(t,k)} \rightarrow \ger{f}_{(i,k)}$    if   $B \circ \ger{f}_{(t,k)} \subset B \circ \ger{f}_{(i,k)}$ (statement \textit{(1)}) and an arrow 
  $\ger{f}_{(i,k)} \dasharrow \ger{f}_{(i,k')}$ if   $B \circ \ger{f}_{(i,k)} \cong B \circ \ger{f}_{(i,k')}$  (statement \textit{(2)}). 
 We focus on the  $B$-submodules of $L(k)$  generated by $\ger{f}_{(4,k)}$   for   $k\in \Lambda_{(4)}=\left\{1,2,3,4\right\}$.  Statement \textit{(1)} yields  $B\circ \ger{f}_{(4,k)}=\spann_K\left\{\ger{f}_{(4,k)}, \ger{f}_{(6,k)}\right\}$  (visualized in the  parallelogram).  Statement \textit{(2)} implies   $L(4) \cong B \circ \ger{f}_{(4,k)} $ for all $k=\Lambda_{(4)}$.  The      $B$-module $L(4)$ is isomorphic to  the two-sided  ideal   $B\circ \ger{f}_{(4,1)} = \ger{f}_{(4,1)} \circ B$ of  $B$ since $\ger{f}_{(4,1)} \in B$ (statement~ \textit{(3)}).
 \end{small}}  
 \psset{xunit=0.52mm,yunit=0.53mm,runit=6mm}
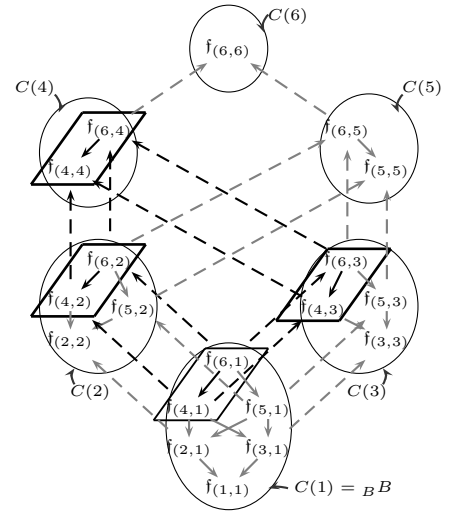
\begin{pspicture}(0,4)(0,0)
\rput(64,-20){
\begin{tiny}
\psellipse[linewidth=0.4pt](0,80)(10,11)
\psellipse[linewidth=0.4pt](35.5,55)(12.5,13.8)
\psellipse[linewidth=0.4pt](-35.5,54)(12.5,13.8)
\psellipse[linewidth=0.4pt](33,15.5)(15,16.9)
\psellipse[linewidth=0.4pt](-33,15.5)(15,16.9)
\psellipse[linewidth=0.4pt](0,-14.5)(16,21)
\rput(0,80){\rnode{66o}{}}
\rput(-5,77){\rnode{66v}{}}
\rput(0,80){\rnode{666}{$\ger{f}_{(6,6)}$}}

\rput(15,88){\rnode{P6}{$C(6)$}}
\rput(10,88.5){\rnode{P6p}{}}
\rput(9,84){\rnode{P6g}{}}
\ncarc[arcangle=-35,linecolor=black!80]{->}{P6p}{P6g}
 
\rput(30,60){
\rput(4,0){
\rput(15,10){\rnode{P5}{$C(5)$}}
\rput(10,8.5){\rnode{P5p}{}}
\rput(9,4){\rnode{P5g}{}}
\ncarc[arcangle=-35,linecolor=black!80]{->}{P5p}{P5g}}
             \rput(0,0){\rnode{65o}{}}
             \rput(10,-10){\rnode{55o}{}}
             \rput(0,0){\rnode{65}{$\ger{f}_{(6,5)}$}}
             
             \rput(10,-10){\rnode{55}{$\ger{f}_{(5,5)}$}}

}
\rput(-30,60){
\psline[linewidth=1pt]{-}(-7,4)(9,4)
\psline[linewidth=1pt]{-}(-20,-14)(-4,-14)
\psline[linewidth=1pt]{-}(-7,4)(-20,-14)
\psline[linewidth=1pt]{-}(9,4)(-4,-14)
\rput(-4,0){
\rput(-15,10){\rnode{P4}{$C(4)$}}
\rput(-10,8.5){\rnode{P4p}{}}
\rput(-9,4){\rnode{P4g}{}}
\ncarc[arcangle=35,linecolor=black!80]{->}{P4p}{P4g}}
              \rput(0,0){\rnode{64o}{}}
              \rput(-10,-10){\rnode{44o}{}}
              \rput(0,0){\rnode{64}{$\ger{f}_{(6,4)}$}}              
              \rput(-10,-10){\rnode{44}{$\ger{f}_{(4,4)}$}} 
}
\rput(-20,17){
\rput(-10,10){\psline[linewidth=1pt]{-}(-7,4)(9,4)
\psline[linewidth=1pt]{-}(-20,-14)(-4,-14)
\psline[linewidth=1pt]{-}(-7,4)(-20,-14)
\psline[linewidth=1pt]{-}(9,4)(-4,-14)}
\rput(-10,-5){
\rput(-5,-18){\rnode{P2}{$C(2)$}}
\rput(-9,-16){\rnode{P2p}{}}
\rput(-9,-12){\rnode{P2g}{}}
\ncarc[arcangle=35,linecolor=black!80]{->}{P2p}{P2g}}
              \rput(-10,12){\rnode{62o}{}}
              \rput(-20,0){\rnode{42o}{}}
              \rput(-4,-2){\rnode{52o}{}}
              \rput(-4,0){\rnode{52u}{}}
              \rput(-20,-10){\rnode{22o}{}}                
              \rput(-10,10){\rnode{62}{$\ger{f}_{(6,2)}$}}             
              \rput(-20,0){\rnode{42}{$\ger{f}_{(4,2)}$}}
              \rput(-4,-2){\rnode{52}{$\ger{f}_{(5,2)}$}}
              \rput(-20,-10){\rnode{22}{$\ger{f}_{(2,2)}$}}  
               }
\rput(20,17){
\rput(12,8.9){\psline[linewidth=1pt]{-}(-7,4)(9,4)
\psline[linewidth=1pt]{-}(-20,-14)(-4,-14)
\psline[linewidth=1pt]{-}(-7,4)(-20,-14)
\psline[linewidth=1pt]{-}(9,4)(-4,-14)}
\rput(10,-5){
\rput(5,-18){\rnode{P2}{$C(3)$}}
\rput(9,-16){\rnode{P2p}{}}
\rput(9,-12){\rnode{P2g}{}}
\ncarc[arcangle=-35,linecolor=black!80]{->}{P2p}{P2g}
}
              \rput(10,10){\rnode{63o}{}}
              \rput(10,12){\rnode{63u}{}}
              \rput(20,0){\rnode{53o}{}}
              \rput(4,-2){\rnode{43o}{}}
              \rput(4,0){\rnode{43u}{}}
              \rput(20,-10){\rnode{33o}{}}  
              
              \rput(10,10){\rnode{63}{$\ger{f}_{(6,3)}$}}              
              \rput(20,0){\rnode{53}{$\ger{f}_{(5,3)}$}}
              \rput(4,-2){\rnode{43}{$\ger{f}_{(4,3)}$}}
              \rput(20,-10){\rnode{33}{$\ger{f}_{(3,3)}$}}   
        }
       
       \rput(0,2){\rnode{6l}{}} 
       \rput(0,2){\rnode{6r}{}} 
        \rput(-10,-10){\rnode{4l}{}} 
        \rput(-9,-13){\rnode{4r}{}} 
        \rput(10,-10){\rnode{5l}{}} 
        \rput(9,-13){\rnode{5r}{}} 
        \rput(-10,-20){\rnode{2l}{}} 
        \rput(10,-20){\rnode{3l}{}} 
        \rput(0,-31){\rnode{1l}{}}  
        
      \rput(10,-15){
\rput(19,-15){\rnode{P1}{$C(1)={}_BB$}}
\rput(5,-15){\rnode{P1p}{}}
\rput(0,-15){\rnode{P1g}{}}
\ncarc[arcangle=-15,linecolor=black!80]{->}{P1p}{P1g}}

\rput(1,1){\psline[linewidth=0.7pt]{-}(-7,4)(9,4)
\psline[linewidth=0.7pt]{-}(-20,-14)(-4,-14)
\psline[linewidth=0.7pt]{-}(-7,4)(-20,-14)
\psline[linewidth=0.7pt]{-}(9,4)(-4,-14)}

       \rput(0,2){\rnode{61}{$\ger{f}_{(6,1)}$}}       
        \rput(-10,-10){\rnode{41}{$\ger{f}_{(4,1)}$}} 
        \rput(10,-10){\rnode{51}{$\ger{f}_{(5,1)}$}}  
       \rput(-10,-20){\rnode{21}{$\ger{f}_{(2,1)}$}} 
        \rput(10,-20){\rnode{31}{$\ger{f}_{(3,1)}$}} 
        \rput(0,-30){\rnode{11}{$\ger{f}_{(1,1)}$}}      
              \end{tiny}}
              \psset{nodesep=0.5pt,offset=0pt,linecolor=black!50,arrows=<-}
              \ncline[linecolor=black]{44}{64}
              \ncline{55}{65}
              
              \ncline{52}{62}
              \ncline[linecolor=black]{42}{62}
              \ncline{22}{52}
              \ncline{22}{42}
              \ncline{33}{53}
              \ncline{33}{43}
              \ncline[linecolor=black]{43}{63}
              \ncline{53}{63}
              \ncline[linecolor=black]{41}{61}
              \ncline{51}{61}
              \ncline{21}{51}
              \ncline{21}{41}
              \ncline{31}{51}
              \ncline{31}{41}
              \ncline{11}{21}
              \ncline{11}{31}
              \psset{nodesep=0.5pt,offset=0pt,linecolor=black!50,arrows=->}
              \ncline{pi1}{pi12}
              
              \psset{nodesep=8pt,linecolor=black!50,linestyle=dashed,arrows=->}
              \ncline[linecolor=black]{6l}{62o}
              \ncline{5r}{52u}
              \ncline[linecolor=black]{4l}{42o}
              \ncline{2l}{22o}
                \ncline[linecolor=black]{6r}{63u}
              \ncline{5l}{53o}
              \ncline[linecolor=black]{4r}{43u}
              
              \ncline{3l}{33o}
              \ncline[linecolor=black]{62o}{64o}
              \ncline[linecolor=black]{42o}{44o}
              \ncline{62o}{65o}
              \ncline{52o}{55o}

              \ncline{63o}{65o}
              \ncline{53o}{55o}
              \ncline[linecolor=black]{63o}{64o}
              \ncline[linecolor=black]{43o}{44o}
              
              \ncline{64o}{66o}
              \ncline{65o}{66o}
              \ncarc[arcangle=-35]{5}{57}
              \ncarc[arcangle=35]{v5}{v7}
              
\end{pspicture}
\\
\\

\textit{Proof.} \textit{(1)} Since $P(1)\cong I(1)$ and  $\ger{f}_{(1,i,i)}:P(i)\hookrightarrow P(1)$,  the universal property  of injective modules imply   that  for
 any $t\in \Lambda^{(i)}$ there exists 
   $F(t)\in B$   such that  $F(t)\circ \ger{f}_{(1,i,i)} = \ger{f}_{(1,t,i)}$. Let  $k\in \Lambda_{(i)}$, then   $\ger{f}_{(1,i,i)}\circ \ger{f}_{(i,i,k)} = \ger{f}_{(1,i,k)}$  provides  the commutative diagram
\begin{center}
\begin{small}
$
\begin{array}{rcclccc}
\ger{f}_{(1,t,k)}:& P(k)& \stackrel{\ger{f}_{(i,i,k)}}{\longrightarrow} & P(i)& \stackrel{\ger{f}_{(1,t,i)}}{\longrightarrow} & P(1) &   \\
& || & & \ \   \downarrow \ger{f}_{(1,i,i)}& & || &  \\
F(t)\circ\ger{f}_{(1,i,k)}:& P(k)& \stackrel{\ger{f}_{(1,i,k)}}{\longrightarrow} & P(1) & \stackrel{F(t)}{\longrightarrow}& P(1) & 
\end{array}
$
 \end{small}
\end{center}
We obtain   $F(t)\circ \ger{f}_{(1,i,k)} = \ger{f}_{(1,t,i)}\circ \ger{f}_{(i,i,k)}=\ger{f}_{(1,t,k)}$. 
In other words, we have   $B\circ \ger{f}_{(1,i,k)}\supseteq \spann_K\left\{ \ger{f}_{(1,t,k)}\mid t\in \Lambda^{(i)}\right\}$.
 The previous remark  yields $B\circ \ger{f}_{(1,i,k)}\subseteq \spann_K\left\{ \ger{f}_{(1,t,k)}\mid t\in \Lambda^{(i)}\right\}$.

 \textit{(2)} Let $k\in \Lambda_{(i)}$.  We consider the map $\left(-\circ \ger{f}_{(k,k,1)}\right) : B \circ \ger{f}_{(1,i,k)} \rightarrow  B \circ  \ger{f}_{(1,i,1)}$, and  we have   $F\circ \ger{f}_{(1,i,k)} \mapsto F \circ \ger{f}_{(1,i,k)} \circ \ger{f}_{(k,k,1)} = F \circ \ger{f}_{(1,i,1)}$ (see ~\ref{Proofhin}). Obviously,  this  map   is   a  surjective  $B$-map. Moreover,    $\dimm_K\left(B \circ \ger{f}_{(1,i,1)}\right) \stackrel{\textit{(1)}}{=} \dimm_K\left(B \circ \ger{f}_{(1,i,k)}\right)$ since  $\left\{\ger{f}_{(1,t,k)} \mid t\in \Lambda^{(i)}\right\}$  is $K$-independent. Thus $B \circ \ger{f}_{(1,i,1)} \cong B \circ \ger{f}_{(1,i,k)}$ for all $k\in \Lambda_{(i)}$.


\textit{(3)} It is
enough to show  $\ger{f}_{(1,i,1)} \circ B = \spann_{K}\left\{\ger{f}_{(1,t,1)} \mid t\in \Lambda^{(i)}\right\}$   [the rest  follows from \textit{(1)}]:
For any $i\in \Lambda$ we have  $\ger{f}_{(1,i,1)}=\ger{f}_{(1,i,i)}\circ \ger{f}_{(i,i,1)}$  and
 $\ger{f}_{(i,i,1)}\circ F\in \Hom_A\left(P(1),P(i)\right)=\spann_K\left\{\ger{f}_{(i,t,1)}\mid t\in \Lambda^{(i)}\right\}$ for all $F\in B$.  With  similar arguments   as in the  previous remark we obtain  $\ger{f}_{(1,t,1)}=\ger{f}_{(1,i,i)}\circ \ger{f}_{(i,t,1)}$. This  provides  $ \ger{f}_{(1,i,1)}\circ  F \in   \spann_K\left\{\ger{f}_{(1,t,1)}\mid t\in \Lambda^{(i)}\right\}$  for every $F\in B$. Therefore, we obtain    $\ger{f}_{(1,i,1)} \circ B \subseteq \spann_{K}\left\{\ger{f}_{(1,t,1)} \mid t\in \Lambda^{(i)}\right\}$.

 For any $t\in \Lambda^{(i)}$
	 the $A$-module generated by $p(1,t,1)$ is a submodule of $A\cdot p(1,i,1)$ (see \cite[2.2 \textit{(a)}]{P1}).  There exists  $p(t)\in P(1)_1$  with $p(1,t,1)=p(t)\cdot p(1,i,1)$. Let $F(t)\in B$  be  given by  $F(t)(e_1)=p(t)$, then  $\ger{f}_{(1,t,1)} = \ger{f}_{(1,i,1)}\circ F(t)$.  Thus  $\ger{f}_{(1,i,1)} \circ B \supseteq \spann_{K}\left\{\ger{f}_{(1,t,1)} \mid t\in \Lambda^{(i)}\right\}$.\hfill $\Box$ 
 \\

\textit{Proof of} \ding{183}. The  parts \textit{(2)} and \textit{(3)}  of the previous lemma  imply  $L(i)=B\circ \ger{f}_{(1,i,i)} \cong B\circ \ger{f}_{(1,i,1)} = \ger{f}_{(1,i,1)} \circ B$. Thus  $L(i)$ is a two-sided local ideal of $B$ for all $i\in \Lambda$.

The $B$-module $B\circ \ger{f}_{(1,i,1)}$ is local, thus  $B\circ \ger{f}_{(1,t,1)}\stackrel{\textit{(1)}}{\subset} B\circ \ger{f}_{(1,i,1)} $ for all  $t\in \Lambda^{(i)}\backslash \left\{i\right\}=\left\{t\in \Lambda \mid i<t\right\}$ implies $\sum_{i<t}\left(B\circ \ger{f}_{(1,t,1)}\right)\subseteq \rad \left(B\circ \ger{f}_{(1,i,1)}\right) $.
 Since  $\left\{\ger{f}_{(1,t,1)}\mid i<t\right\} $ is linearly independent and   $\left\{\ger{f}_{(1,t,1)}\mid i<t\right\}\subseteq \sum_{i<t}\left(B\circ \ger{f}_{(1,t,1)}\right)$,  we have  $\dimm_K\left(\sum_{i<t}\left(B\circ \ger{f}_{(1,t,1)}\right)\right) \geq \left|\Lambda^{(i)}\right|-1 \stackrel{\textit{(1)}}{=} \dimm_K \left(B\circ \ger{f}_{(1,i,1)}\right)-1=\dimm_K \rad \left(B\circ \ger{f}_{(1,i,1)}\right)$. We obtain  $\sum_{i<t}\left(B\circ \ger{f}_{(1,t,1)}\right)= \rad \left(B\circ \ger{f}_{(1,i,1)}\right) $ for all $i\in \Lambda$.  \hfill $\Box$

\subsection{Proof of the  Theorem  A   $(ii)\Rightarrow (i)$} 
Let $B$ be a local,  self-injective $K$-algebra  with  $\dimm_KB=n$. Let  the set   $(\Lambda=\left\{1,\ldots , n\right\},\ma )$  be     partially ordered   and let   $L\cong \bigoplus_{i\in \Lambda}L(i)$ 
 be  a $B$-module, such that $(B,L)$ satisfies the condition \fbox{$\ma$} (see Definition~\ref{property_ma}).
 
 It is easy to see that there exists a uniquely determined minimal and maximal element in $(\Lambda,\ma)$:
  Since
 $B(=L(1))$ is  a projective cover of  any local $B$-module, we have     $L(1)\sur L(i)$, thus \fbox{$\ma$}(2)(a) implies  $1\ma i$  for all  $i\in \Lambda$.     Let $k\in \Lambda$ be maximal, then      \fbox{$\ma$}(2) (b)   yields     $\rad L(k)=0$. The local submodule   $L(k)$ of $B$ is simple, thus   $L(k)=\soc (B)$, because the socle  of a local self-injective algebra is simple. 
   Since   $L(j)\cong L(i)$ iff  $j= i$ (see \fbox{$\ma$}(2)(a)), we have    $L(k)\subseteq \rad L(i)\stackrel{\text{(b)}}{=}\sum_{i<j}L(j)$   for every  $i\in \Lambda$ with $i\neq k$. Consequently,  $i\ma k$ for all  $i\in \Lambda$.   The  maximal element in $(\Lambda,\ma)$ will be denoted  by   $n=\left|\Lambda\right|$.  We have $\left\{1\right\}=\min(\Lambda,\ma)$ and $\left\{n\right\}=\max(\Lambda,\ma)$ with  $L(1)=B$ and  $L(n)=\soc (B)$.
 
  We consider the  algebra   $A:=\End_B(L)^{op}$. Because   $L$    is multiplicity-free and has $n$ direct summands,  $A$       is  basic    and  the quiver   $Q(A)$ has $n$ vertices.  We identify   these  with the elements in  $\Lambda$. 
In order to  prove   that $(A,\ma)$ is a 1-quasi-hereditary algebra (see Definition ~\ref{def1qh})    we have to show that for all $j\in \Lambda$ the following holds:
 \begin{itemize}
	\item[\ding{182}]  $[\De(j):S(i)]=1$   for all  $i\in \Lambda_{(j)}$,
	\item[\ding{183}]  $P(j)$ has a $\De$-good filtration with $\left(P(j):\De(i)\right)=\left\{
\begin{array}{cl}
1 & \text{if } \  i\in \Lambda^{(j)} ,\\
0& \text{else},
\end{array}
\right.$  
	\item[\ding{184}] $\soc P(j)\cong \topp I(j) \cong  S(1)$, 
	\item[\ding{185}] $\De(j)\hookrightarrow \De(n)$ and $\nabla(n)\sur \nabla(j)$.
\end{itemize}

Recall that  an (left) ideal $J$ of a local, self-injective (basic) algebra $\B$ is local if and only if it is generated by some non-zero  element $x\in \B$ (i.e.,  $J=\B\cdot x$). Moreover, $\dimm_KJ=\dimm_K(\rad J)+1$ and $J= \spann_{K}\left\{x \cup \rad J\right\}$. In addition, for an ideal $M$ of $\B$,  any $\B$-map $f:J\to M$ is induced by right-multiplication by an element $b_f\in \B$, more precisely, $f(y\cdot x)=y\cdot x\cdot b_f$ for all $y\in \B$  (in this case we write $J\stackrel{\cdot b_f}{\rightarrow}M$). The element $f(x)=x\cdot b_f$ generates $\im (f)=\B\cdot f(x)$.   In particular, we have $\Hom_{\B}(J,\B)=\left\{J\stackrel{\cdot b}{\rightarrow}\B\mid b\in \B\right\}$. 

 The annihilator of $L(j)$ is  $\Ann (L(j)):=\left\{b\in B\mid b\cdot L(j)=0\right\}$.

 \begin{subnum}\begin{normalfont}\textbf{Lemma.}\end{normalfont} 
 Let $B$ be a local, self-injective algebra  and let  $(\Lambda, \ma)$  be a poset with  $\dimm_KB=\left|\Lambda\right|=n$.  Let   $\displaystyle  L\cong \bigoplus_{i\in \Lambda}L(i)$  be a $B$-module   such that $(B,L)$ satisfies the condition \fbox{$\ma$}.      Then   for all $i,j,k\in \Lambda$  the following properties are satisfied:
\begin{itemize}
	\item[(1)] Let  $x_i\in B$ be a generator of $L(i)$    and  $X(i):=\left\{x_k\mid k\in \Lambda^{(i)}\right\}$. Then  we have:
\begin{itemize}
	\item[(1.1)] The set  $X(i)$ is a $K$-basis of   $L(i)$. In particular, for  any subset   $\Gamma \subseteq \Lambda$ the set  $\displaystyle 
	\bigcup_{j\in \Gamma}X(j)$ is a $K$-basis of  \  $\displaystyle \sum_{j\in \Gamma}L(j)$ \   and \   $\displaystyle \bigcap_{j\in \Gamma}X(j)$ \   is a $K$-basis of  \  $\displaystyle \bigcap_{j\in \Gamma}L(j)$.
	\item[(1.2)]  $L(i)=B\cdot x_i=x_i\cdot B$.
\end{itemize}
\item[(2)] Let $N$ be a submodule of $B$, then $\im(f)\subseteq L(i)\cap  N$ for all  $f\in \Hom_B(L(i),N)$.
 	\item[(3)] We have  $L(i)\hookrightarrow L(j)$ resp. $L(j)\sur L(i)$  if and only if $i\in \Lambda^{(j)}$.   Moreover, 
\begin{itemize}
	\item[(3.1)] $\im\left(L(i)\hookrightarrow L(j)\right) =L(i)$   for any injective $B$-map  from $L(i)$ to $L(j)$,
	\item[(3.2)] $\ker \left(L(j)\sur L(i)\right)= \Ann (L(i))\cdot L(j)$  for any surjective $B$-map from $L(j)$ to $L(i)$.
\end{itemize}
	\item[(4)]  For  $i\in \Lambda^{(j)}\cap \Lambda^{(k)}$ let  $f_{(k,i,i)}:L(i)\hookrightarrow L(k)$  be an injective,   $f_{(i,i,j)}:L(j)\sur L(i)$   a  surjective $B$-map  and   $f_{(k,i,j)}:= f_{(k,i,i)}\circ f_{(i,i,j)}(\in  \Hom_{B}(L(j),L(k)))$. The set 
\begin{center}
$\ger{B}(k,j):= \left\{f_{(k,i,j)} \mid i\in \Lambda^{(j)}\cap \Lambda^{(k)}\right\}$ \  is a $K$-basis of \   $\Hom_B(L(j), L(k))$.
\end{center}
In particular,
every map $f\in \Hom_B(L(j), L(k))$  factors through $\displaystyle \bigoplus_{i\in \Lambda^{(j)}\cap \Lambda^{(k)}} L(i)$. 
\end{itemize}
\label{L5}
\end{subnum}

\textit{Proof.} 
 \textit{(1.1)} 
 By  induction on $\dimm_K (L(j))$ we show  that  the   $K$-space  $L(j)$ is spanned by  $X(j)$  for every  $j\in \Lambda$:   If $\dimm_K(L(j))=1$,  then $L(j)$ is simple, thus $j=n$  and    $X(n)=\left\{x_n\right\}$ is a $K$-basis of $L(n)=\soc B$.  Let    $j\in \Lambda$ with $\dimm_K(L(j))=m+1$, then    $L(l)=\spann_KX(l)$   for any $l\in \Lambda^{(j)}\backslash\left\{j\right\}$,    because $L(l)\subseteq \rad L(j)$ (see ~\ref{property_ma}(2)(b)) implies    $\dimm_KL(l)\leq \dimm_K\left(\rad L(j)\right)=m$.   Thus     $\rad L(j)=\sum_{j<l}L(l)$ is spanned by  $\bigcup_{j<l}X(l)=\left\{x_k\mid j <k \right\}$  and  consequently  $X(j)=\left\{x_j\right\}\cup \bigcup_{j<l}X(l)$ spans  the $K$-space   $L(j)$.  Furthermore, 
 $B =L(1)= \spann_K X(1)=\left\{x_1,\ldots , x_n\right\}$  and  $\dimm_KB=n$ implies  that $X(1)$ is a $K$-basis of $B$ and consequently  the  subset $X(j)$  of $X(1)$  is linearly independent.
 
 For  $\Gamma\subseteq \Lambda$ the subsets $	\bigcup_{j\in \Gamma}X(j)	=\left\{x_k\mid k\in \Lambda, \  k\geqslant i \text{ for some } i\in \Gamma\right\}$ and  $	\bigcap_{j\in \Gamma}X(j)	=\left\{x_k\mid k\in \bigcap_{j\in \Gamma}\Lambda^{(j)} \right\}$ of $X(1)$ generate $ \sum_{j\in \Gamma}L(j)$  and $\bigcap_{j\in \Gamma}L(j)$ as  $K$-spaces respectively.

\textit{(1.2)}  Since   $L(j)\sur L(i)$ for any  $i\in \Lambda^{(j)}$ (see ~\ref{property_ma}(2)(a)),   there exists    $b_i\in B$  with $L(j)\stackrel{\cdot b_i}{\sur} L(i)$ and   $x_i= x_j \cdot b_i$.  Let     $y\in L(j)=B \cdot x_j$,  then     $ y\stackrel{\textit{(1.1)}}{=}\sum_{i\in\Lambda^{(j)}}c_i\cdot x_i = \sum_{i\in\Lambda^{(j)}}c_i\cdot x_j\cdot b_i = x_j\cdot \left(\sum_{i\in\Lambda^{(j)}}c_i\cdot b_i\right)$  (here  $c_i\in K$). We obtain  $B \cdot x_j\subseteq x_j \cdot B$. 
Since    $L(j)$ is a two-sided ideal, we have     $B\cdot x_j \cdot B\subseteq B\cdot x_j$, thus     $x_j \cdot B\subseteq B \cdot x_j$   and   consequently $L(j)= B\cdot x_j =  x_j \cdot B$.

\textit{(2)}  For  $f\in \Hom_B(L(i),N)$    there exists some $b\in B$ with  $f:L(i)\stackrel{\cdot b}{\rightarrow}N$. Let $x_i$ be a generator of $L(i)$, then  $f(x_i)=x_i\cdot b=\widetilde{b} \cdot x_i$ for some $\widetilde{b}\in B$  by \textit{(1.2)}. Thus $f(x_i)\in L(i)$ and consequently $\im (f)\subseteq L(i)\cap N$.

\textit{(3.1)} According to   ~\ref{property_ma}(2)(b) we have     $L(i)\subseteq L(k)$ if and only if  $i\in \Lambda^{(k)}$. Let  $N$ be  a submodule of $L(k)$ with $N\cong L(i)$ and  $f\in \Hom_B(L(i),N)$  be an  isomorphism, then    $f(L(i))= N \subseteq  L(i) \cap  N$ (see \textit{(2)}). We obtain     $f(L(i))=L(i)=N$.

\textit{(3.2)}  Let    $\pi_i:L(j)\sur L(i)$ be  a surjection  and $x_j$  a generator of $L(j)$.  Then   $x_i:=\pi_i(x_j)$  generates $L(i)$. 
 Let  $x\in L(j)$, then      $x=b\cdot x_j$  for some $b\in B$.  Obviously,    $x\in \ker \left(\pi_i\right)$  if and only if  $\pi_i(x)= b\cdot x_i =0$  if and only if    $ b\cdot x_i\cdot B \stackrel{\textit{(1.2)}}{=} b\cdot L(i)=0$.  We obtain     $\ker (\pi_i)=\left\{b\cdot x_j\in L(j)\mid b\in \Ann (L(i))\right\}= \Ann (L(i))\cdot L(j)$.

 \textit{(4)} Let   $x_j$  be some  generator of   $L(j)$.  Then   $x_i:= f_{(k,i,j)}(x_j)$   generates   the submodule $L(i)$ of $L(j)\cap L(k)$. The set  $\left\{x_i\mid i\in \Lambda^{(j)}\cap \Lambda^{(k)}\right\}\stackrel{\textit{(1.1)}}{=} X(j)\cap X(k)$   is a $K$-basis of $L(j)\cap L(k)$. 
 Let  $f\in \Hom_B(L(j),L(k))$, then    $\im (f)\subseteq L(j)\cap L(k)$ (see  \textit{(2)}). Thus $f(x_j)=\sum_{i\in \Lambda^{(j)}\cap \Lambda^{(k)}}c_i\cdot x_i$ and consequently   $f=\sum_{i\in \Lambda^{(j)}\cap \Lambda^{(k)}}c_i\cdot f_{(k,i,j)}$.
 Let   $f_1, f_2$  be the  $B$-maps  given by   $f_1:L(j)\rightarrow \bigoplus_{i\in \Lambda^{(j)}\cap \Lambda^{(k)}}L(i)$ with $f_1(x_j)=\left( \pi_i(x_j)\right)_{i=1,\ldots,n}$ and $f_2:\bigoplus_{i\in \Lambda^{(j)}\cap \Lambda^{(k)}}L(i) \rightarrow L(k)$ with $f_2\left(x_i\right)_{i=1,\ldots,n}=\sum_{i\in \Lambda^{(j)}\cap \Lambda^{(k)}} c_i\cdot x_i$.   We have  $f=f_2\circ f_1$ and  thus $f$ factors through $\bigoplus_{i\in \Lambda^{(j)}\cap \Lambda^{(k)}}L(i)$.
 \hfill $\Box$

\begin{subnum}\begin{normalfont}\textbf{Remark.}
For a pair $(B,L)$  satisfying   the property \fbox{$\ma$} with $L=\bigoplus_{i\in \Lambda}L(i)$   denote by  $\texttt{B}$ the algebra $B^{op}$ and by $\texttt{L}(i)$ the  $\texttt{B}$-module $\DD(L(i))$, where $\DD:\modd B \to \modd \texttt{B}$  is the standard duality functor.   Since $\topp L(i)$ and  $\soc L(i)$  are   simple,  we obtain that  $\soc \texttt{L}(i)$  and  $\topp \texttt{L}(i)$ are   simple  for all $i\in \Lambda$.  In particular, $\texttt{B}(=\texttt{L}(1))$  is a local,  self-injective algebra with  $\dimm_K\texttt{B}=\dimm_KB=n$  and  $\texttt{L}(i)$  can by considered   as a  local (left) ideal of  $\texttt{B}$ for any $i\in \Lambda$.  

We denote by   $\texttt{f}_{(j,i,k)}$  the $\texttt{B}$-map  $\DD\left(f_{(k,i,j)}\right):\texttt{L}(k)\to \texttt{L}(j)$  for all $i,j,k\in \Lambda$ with $i\in \Lambda^{(j)}\cap \Lambda^{(k)}$,  
 where $f_{(k,i,j)}$ is the $B$-map described in ~\ref{L5}\textit{(4)}. 
As $\DD$ is duality, we have $\texttt{f}_{(j,i,k)}:  \texttt{L}(k) \stackrel{\texttt{f}_{(i,i,k)}}{\sur} \texttt{L}(i) \stackrel{\texttt{f}_{(j,i,i)}}{\hookrightarrow} \texttt{L}(j)$. 
\end{normalfont}
\label{duality}
\end{subnum}

\begin{subnum}\begin{normalfont}\textbf{Proposition.}\end{normalfont} 
A pair $(B,L)$ satisfies the condition \fbox{$\ma$} if and only if $\left(B^{op},\DD(L)\right)$ satisfies the condition \fbox{$\ma$}.
\label{maop}
\end{subnum}

\textit{Proof.} (We use the  notations 	introduced in ~\ref{duality}.)
According to ~\ref{L5} \textit{(3)}   we have    $\texttt{L}(j)\sur \texttt{L}(i)$    resp.  $\texttt{L}(i)\hookrightarrow \texttt{L}(j)$   if and only if  $i\in \Lambda^{(j)}$.  Moreover, 
any two   surjections  $\widetilde{\pi}_1$,  $\widetilde{\pi}_2$ from  $\texttt{L}(j)$ to $\texttt{L}(i)$  are  induced   by  some  injective maps $\iota_1, \iota_2:L(i)\hookrightarrow L(j)$  and $\ker (\widetilde{\pi}_{k}) = \left\{\xi\in \Hom_K(L(j),K) \mid  \xi|_{\im ( \iota_{k})}=0\right\}$, here $k=1,2$. Since $\im(\iota_1)\stackrel{\ref{L5}\textit{(3)}}{=} \im(\iota_2)$, we obtain   $\ker(\widetilde{\pi}_1)=\ker(\widetilde{\pi}_2)$.  Similarly,  $\im \left(\widetilde{\iota}\right) = \texttt{L}(i)$ for all injections $\widetilde{\iota}:\texttt{L}(i)\hookrightarrow \texttt{L}(j)$. 


 Let $\texttt{1}:=\texttt{1}_{\texttt{B}}$  and $\texttt{y}_i:=\texttt{f}_{(j,i,1)}(\texttt{1}) $ for any $i\in \Lambda^{(j)}$. Obviously,  $\texttt{y}_i$ is a generator of the submodule $\texttt{L}(i)$ of  $\texttt{L}(j)$.  
   According to ~\ref{L5}\textit{(4)} the set  $\left\{\texttt{f}_{(j,i,1)} \mid i\in  \Lambda^{(j)}\right\}$ is a $K$-basis of $\Hom_{\texttt{B}}(\texttt{L}(1),\texttt{L}(j))$.   Thus  
 the set  $\texttt{X}(j):=\left\{\texttt{y}_i\mid i\in \Lambda^{(j)}\right\}$ is a $K$-basis of $\texttt{L}(j)$, since  $\dimm_K\texttt{L}(j)=\dimm_KL(j)\stackrel{~\ref{L5}\textit{(1)}}{=}\left|\texttt{X}(j)\right|$.

 Now we  show  $\texttt{L}(j)=\texttt{B}\cdot \texttt{y}_j=\texttt{y}_j\cdot \texttt{B}$ (this implies, that $\texttt{L}(j)$ is a two-sided  ideal of $\texttt{B}$): 
 Let $\texttt{f}_{(i,i,j)}:\texttt{L}(j)\stackrel{\cdot b_i}{\sur}\texttt{L}(i)$  such that  $\texttt{y}_i=\texttt{y}_j\cdot b_i$  for any $i\in \Lambda^{(j)}$.    
 Let   $\texttt{y}\in \texttt{B}\cdot \texttt{y}_j$,  then    $\texttt{y}=\sum_{i\in \Lambda^{(j)}}c_i\cdot \texttt{y}_i=\texttt{y}_j\cdot \left(\sum_{i\in \Lambda^{(j)}}c_i\cdot b_i\right)\in \texttt{y}_j\cdot\texttt{B}$,  thus   $\texttt{B}\cdot \texttt{y}_j\subseteq \texttt{y}_j\cdot\texttt{B}$. On the other hand, for some  $\texttt{y}\in \texttt{y}_j\cdot  \texttt{B}$ with $\texttt{y}=\texttt{y}_j \cdot b$ we have   $\texttt{y}\in \im(f)$, where $f\in \Hom_{\texttt{B}}\left(\texttt{L}(j), \texttt{L}(1)\right)= \spann_K\left\{\texttt{f}_{(1,i,j)}\mid i\in \Lambda^{(j)}\right\}$ is given by $f:\texttt{L}(j)\stackrel{\cdot  b}{\rightarrow}\texttt{L}(1)=\texttt{B}$.  Since  $\im\left(\texttt{f}_{(1,i,j)}\right)= \texttt{L}(i) \subseteq \texttt{L}(j)$ for each  $i\in \Lambda^{(j)}$, we obtain   $\im (f)\subseteq \texttt{L}(j)$ and consequently $\texttt{y}\in \texttt{L}(j)=\texttt{B}\cdot \texttt{y}_j$. Thus  we obtain  $\texttt{B}\cdot \texttt{y}_j \supseteq \texttt{y}_j\cdot \texttt{B}$.

If  $i\neq j$, then  $\texttt{L}(i)\neq \texttt{L}(j)$, thus  $\texttt{L}(i)\subseteq \rad (\texttt{L}(j))$ for  any $i\in \Lambda^{(j)}\backslash\left\{j\right\}$. Consequently,  $\sum_{j<i}\texttt{L}(i)\subseteq \rad (\texttt{L}(j))$.  The set $\bigcup_{j<i}\texttt{X}(i)= \left\{\texttt{y}_i\mid i\in \Lambda^{(j)}\backslash \left\{j\right\}\right\}$  is a $K$-basis of $\sum_{j<i}\texttt{L}(i)$,  since  $\texttt{X}(i)$ is a $K$-basis of $\texttt{L}(i)$.
Thus 
 $\dimm_K\left(\sum_{j<i}\texttt{L}(i)\right) = \left|\Lambda^{(j)}\backslash\left\{j\right\}\right| = \left|\Lambda^{(j)}\right|-1= \dimm_K\texttt{L}(j)-1 = \dimm_K\left(\rad \texttt{L}(j)\right)$  implies  $\sum_{j<i}\texttt{L}(i)= \rad (\texttt{L}(j))$  for all $j\in \Lambda$.    \hfill $\Box$
 \\

Furthermore, for  a pair $(B,L)$ which  satisfies the condition \fbox{$\ma$}  we consider the algebra  $A=\End_{\texttt{B}}\left(\texttt{L}\right)\cong \End_B(L)^{op}$.   
    The evaluation functor  $\Hom_{\texttt{B}}(\texttt{L},-):\modd \texttt{B}\rightarrow \modd A$   provides  an isomorphism  $\Hom_{\texttt{B}}(\texttt{L}(i),\texttt{L}(j))\cong \Hom_A(P(i),P(j))$ (see   \cite[Proposition 2.1]{ARS}).  Moreover,   an injective $\texttt{B}$-map $\texttt{L}(i)\stackrel{\texttt{f}}{\hookrightarrow} \texttt{L}(j)$ induces an injective $A$-map $P(i)\stackrel{f}{\hookrightarrow} P(j)$, since $\Hom_{\texttt{B}}(\texttt{L},-)$ is left exact. 

 The previous  Lemma shows that   the properties  described in Lemma~\ref{L5} are  also  satisfied    for the  $\texttt{B}$-ideals  $\texttt{L}(i)$. 
 The part  \textit{(3)}   implies       $P(i)\hookrightarrow P(j)$  for all  $i\in \Lambda^{(j)}$, moreover, since $\im\left(\texttt{L}(i)\hookrightarrow \texttt{L}(j)\right)= \texttt{L}(i)$ for any injective  $\texttt{B}$-map,  we obtain  that   a submodule of $P(j)$ isomorphic to $P(i)$ is  uniquely determined. For  $i\in \Lambda^{(j)}$ we consider   $P(i)$  as a submodule~ of ~$P(j)$.

For all $i\in \Lambda^{(j)} \cap \Lambda^{(k)}$   let  $\ger{f}_{(j,i,k)} \in  \Hom_A(P(k),P(j))$  be the  map  induced by   $\texttt{f}_{(j,i,k)} \in \Hom_{\texttt{B}}(\texttt{L}(k),\texttt{L}(j))$   (described in  ~\ref{duality}).    Since  $\texttt{f}_{(j,i,k)}= \texttt{f}_{(j,i,i)}\circ \texttt{f}_{(i,i,k)}$  and   $\texttt{f}_{(j,i,i)}$  is injective,  we obtain   $\ger{f}_{(j,i,k)}=\ger{f}_{(j,i,i)}\circ \ger{f}_{(i,i,k)}:\left(P(k)\stackrel{\ger{f}_{(i,i,k)}}{\longrightarrow}P(i)\stackrel{\ger{f}_{(j,i,i)}}{\hookrightarrow}P(j)\right)$. 
 Obviously,  $\im \left(\ger{f}_{(j,i,k)}\right)$ belongs to the submodule  $ P(i)$ of $ P(j)$.

\begin{subnum}\begin{normalfont}\textbf{Remark-Notations.}
  According to Lemma~\ref{L5}\textit{(4)}
 the set  $\left\{\ger{f}_{(j,i,k)} \mid  i\in \Lambda^{(j)}\cap \Lambda^{(k)}\right\}$   is a $K$-basis of    $\Hom_A(P(k),P(j))$   for all  $j,k\in \Lambda$. 
  Thus   $\im(\ger{f})\subseteq \sum_{i\in \Lambda^{(j)}\cap \Lambda^{(k)}}  \im \left(\ger{f}_{(j,i,k)}\right)\subseteq   \sum_{i\in \Lambda^{(j)}\cap \Lambda^{(k)}}P(i) $ for all $\ger{f}\in \Hom_{A}(P(k),P(j))$.  Let $p(j,i,k):= \ger{f}_{(j,i,k)}(e_j)$, then we obtain that   $\left\{p(j,i,k)\mid i\in \Lambda^{(j)}\cap \Lambda^{(k)}\right\}$ is a $K$-basis of $P(j)_k$.  
  
  Let $\Gamma$ be some subset of $\Lambda^{(j)}$. Then   the following hold:
\\[5pt]
\hspace*{3mm}	(1) For any $l\in \Lambda$ the set  $\textbf{B}_j(\Gamma, l):=\left\{p(j,t,l)\mid t\in \bigcup_{i\in \Gamma}\left(\Lambda^{(i)}\cap \Lambda^{(l)}\right)\right\} $  is a  $K$-basis of the subspace  $\left(\sum_{i\in \Gamma}P(i)\right)_l$ of $P(j)_l$  for  the submodule $\sum_{i\in \Gamma}P(i)$ of $P(j)$:  Let $i\in \Gamma$, then $\left\{\ger{f}_{(i,t,l)}\mid t\in \Lambda^{(i)}\cap \Lambda^{(l)}\right\}$ is a $K$-basis of $\Hom_A(P(l),P(i))$.  By applying  $\Hom_A(P(l),-)$  to  $\ger{f}_{(j,i,i)}:P(i)\hookrightarrow P(j)$,  $e_i\mapsto p(j,i,i)$, we obtain  $\Hom_A(P(l),P(i))\hookrightarrow \Hom_A(P(l),P(j))$ with $\ger{f}_{(i,t,l)}\mapsto \ger{f}_{(j,i,i)}\circ \ger{f}_{(i,t,l)} =\ger{f}_{(j,t,l)}$ (or, equivalently, $P(i)_l\hookrightarrow P(j)_l$ with $p(i,t,l)\mapsto p(j,t,l)$). The set  $\left\{p(j,t,l)\mid t\in \Lambda^{(i)}\cap\Lambda^{(l)}\right\}$ is a $K$-basis of $P(i)_l\subseteq P(j)_l$ for any $i\in \Gamma$. Thus $\textbf{B}_j(\Gamma, l)=\bigcup_{i\in \Gamma}\left\{p(j,t,l)\mid t\in \Lambda^{(i)}\cap\Lambda^{(l)}\right\}$ is a $K$-basis of the subspace  $\left(\sum_{i\in \Gamma}P(i)\right)_l$ of $P(j)_l$. 
\\[5pt]	 
\hspace*{3mm}	(2) The set  $\textbf{B}_j(\Gamma):= \left\{p(j,t,l)\mid t\in \bigcup_{i\in \Gamma}\Lambda^{(i)}, \ l\in \Lambda_{(t)}\right\}$  is a $K$-basis of the submodule $\sum_{i\in \Gamma}P(i)$ of $P(j)$: Obviously,  $\textbf{B}_j(\Gamma) $ is  the disjoint union of    $\textbf{B}_j(\Gamma, l) $ for $l\in\Lambda$, thus  $\textbf{B}_j(\Gamma)$  is a $K$-basis of  $\sum_{i\in \Gamma} P(i)$. 
\end{normalfont} 
\label{basisrem}
\end{subnum}

We can now prove the   	four statements formulated 	at the beginning of this subsection, so that  the algebra $A$  with $(\Lambda,\ma)$ is 1-quasi-hereditary (recall that $1\ma i\ma n$ for all $i\in \Lambda$).   
\\

\textit{Proof. }\ding{182}  
 The definition of the standard modules  provides   $\De(j)=P(j)/N(j)$, where $N(j)= \sum_{k\not\ma j}\sum_{\ger{f}\in \Hom_{A}(P(k),P(j))}\im(\ger{f})$. The previous   	deliberations  imply  $N(j)\subseteq \sum_{j<i}P(i) $.  Since   $ P(i)=\im \left(\ger{f}_{(j,i,i)}\right) \subseteq N(j)$ for any $i\in \Lambda^{(j)}\backslash\left\{j\right\}$,  we obtain $N(j)\supseteq\sum_{j<i}P(i) $, thus  $\De(j)=P(j)/\left(\sum_{j<i}P(i) \right)$.

Let $\Gamma=\Lambda^{(j)}\backslash\left\{j\right\}$ and $k\in \Lambda_{(j)}$, then $\textbf{B}_j\left(\Lambda^{(j)},k\right)= \left\{p(j,j,k)\right\}\dot{\cup }\textbf{B}_j(\Gamma,k)$, using the notation of ~\ref{basisrem}(1). Since   $P(j)_k=\spann_K\textbf{B}\left(\Lambda^{(j)},k\right)$  and $\left(\sum_{j<i}P(i) \right)_k=\spann_K\textbf{B}_j(\Gamma,k)$, we obtain $1=\dimm_K\De(j)_k=\dimm_K\left(P(j)_k/\left(\sum_{j<i}P(i)\right)_k\right)=[\De(j):S(k)]$.

\ding{183} 
Let   $\textbf{i}\in \mathcal{L}(j):=\left\{(i_1, i_2, \ldots , i_r) \mid i_m\in \Lambda^{(j)} , \   i_{l}\not\geqslant i_{t}, \ 1\leq  l < t \leq r :=\left|\Lambda^{(j)}\right|\right\}$ 
 (see \cite[4.2]{P}). Obviously, $P(i_t)\subseteq P(j)=P(i_1)$ for all $1\leq t\leq r$.    Denote by $\mathscr{D}(\textbf{i})$  the filtration  $0=D(r+1)\subset D(r)\subset \cdots \subset D(t) \subset \cdots \subset D(1) $      with   $D(t):=\sum_{m=t}^{r}P(i_m)$. 
 It is easy to check that 
  $\textbf{B}_j\left(\Lambda^{(i_t)}\backslash \left\{i_t\right\}\right)= \textbf{B}_j\left(\left\{i_t\right\}\right)\cap\textbf{B}_j\left(\left\{i_{t+1},\ldots, i_r\right\}\right)$ for all $1\leq t\leq r-1$. This  implies $\sum_{i_t<k}P(k)= P(i_t)\cap\left(\sum_{m=t+1}^rP(i_m)\right)$  (see ~\ref{basisrem}(2)) and consequently $D(t)/D(t+1)\cong P(i_t)/\left(\sum_{i_t<k}P(k)\right) \cong \De(i_t)$ for all $1\leq t\leq r$. The filtration $\mathscr{D}(\textbf{i})$ of $P(j)$ is $\De$-good. Since $\left\{i_1, i_2, \ldots , i_r\right\}=\Lambda^{(j)}$ and   $l\neq t$ implies $i_l\neq i_t$, we obtain $\left(P(j):\De(i)\right)=1$ for any $i\in \Lambda^{(j)}$ and $\left(P(j):\De(i)\right)=0$ if $i\in \Lambda\backslash \Lambda^{(j)}$.

 \ding{184} Since $P(i)\hookrightarrow P(1)$ for all $i\in \Lambda$,   it is enough to show that  $\soc P(1)\cong S(1)$.  We consider the map  $\ger{f}_{(1,n,1)}:P(1)\rightarrow P(1)$  induced by  $\texttt{f}_{(1,n,1)}:\texttt{L}(1)\sur \texttt{L}(n)\hookrightarrow \texttt{L}(1)$, here  $\texttt{L}(n)=\soc \texttt{L}(1)$ since $n$ is maximal. 
 
  We show that 
  $\im \left(\ger{f}_{(1,n,1)}\right) \subseteq \im \left(\ger{f}\right)$ for all  $\ger{f}\in \Hom_A(P(i),P(1))$ with $\ger{f}\neq 0$ and  all $i\in \Lambda$.
 This implies, that $\im \left(\ger{f}_{(1,n,1)}\right)$  is contained  in  every local submodule of $P(1)$ and therefore in  every non-zero submodule of $P(1)$. Thus  $\im \left(\ger{f}_{(1,n,1)}\right)$ is the uniquely determined simple submodule of $P(1)$ and $\topp \left(\im \left(\ger{f}_{(1,n,1)}\right)\right)\cong S(1)$ implies $\im \left(\ger{f}_{(1,n,1)}\right)\cong S(1)$: 
Let $\ger{f}\in \Hom_A(P(i),P(1))\backslash\left\{0\right\}$ be induced by $\texttt{f}\in \Hom_{\texttt{B}}\left(\texttt{L}(i),\texttt{L}(1)\right)$, then  $\texttt{f}\neq 0$   and consequently  $\texttt{L}(n)\subseteq \im \left(\texttt{f}\right)$.  There exists some  $\texttt{x}\in \texttt{L}(i)$ with $\texttt{f}(\texttt{x})=\texttt{b}_n$, where $\texttt{b}_n$ is a generator of $\texttt{L}(n)$. Let $\texttt{g}$  be in $ \Hom_{\texttt{B}}(\texttt{L}(1),\texttt{L}(i))$   given by  $\texttt{g}:\texttt{L}(1)\stackrel{\cdot \texttt{x}}{\rightarrow}\texttt{L}(i)$     and $\ger{g}\in \Hom_A(P(1),P(i))$ is induced by $\texttt{g}$.   We have $\texttt{f}\circ \texttt{g}=\texttt{f}_{(1,n,1)}$. This implies $\ger{f}\circ \ger{g}=\ger{f}_{(1,n,1)}$ and consequently $\im \left(\ger{f}_{(1,n,1)}\right)\subseteq \im \left(\ger{f}\right)$.

   According to ~\ref{maop} for the algebra $A^{op}\cong \End_{B}(L)\cong \End_{\texttt{B}^{op}}(\DD(\texttt{L}))$ we have  $\soc P_{A^{op}}(i)\cong S_{A^{op}}(1)$, thus $\topp I(i)\cong S(1)$ for all $i\in \Lambda$.

\ding{185}  Let   $\ger{f}_{(n,n,j)}:P(j)\rightarrow P(n)$  be the $A$-map  induced by the   $\texttt{B}$-map $\texttt{f}_{(n,n,j)}:\texttt{L}(j)\sur \texttt{L}(n)$.   It is enough 
to show the equation $\sum_{j<i}P(i)= \Kern \left(\ger{f}_{(n,n,j)}\right)$.  This  implies       $P(j)/\left(\sum_{j<i}P(i)\right) \stackrel{\text{\ding{182}}}{=} \De(j)\hookrightarrow P(n)=\De(n)$ for any $j\in \Lambda$: 
 Let  $i\in \Lambda^{(j)}\backslash\left\{j\right\}$, then  the functor $\Hom_{\texttt{B}}(\texttt{L},-)$ maps a $\texttt{B}$-map  $\texttt{f}:\left(\texttt{L}(i)\stackrel{\texttt{f}_{(j,i,i)}}{\hookrightarrow}\texttt{L}(j)\stackrel{\texttt{f}_{(n,n,j)}}{\sur}\texttt{L}(n)\right)$ to the $A$-map  $\ger{f}:\left(P(i)\stackrel{\ger{f}_{(j,i,i)}}{\hookrightarrow}P(j)\stackrel{\ger{f}_{(n,n,j)}}{\rightarrow}P(n)\right)$. Since  $\texttt{L}(n)=\soc B$ is simple,  we have $\Kern\left(\texttt{f}_{(n,n,j)}\right)=\rad \texttt{L}(j)$, thus  $\texttt{L}(i)\subseteq \rad \texttt{L}(j)$ since $i\in \Lambda^{(j)}\backslash \left\{j\right\}$. Hence $\texttt{f}$ and therefore $\ger{f}$   are  zero-maps. Consequently,  the submodule $P(i)$  of $P(j)$ belongs to  $\ker \left(\ger{f}_{(n,n,j)}\right)$ for any $i\in \Lambda^{(j)}\backslash\left\{j\right\}$.  We obtain  $\sum_{j<i}P(i)\subseteq \ker \left(\ger{f}_{(n,n,j)}\right)$. 
 
 Because  $\De(j)=P(j)/\left(\sum_{j<i}P(i)\right)$,  there exists a submodule $U$ of $\De(j)$  such that   $P(j)/\left(\ker \ger{f}_{(n,n,j)}\right) \cong \De(j)/U$.  
 For the $K$-subspace of $\im( \ger{f}_{(n,n,j)})$ corresponding to some  $k\in \Lambda$ we have 
 \    $\dimm_K \left(\im \ger{f}_{(n,n,j)}\right)_k = \dimm_K\left(P(j)/\ker\ger{f}_{(n,n,j)}\right)_k = \dimm_K\left(\De(j)/U \right)_k \leq \dimm_K\left(\De(j)\right)_k. $

 Let $k\in \Lambda_{(j)}$,   then the $\texttt{B}$-map   $\texttt{g}:\left(\texttt{L}(k) \stackrel{\texttt{f}_{(j,j,k)}}{\sur}\texttt{L}(j)\stackrel{\texttt{f}_{(n,n,j)}}{\sur}\texttt{L}(n)\right)$  is non-zero, thus the induced $A$-map $\ger{g}:\left(P(k) \stackrel{\ger{f}_{(j,j,k)}}{\longrightarrow}P(j)\stackrel{\ger{f}_{(n,n,j)}}{\longrightarrow}P(n)\right)$  is non-zero. Hence  $0\neq \ger{g}(e_k)\in \left(\im \ger{f}_{(n,n,j)}\right)_k$, thus $\dimm_K\left(\im \ger{f}_{(n,n,j)}\right)_k\neq 0$ and consequently $\left(\im \ger{f}_{(n,n,j)}\right)_k=\left(\De(j)\right)_k$ for all $k\in \Lambda$, because  $\dimm_K\left(\De(j)\right)_k\stackrel{\text{\ding{182}}}{=}\begin{small}\left\{
\begin{array}{cl}
1 & \text{if  } k\in \Lambda_{(j)},\\
0 & \text{else}
\end{array}
\right.\end{small}$.  We obtain $U=0$ and  $\sum_{j<i}P(i)= \Kern \left(\ger{f}_{(n,n,j)}\right)$.

Proposition~\ref{maop} implies that the  functor  $\Hom_{B}(L,-):\modd B\rightarrow \modd A^{op}$ 
 yields  $\De_{A^{op}}(j)\hookrightarrow \De_{A^{op}}(n)$. By applying the duality  we get  $\nabla(n)\sur \nabla(j)$ for all $j\in \Lambda$. 
 
 This finishes the proof of Theorem A. \hfill $\Box$ 
\\


The features of the  $B$-module $L$, for a pair  $(B,L)$ with  \fbox{$\ma$}, implies some properties  for the quiver and relations of the corresponding  1-quasi-hereditary algebra $\End_B(L)^{op}$.

 \begin{subnum}\begin{normalfont}\textbf{Remark.} Let $(B,L)$ be a pair with the property \fbox{$\ma$}.  For every $i\in \Lambda$ we fix a generator $x_i\in B$ of the direct summand  $L(i)$ of $L$   and  $x_1=\texttt{1}:=\texttt{1}_B$, here $\left\{1\right\}=\min\left(\Lambda,\ma\right)$.  For all $i,j\in \Lambda$ with 
  $j<i$ there exists   $b_{ij}\in B$ with $x_j\cdot b_{ij}=x_i$, because    $x_i\in L(i)\subset L(j)=B\cdot x_j=x_j\cdot B$ (see ~\ref{L5}\textit{(1)}and \textit{(3)}).
  
 We define  an injective and a surjective $B$-map between $L(i)$ and $L(j)$ by 
\begin{center}
$f_{(i\to j)}:L(i)\stackrel{\cdot \texttt{1}}{\hookrightarrow} L(j)$\   \    \   \    and   \    \  \   \     $f_{(i\leftarrow j)}:L(j)\stackrel{\cdot b_{ij}}{\sur} L(i)$
\end{center}
For  any $l,t\in \Lambda$ 
let $\ger{X}(L(l),L(t))$  be the space of non-invertible maps  $f\in \Hom_{B}\left(L(l),L(t)\right) $    with the property  \dq \ if $  f =f_2\circ f_1 $   factors through  $  \add L$,   then  either  $ f_1$ is a split monomorphism  or   $f_2$   is a  split epimorphism\dq.   The number of arrows from $l$ to $t$  in the quiver of the algebra  $A=\End_B(L)^{op}$ is  $\dimm_K\ger{X}(L(l),L(t))$ (see \cite{ARS} or \cite{ASS}).

 According to ~\ref{L5}\textit{(4)}, any map $f\in \Hom_{B}\left(L(l),L(t)\right)$ factors through $\bigoplus_{i\in \Lambda^{(l)}\cap \Lambda^{(t)}}L(i)$. If $l$ and $t$ are incomparable, then $l,t\not\in  \Lambda^{(l)}\cap \Lambda^{(t)}$, thus $\ger{X}(L(l),L(t))=0$. Assume $l<t$ and $f:L(l) \rightarrow  L(t)$, then there exists $b\in B$ with $x_t\cdot b = f(x_l)$,  since $L(t)=B\cdot x_t\stackrel{\ref{L5}\textit{(1)}}{=}x_t\cdot B$ and hence   $f:\left(L(l) \sur L(t) \stackrel{ \cdot b}{\rightarrow}L(t) \right)$.  If $f$ is not surjective, then $b$ is not invertible and consequently       $f\not\in  \ger{X}(L(l),L(t))$.   If $f$ is surjective but $l$ and $t$ are not the adjacent, then $f:L(l)\sur L(j) \sur L(t)$  for some $j\in \Lambda$ with $l<j<t$ and therefore $f\not \in  \ger{X}(L(l),L(t))$. Let $l\triangleleft t$ and  $g: L(l)\sur L(t)$ with $g:\left(L(l)\sur L(t) \stackrel{\cdot b}{\sur}L(t)\right)$, then $L(t) \stackrel{\cdot b}{\sur}L(t)$ is a  split  epimorphism  if $b= c\cdot \texttt{1}$ for some $c\in K\backslash \left\{0\right\}$, in other words  $g=c\cdot f_{(t\leftarrow l)}$. Using  analogous arguments also  for  $j> t$,  we obtain      $\ger{X}(L(l),L(t))= \begin{small}\left\{
\begin{array}{cl}
\spann_K\left\{f_{(l\to t)}\right\} & \text{if  } l\triangleright  t, \\
\spann_K\left\{f_{(t\leftarrow l)}\right\} & \text{if  } l\triangleleft t, \\
0 & \text{else}.
\end{array}
\right. \end{small}$\\

(1)  In the quiver $Q_0(A)=\Lambda$ of the 1-quasi-hereditary algebra  $A$  two vertices $i$ and $j$ are connected  by an arrow if they are  neighbours   with respect to  $\ma$, more precisely,  we have   $i \rightleftarrows j$. Assume $j\triangleleft i$, then   the $B$-maps  $f_{(i\to j)}$ and  $f_{(i\leftarrow j)}$   can by considered as the maps corresponding to the  arrows $i\to j$ and $j\to i$ respectively. In this case we  use the  notation    
 $L(i) \overset{\cdot \texttt{1}}{\underset{\cdot b_{ij}}{\rightleftarrows}} L(j)$. In general the notation  $L(l) \overset{\cdot b}{\underset{\cdot d}{\rightleftarrows}} L(t)$  means that $l$ and $t$ are neighbours  and $(b,d)=\begin{small}\left\{\begin{array}{ll}
(\texttt{1},b_{lt}) & \text{if }  t\triangleleft l, \\
(b_{tl}, \texttt{1}) & \text{if }  t\triangleright l. 
\end{array}
\right. \end{small}$   We    always     have     $(x_l \cdot b, x_t\cdot d)=\begin{small} \left\{
\begin{array}{ll}
(x_l, x_l) & \text{if   }  t\triangleleft l, \\
(x_t, x_t) & \text{if   }  t\triangleright l. 
\end{array}
\right. \end{small}$

(2)  Let $p_t=(i,i_1^{(t)},\ldots , i_{m_t}^{(t)},j)$ for $1\leq t\leq r$  be some paths in $Q(A)$  (obviously, $i_{k}^{(t)}$ and $i_{k+1}^{(t)}$ are neighbours). Then
\\[5pt]
 $\displaystyle\sum_{t=1}^r c_t\cdot p_t\in \II(A)$  \  if and only if \   $\displaystyle \sum_{t=1}^{r}c_t\cdot \left(L(i)\stackrel{\cdot b_0^{(t)}}{\longrightarrow} L(i_1^{(t)}) \stackrel{\cdot b_1^{(t)}}{\longrightarrow} \cdots \stackrel{\cdot b_{m_t-1}^{(t)}}{\longrightarrow} L(i_{m_t}^{(t)}) \stackrel{\cdot b_{m_t}^{(t)}}{\longrightarrow} L(j) \right) = 0$,
\\[5pt]
    here the maps  $L(i_{k}^{(t)}) \stackrel{\cdot b_{k}^{(t)}}{\longrightarrow}L(i_{k+1}^{(t)})$ are of the form $f_{(l\to t)}$ or $f_{(l\leftarrow t)}$.
\end{normalfont} 
\label{streles}
\end{subnum}

 \begin{subnum}\begin{normalfont}\textbf{Lemma.}\end{normalfont} Let $A\cong \End_B\left(L\right)^{op}\cong KQ/\II$ with $(\Lambda,\ma)$  be a 1-quasi-hereditary algebra, where $\left(B,L=\bigoplus_{i\in \Lambda}L(i)\right)$ satisfies the property \fbox{$\ma$}.     
\begin{itemize}
	\item[(1)] If $p$ and $q$ are some paths in $Q$ of the form $p(j,i,k)$ (see Subsection 2.1),  then $p-q\in \II$.
	In particular, $\ger{o}(p(j,i,k))=p(k,i,j)$. 
	\item[(2)] Let $\Gamma$ be the set of   (large) neighbours of $1$, where $\left\{1\right\}=\min(\Lambda,\ma)$ and $x_i$ be a generator of $L(i)$ for any $i\in \Gamma$. Then the 
	 set   $\left\{x_i\mid i\in \Gamma\right\}$  is   generating system of $B$. 
	 In particular,  
	 $B$ is   a factor algebra of $K\left\langle y_1,\ldots ,y_m\right\rangle$, where $m=\left|\Gamma\right|$. 
\end{itemize}
\label{polynom}
\end{subnum}

\textit{Proof.} \textit{(1)}  Let $j<i$ and  $p$,  $q$  be some increasing paths from $j$ to $i$ as well as $\ger{o}(p)$, $\ger{o}(q)$ be the corresponding decreasing paths from $j$ to $i$ in $Q$, i.e., there exists   $i=i_0\triangleleft i_1\triangleleft \cdots \triangleleft i_m=j$ and $i=j_0\triangleleft j_1 \triangleleft \cdots \triangleleft j_r=j$  with $p=(i,i_1,\ldots, i_m,j)$ and $q=(i,j_1,\ldots, j_r,j)$
 as well as $\ger{o}(p)=(j,i_m,\ldots, i_1,i)$ and $\ger{o}(q)=(j,j_r,\ldots,j_1,i)$.  For the corresponding $B$-maps $f_{(p)}=f_{(i_m\to j)} \circ  \cdots \circ f_{(i_1\to i_2)} \circ f_{(i\to i_1)}$ and  $f_{(q)}=f_{(j_r\to j)} \circ \cdots \circ f_{(j_1\to j_2)} \circ f_{(i\to j_1)}$ as well as $f_{\ger{o}(p)}=f_{(i\leftarrow i_1)} \circ  f_{(i_1\leftarrow i_2)} \circ \cdots  \circ f_{(i_m\leftarrow j)}  $
 and $f_{\ger{o}(q)}=f_{(i\leftarrow j_1)} \circ  f_{(j_1\leftarrow j_2)} \circ \cdots  \circ f_{(j_r\leftarrow j)}$ we obtain $f_{(p)}- f_{(q)}=0$ and $f_{\ger{o}(p)}-f_{\ger{o}(q)}=0$.

For some $i,j,k\in \Lambda$ with $i \ma j,k$ let  $p$ and $q$ be some paths in $Q$ of the form $p(j,i,k)$, then $p=\ger{o}(p_1)\cdot p_2$ and $q=\ger{o}(q_1)\cdot q_2$ with some increasing paths  $p_1,q_1$ from $k$ to $i$ and $p_2,q_2$ from $j$ to $i$. 
For the corresponding $B$-maps we have $f_{(p)}=f_{(\ger{o}(p_1))}\circ f_{(p_2)}$ and   $f_{(q)}=f_{(\ger{o}(q_1))}\circ f_{(q_2)}$. Since $f_{(p_1)}=f_{(q_1)}$ and  $f_{(\ger{o}(p_1))}=f_{(\ger{o}(q_1))}$,  we obtain $f_{(p)}= f_{(q)}$. This implies $p-q\in \II$.

\textit{(2)} 
For any $i\in \Lambda$ with $i\neq 1$ there exists $j\in \Gamma$ with $j\ma i$.  Thus $L(i)\subseteq L(j)$ and consequently $\rad B= \sum_{i\in \Lambda\backslash\left\{1\right\}}L(i)=\sum_{i\in \Gamma}L(i)$ (see ~\ref{property_ma}(2)(b)). The set $\left\{x_i\mid i\in \Gamma\right\}$ generates $\rad B$. Since $B$ is local, we obtain that $B$ is a factor algebra of $K$$\left\langle y_1,\ldots , y_{\left|\Gamma
\right|}\right\rangle$. \hfill $\Box$

\section{Proof of Theorem  B}

In this section $A=KQ/\II$ with $(\Lambda,\ma)$ is a 1-quasi-hereditary algebra and $(B,L=\bigoplus_{i\in \Lambda}L(i))$  is the corresponding pair with the property \fbox{$\ma$}, i.e., $A\cong \End_B(L)^{op}$. For a relation $\rho = \sum_{t=1}^{r}c_t\cdot p_t$ in $\II$ we define  $\ger{o}(\rho)=\sum_{t=1}^{r}c_t\cdot \ger{o}(p_t)$. The definition of a $\BGG_{(\leftrightarrows)}$-algebra is given in Subsection 1.3.

  For the proof of Theorem B we  have to show the  equivalence of the    following statements:
\begin{center}
\ding{182} $B$  is commutative, \   \  
	\ding{183}  $A$ is a $\BGG_{(\leftrightarrows)}$-algebra, \   \  
	\ding{184}  $\rho \in \II$ if and only if $\ger{o}(\rho ) \in \II$
\end{center}
Let $x_i$ be a generator of $L(i)$ for any $i\in \Lambda$ and $x_1=\texttt{1}:= \texttt{1}_{B}$ where $\left\{1\right\}=\min(\Lambda,\ma)$. For  $i,j\in \Lambda$ with $i\triangleleft j$ or $i\triangleright j$  we denote   by $L(i) \overset{\cdot b}{\underset{\cdot d}{\rightleftarrows}}   L(j)$ the $B$-maps  described in ~\ref{streles}.

\begin{num}\begin{normalfont}\textbf{Lemma.}\end{normalfont} 
The
following statements are equivalent:
\begin{itemize}
	\item[(i)] $B$ is commutative.
	\item[(ii)] Let $p=(i_0,i_1,\ldots, i_m)$ be some path in $Q$ and $L(i_0) \overset{\cdot \texttt{b}_1}{\underset{\cdot \texttt{d}_1}{\rightleftarrows}}   L(i_1) \overset{\cdot \texttt{b}_2}{\underset{\cdot \texttt{d}_2}{\rightleftarrows}} \cdots \overset{\cdot \texttt{b}_{m}}{\underset{\cdot \texttt{d}_{m}}{\rightleftarrows}} L(i_{m})$  be  the corresponding   $B$-maps with $x_0:=x_{i_o}$ and   $x_m:=x_{i_m}$.  Then 
\begin{center}
$x_0\cdot \texttt{b}_1\cdot \texttt{b}_2\cdots \texttt{b}_{m} = x_m\cdot \texttt{d}_{m} \cdots \texttt{d}_2\cdot \texttt{d}_1$
\end{center}
	\item[(iii)]  Statement $(ii)$ holds for $m=4$.
\end{itemize}
\label{komutativ}
\end{num}

\textit{Proof.}  $(i)\Rightarrow (ii)$  We show this by induction on $m$: If  $m=1$, then  
   for  $L(i_0)\overset{\cdot \texttt{b}_1}{\underset{\cdot \texttt{d}_1}{\rightleftarrows}} L(i_1)$ 
 we have  $\left(x_0\cdot \texttt{b}_1,  x_1\cdot  \texttt{d}_1\right) \in \Big\{\left(x_0,x_0\right), \left(x_1,x_1\right)\Big\}$ (see  ~\ref{streles}(1)), thus $x_0\cdot \texttt{b}_1 =  x_1\cdot  \texttt{d}_1$. 
   Assume  $x_0\cdot \texttt{b}_1\cdot \texttt{b}_2\cdots \texttt{b}_{m-1} = x_{m-1}\cdot \texttt{d}_{m-1} \cdots \texttt{d}_2\cdot \texttt{d}_1$, then by multiplication with $\texttt{b}_m$ we obtain
\begin{center}
\hspace*{30mm}$x_0\cdot \texttt{b}_1\cdot \texttt{b}_2\cdots \texttt{b}_{m-1}\cdot \texttt{b}_m = x_{m-1}\cdot \texttt{b}_m \cdot \texttt{d}_{m-1} \cdots \texttt{d}_2\cdot \texttt{d}_1$, \hspace{30mm} $(\circledast)$
\end{center}
 because $B$ is commutative. 
  If $i_m \triangleleft i_{m-1}$, then $\texttt{b}_m=\texttt{1}$ and $x_m\cdot \texttt{d}_m=x_{m-1}$. The equation  $(\circledast)$  implies $x_0\cdot \texttt{b}_1\cdot \texttt{b}_2\cdots \texttt{b}_{m-1}\cdot \texttt{b}_m = x_{m}\cdot \texttt{d}_m\cdot \texttt{d}_{m-1} \cdots \texttt{d}_2\cdot \texttt{d}_1$.  If  $i_m \triangleright i_{m-1}$,  we obtain  $x_{m-1}\cdot \texttt{b}_m= x_m$ and $\texttt{d}_m=\texttt{1}$, thus the equation $(\circledast)$ \  is \  $x_0\cdot \texttt{b}_1\cdot \texttt{b}_2\cdots \texttt{b}_{m-1}\cdot \texttt{b}_m = x_{m}\cdot \texttt{d}_m\cdot \texttt{d}_{m-1} \cdots \texttt{d}_2\cdot \texttt{d}_1$.

  $(ii)\Rightarrow (iii)$ This is trivial.

$(iii)\Rightarrow (i)$  Let $i,j\in \left\{l\in \Lambda\mid 1\triangleleft l\right\}$, then for   $L(1)  \overset{\cdot \texttt{b}_1}{\underset{\cdot \texttt{d}_1}{\rightleftarrows}} L(i)  \overset{\cdot \texttt{b}_2}{\underset{\cdot \texttt{d}_2}{\rightleftarrows}} L(1) \overset{\cdot \texttt{b}_3}{\underset{\cdot \texttt{d}_3}{\rightleftarrows}} L(j)  \overset{\cdot \texttt{b}_4}{\underset{\cdot \texttt{d}_4}{\rightleftarrows}} L(1)$,   since $1\triangleleft i \triangleright 1 \triangleleft j \triangleright 1$,   we have $(\texttt{b}_1, \texttt{d}_1)=(x_i,\texttt{1})$,  $(\texttt{b}_2, \texttt{d}_2)=(\texttt{1},x_i)$, $(\texttt{b}_3, \texttt{d}_3)=(x_j,\texttt{1})$, $(\texttt{b}_4, \texttt{d}_4)=(\texttt{1},x_j)$ (see ~\ref{streles}(1)). By the assumption,   we have $x_1\cdot \texttt{b}_1 \cdot \texttt{b}_2 \cdot \texttt{b}_3 \cdot \texttt{b}_4 = x_1\cdot \texttt{d}_4 \cdot \texttt{d}_3 \cdot \texttt{d}_2\cdot \texttt{d}_1$, thus 
  we obtain $x_i\cdot x_j=\texttt{1} \cdot x_i \cdot \texttt{1} \cdot x_j\cdot \texttt{1} =\texttt{1}\cdot x_j\cdot \texttt{1} \cdot x_i\cdot \texttt{1} = x_j\cdot x_i$. Thus  $B$ is commutative,  because $\left\{x_i\mid 1\triangleleft i\right\}$ is  a  	 generating system of $B$ (see ~\ref{polynom}(2)).
 \hfill $\Box$ 
\\

\textit{Proof} \ding{182}$\Rightarrow$\ding{184}
Let  $ \rho=\sum_{t=1}^{r}c_t\cdot \left(i_0^{(t)}, i_1^{(t)}, \ldots,  i_{m_t}^{(t)} \right)$ be a relation of $A$   with $i=i_0^{(t)}$,  $j = i_{m_t}^{(t)} $  and  $L(i_{v-1}^{(t)}) \overset{\cdot \texttt{b}_v^{(t)}}{\underset{\cdot \texttt{d}_v^{(t)}}{\rightleftarrows}}   L(i_{v}^{(t)})$   the corresponding  $B$-maps for all $1\leq t\leq r$.  We obtain  $\sum_{t=1}^{r}c_t\cdot \left(L(i)\stackrel{\cdot b_1^{(t)}}{\longrightarrow} L(i_1^{(t)}) \stackrel{\cdot b_2^{(t)}}{\longrightarrow} \cdots  \stackrel{\cdot b_{m_t}^{(t)}}{\longrightarrow} L(j) \right) = 0$ (see ~\ref{streles}(2)). Hence,  
  $\sum_{t=1}^{r}c_t\cdot \left(x_i\cdot b_1^{(t)} \cdot \cdots \cdot b_{m_t}^{(t)}\right) = 0$.
According to  Lemma ~\ref{komutativ} we have    $\sum_{t=1}^{r}c_t\cdot \left(x_j\cdot d_{m_t}^{(t)} \cdot \cdots \cdot d_{1}^{(t)}\right) = 0$.  This  implies $\sum_{t=1}^{r}c_t\cdot \left(L(j)\stackrel{\cdot d_{m_t}^{(t)}}{\longrightarrow}  \cdots \stackrel{\cdot d_{2}^{(t)}}{\longrightarrow} L(i_{1}^{(t)}) \stackrel{\cdot d_1^{(t)}}{\longrightarrow} L(i) \right) = 0$  and consequently $\ger{o}(\rho)=\sum_{t=1}^{r}c_t\cdot \left(i_{m_t}^{(t)},  \cdots , i_{1}^{(t)},  i_{0}^{(t)} \right) \in \II$. 

\ding{184}$\Rightarrow$\ding{182} Let $\Gamma:=\left\{k\in \Lambda \mid 1\triangleleft k\right\}$.  It is enough  to show $x_k\cdot x_j=x_j\cdot x_k$ for all  $k,j \in \Gamma$ (see ~\ref{polynom}):  Let $k,j\in \Gamma$, then 
for the $B$-maps  $f$$:$$\left(L(1)\stackrel{\cdot x_k}{\sur} L(k)\stackrel{\cdot 1}{\hookrightarrow} L(1)\stackrel{\cdot x_j}{\sur} L(j)\stackrel{\cdot 1}{\hookrightarrow} L(1)\right)$ and $g$$:$$\left(L(1)\stackrel{\cdot x_j}{\sur} L(j)\stackrel{\cdot 1}{\hookrightarrow} L(1)\stackrel{\cdot x_k}{\sur} L(k)\stackrel{\cdot 1}{\hookrightarrow} L(1)\right)$ we obtain $f=g$ if and only if $f(\texttt{1})=x_k\cdot x_j=x_j\cdot x_k=g(\texttt{1})$.  Since $f$  and     $g$ correspond  to the paths   $(1,k,1,j,1)$ and   $(1,j,1,k,1)$ respectively,  we have to show $(1,k,1,j,1)-(1,j,1,k,1)\in \II$ (or $(1,j,1,k,1)=(1,k,1,j,1)$ in $A$).
 
According to  \cite[Theorem 3.2]{P} for the path $(j,1,k)$ there exists  some $c_i\in K$ with  $\rho=(j, 1, k)-\sum_{i\in \Lambda^{(j)}\cap \Lambda^{(k)}} c_i \cdot p(j,i,k)\in \II$.  Since $\ger{o}(j,1,k)=(k,1,j)$ and $\ger{o}(p(j,i,k))=p(k,i,j)$, by our   assumption   we obtain   
$\ger{o}(\rho)=(k,1,j) - \sum_{i\in \Lambda^{(j)}\cap \Lambda^{(k)}} c_i \cdot p(k,i,j)\in \II$. 
Obviously, $(1,j,1,k,1)= (k\to 1)\cdot (j, 1, k) \cdot (1\to j)$ and $(1,k,1,j,1)= \ger{o}(1,j,1,k,1)= (j\to 1)$$\cdot $$ (k, 1, j)$$ \cdot$$ (1\to k) $. 
The relations $\rho$ and $\ger{o}$$(\rho)$ implies the following  equations  in $A$.
\begin{center}
$
\begin{array}{ccl}
(1,j,1,k,1) & = &  \displaystyle \sum_{i\in \Lambda^{(j)}\cap \Lambda^{(k)}} c_i \cdot \left((k\to 1)\cdot  p(j,i,k)\cdot (1\to j)\right) \\
(1,k,1,j,1) & = & 
 \displaystyle \sum_{i\in \Lambda^{(j)}\cap \Lambda^{(k)}} c_i \cdot \left((j\to 1)\cdot  p(k,i,j)\cdot (1\to k)\right)   \\
\end{array} 
$
\end{center}
\parbox{14cm}{
For every   $i\in \Lambda^{(j)}\cap \Lambda^{(k)}$ the paths  $p_{(i)}= (k\to 1)\cdot  p(j,i,k)\cdot (1\to j) $ and  $q_{(i)}=(j\to 1)\cdot  p(k,i,j)\cdot (1\to k)$  are   of the form  $p(1,i,1)$
    (in the right picture the black and the  gray  path respectively). Thus  Lemma ~\ref{polynom}(1) implies $p_{(i)}=q_{(j)}$ in $A$, hence  $(1,j,1,k,1)=(1,k,1,j,1)$.
 }\psset{xunit=0.48mm,yunit=0.47mm,runit=1mm}
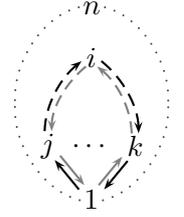
\begin{pspicture}(0,0)(0,0)
\rput(30,30){\begin{small}
\rput(0,10){\rnode{A}{$n$}}
\rput(0,-44){\rnode{B}{$1$}}
\rput(20,-17){\rnode{aa}{}}
\rput(-20,-17){\rnode{bb}{}}
              \rput(0,-4){\rnode{3}{$i$}}
              \rput(-12,-29){\rnode{6}{$j$}}
              \rput(0,-29){\rnode{66}{$\cdots $}}
              \rput(12,-29){\rnode{7}{$k$}}
              \end{small}
                      }
                      
                      \psset{nodesep=1pt,offset=1.3pt,arrows=->}
              \ncarc[arcangle=25,linestyle=dashed]{6}{3}
               \ncarc[arcangle=-25,linestyle=dashed,linecolor=gray]{3}{6}
              \ncarc[arcangle=25,linestyle=dashed]{3}{7}
              \ncarc[arcangle=-25,linestyle=dashed,linecolor=gray]{7}{3}
             \ncline[linecolor=gray]{6}{B}
              \ncline{B}{6}
              \ncline{7}{B}
             \ncline[linecolor=gray]{B}{7}
                      
             \psset{nodesep=1pt,offset=1.3pt}
               \ncline{3}{5}
               \ncline{5}{3}
               \ncline{5}{6}
               \ncline{6}{5}
               \ncline{5}{7}
               \ncline{7}{5}
               \ncline{5}{4}
               \ncline{4}{5}
              \ncline[linestyle=dotted]{21}{22}
              \ncarc[arcangle=40,linestyle=dotted,linecolor=black!70]{-}{A}{aa}
              \ncarc[arcangle=40,linestyle=dotted,linecolor=black!70]{-}{aa}{B}
             \ncarc[arcangle=40,linestyle=dotted,linecolor=black!70]{-}{B}{bb}
              \ncarc[arcangle=40,linestyle=dotted,linecolor=black!70]{-}{bb}{A}
\end{pspicture}
\\ 

\ding{183}$\Leftrightarrow$\ding{184} If the  $K$-map $\ger{o}:A\to A$ with $p\mapsto \ger{o}(p)$  is  an anti-automorphism of $A$ and 
 for some paths  $p_1,\ldots , p_r$, which start in $i$ and end in $j$ we have  $\sum_{t=1}^rc_t\cdot p_t=0$, then  $\sum_{t=1}^rc_t\cdot \ger{o}(p_t)=0$ (in other words, if $\rho\in \II$, then $\ger{o}(\rho)\in \II$ ). 
 
 On the other hand, if  $\II=\left\{\ger{o}(\rho)\mid \rho\in \II\right\}$, then $\ger{o}$ is obviously an  anti-automorphism. \hfill $\Box$

\begin{num}\begin{normalfont}\textbf{Example.} The pair $(B,L_{(i)})$ in the Example ~\ref{exampleB} provides a $\BGG_{(\leftrightarrows)}$-algebra $A_i=\End_B\left(L_{i}\right)^{op}$ for any $i=1,2,3$, because $B$ is commutative.
The quiver and relations to the left and to the right  present the algebra $A_1$ and $A_2$ respectively  (the algebra $A_3$ is presented in ~\cite[Example 1]{P1}). For any relation $\rho$ of $A_i$ also $\ger{o}(\rho)$ is a relation.
\end{normalfont} 
\label{Examp_relat}
\end{num}

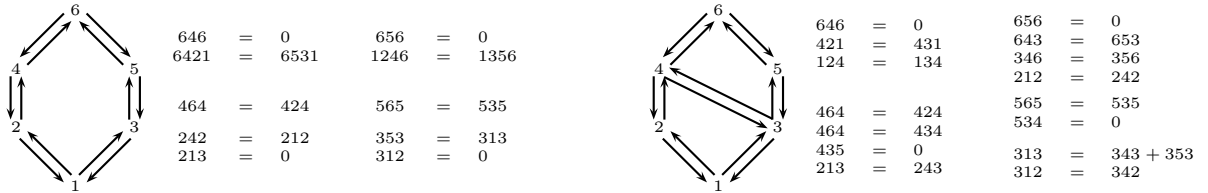
\text{ }\\

\psset{xunit=2.6mm,yunit=2.6mm,runit=5mm}
\begin{pspicture}(9,0)(10,0)
\rput(12,5){

\begin{tiny}
              \rput(0,0){\rnode{0}{6}}
              \rput(-3,-3){\rnode{1}{4}}
              \rput(3,-3){\rnode{2}{5}}
              \rput(-3,-6){\rnode{3}{2}}
              \rput(3,-6){\rnode{4}{3}}
              \rput(0,-9){\rnode{5}{1}}   \end{tiny}        }
\rput[l](16,0.5){\begin{tiny}$
\begin{array}{ccl} 
646& = &0\\
 6421& = & 6531\\[12pt]
 464& = & 424\\[5pt]
 242 & = & 212\\
213 & = & 0\\
\end{array}
$\end{tiny}}  
\rput[l](26,0.5){\begin{tiny}$
\begin{array}{ccl}
656& = &0\\
1246&=& 1356\\[12pt]
 565& = & 535\\[5pt]
 353 & = & 313\\
312 & = & 0\\
\end{array}$\end{tiny}}
              \psset{nodesep=1pt,offset=2pt,arrows=<-}
              \ncline{0}{1}
              \ncline{1}{0}
              \ncline{1}{3}
               \ncline{3}{1}
               \ncline{3}{5}
               \ncline{5}{3}
              \psset{nodesep=1pt,offset=2pt,arrows=->}
              \ncline{2}{4}
              \ncline{4}{2}
              \ncline{0}{2}
              \ncline{2}{0}
              \ncline{4}{5}
              \ncline{5}{4}
              
\end{pspicture}
\psset{xunit=2.6mm,yunit=2.6mm,runit=5mm}
\begin{pspicture}(-22,0)(0,0)
\rput(12,5){

\begin{tiny}
              \rput(0,0){\rnode{0}{6}}
              \rput(-3,-3){\rnode{1}{4}}
              \rput(3,-3){\rnode{2}{5}}
              \rput(-3,-6){\rnode{3}{2}}
              \rput(3,-6){\rnode{4}{3}}
              \rput(0,-9){\rnode{5}{1}}   \end{tiny}        }
\rput[l](16,0.5){\begin{tiny}$
\begin{array}{ccl} 646& = &0\\
421& = & 431\\
124 & = & 134\\[12pt]
 464& = & 424\\
 464& = & 434\\
 435& = & 0\\
213 & = & 243\\
\end{array}
$\end{tiny}}  
\rput[l](26,0.5){\begin{tiny}$
\begin{array}{ccl}
 656& = &0\\
 643& = &653\\
 346&=& 356\\
 212&=&242\\[3pt]
 565& = & 535\\
534 & = & 0\\[5pt]
313 & = & 343+353\\
312 & = & 342\\
\end{array}$\end{tiny}}
              \psset{nodesep=1pt,offset=2pt,arrows=<-}
              \ncline{0}{1}
              \ncline{1}{0}
              \ncline{1}{3}
               \ncline{3}{1}
               
              \ncline{1}{4}
               \ncline{4}{1}
               \ncline{3}{5}
               \ncline{5}{3}
              \psset{nodesep=1pt,offset=2pt,arrows=->}
              \ncline{2}{4}
              \ncline{4}{2}
              \ncline{0}{2}
              \ncline{2}{0}
              \ncline{4}{5}
              \ncline{5}{4}
              
\end{pspicture}
\\[14pt]

\section{Ringel-duality on generators-cogenerators of local self-injective algebras}
\begin{small}
Let  $(\A,\ma)$ be  a quasi-hereditary algebra,  then  for any $i\in Q_0(\A)$ there exists a (up to isomorphism) uniquely determined 
indecomposable module   $T(i)\in \ger{F}(\De)\cap \ger{F}(\nabla)$  with  the  following properties: For all $j\in Q_0(\A)$  with   $j\not\ma i$   we have   $[T(i):S(j)]=0$ and $[T(i):S(i)]=\left(T(i):\De(i)\right)=\left(T(i):\nabla(i)\right)=1$, moreover, 
  $\ger{F}(\De)\cap \ger{F}(\nabla)=\add \left(\bigoplus_{i\in Q_0(\A)}T(i)\right)$,    the module $T:=\bigoplus_{i\in Q_0(\A)}T(i)$ is called  \textit{ characteristic   tilting module}.  
The  \textit{Ringel-dual} $R(\A):= \End_{\A}(T)^{op}$ of $\A$  is a basic algebra with $Q_0(R(\A))=Q_0(\A)$  and quasi-hereditary  with  the opposite order $\geqslant$ (we use the notation $\ma_{(R)}$). Moreover,  $R(R(\A))\cong \A$  as quasi-hereditary algebras     (for more details, see \cite{Rin1}).

The functor  $\mathscr{R}:=\Hom_{A}(T, -): \modd A \longrightarrow  \modd R(A)$  yields an exact equivalence between the subcategories  $\ger{F}(\De_{\A})$ and     $\ger{F}(\nabla_{R(\A)})$. Moreover,   $P_{R(\A)}(i)\cong \mathscr{R}(T(i))$, $T_{R(\A)}(i)\cong \mathscr{R}(I(i))$  and   $\De_{R(\A)}(i)\cong \mathscr{R}(\nabla(i))$ for all  $i\in Q_0(\A)$.  
 \end{small}
\\

The class of 1-quasi-hereditary algebras is not closed under Ringel duality.  Example~ 4 in  \cite{P1} presents a 1-quasi-hereditary algebra $A$ for which $R(A)$ is not 1-quasi-hereditary. 
   However,  the  properties of  $R(A)$ for  a 1-quasi-hereditary algebra $A$ considered in ~\cite[Lemma  6.2]{P}   yield   the following lemma.  For  $R(A)$-modules  we will  use  the  index $(R)$, (note that  $1\ma i\ma n$ implies $n\ma_{(R)} i\ma_{(R)} 1 $ for all $i\in \Lambda:=Q_0(R(A))= Q_0(A)$).  

\begin{num}\begin{normalfont}\textbf{Lemma.}\end{normalfont} 
Let $(A,\ma)$ be a 1-quasi-hereditary algebra  with $1\ma i \ma n$ for any $i\in \Lambda$  and let $(R(A),\ma_{(R)} )$ be   the  Ringel dual of $(A,\ma)$. Then
 $\domdim R(A)\geq 2$ and  $P_{(R)}(n)$ is a minimal faithful $R(A)$-module.
\label{33}
    \end{num}

\textit{Proof.}  Since $\left\{1\right\}=\max\left(\Lambda,\ma_{(R)}\right)$, the definition of standard modules implies $P_{(R)}(1)=\De_{(R)}(1)$. 
 According to  \cite[Lemma 6.2 (c) and (b)]{P} we obtain $\soc\left(\De_{(R)}(1)\right) \cong S_{(R)}(n)$ and  $\De_{(R)}(i) \hookrightarrow \De_{(R)}(1)$, since $T(1)\cong S(1)$ (see \cite[5.3]{P}).  Consequently, for all $i\in \Lambda$ we have   $\soc\left(\De_{(R)}(i)\right)\in \add \left(\soc P_{(R)}(n)\right) $,    because    $P_{(R)}(n) \cong I_{(R)}(n)$ (see \cite[6.2(a)]{P}). 

According to  Remark ~\ref{dd2r},    it is 	enough to show $P_{(R)}(i) \hookrightarrow P_{(R)}(n)^{r_i}$ (for some $r_i\in \NN$)  and $P_{(R)}(n)^{r_i}/P_{(R)}(i) \in \ger{F}(\De_{(R)})$  for any $i\in \Lambda$: 
Since $T(i)\in \ger{F}(\De)$, we have $\soc T(i)\in \add\left(\bigoplus_{j\in \Lambda}\soc \De(j)\right)\stackrel{~\ref{def1qh}}{=}\add \left(S(1)\right)$  for any $i\in \Lambda$. Let $\soc T(i)\cong S(1)^{r_i}$, then $T(i)\hookrightarrow T(n)^{r_i}$ since $T(n)\cong I(1)$ (see \cite[5.3]{P}).  
 The exact sequence $\xi:0\to T(i)\to T(n)^{r_i}\to T(n)^{r_i}/T(i)\to 0$ yields $ T(n)^{r_i}/T(i)\in \ger{F}(\nabla)$,   because  $\ger{F}(\nabla)$   is closed under  cokernels of injective maps (see \cite{Rin1}).   By applying   $\mathscr{R}(-)$  to    $\xi$  we obtain an  exact sequence  $0\to P_{(R)}(i)\to P_{(R)}(n)^{r_i}\to P_{(R)}(n)^{r_i}/P_{(R)}(i) \to 0 $ with $ P_{(R)}(n)^{r_i}/P_{(R)}(i)\in \ger{F}(\De_{(R)})$  for any $i\in \Lambda$. 
\hfill $\Box$

\subsection{Transfer of Ringel duality}
Throughout, we keep the notation  of the sets   $\textbf{X}$,  $\textbf{Y}$, $\textbf{X}(1)$,  $\textbf{Y}(1)$  and  the functions $\textbf{X}  \overset{\Phi}{\underset{\Psi}{\rightleftarrows}}   \textbf{Y}$ used  in Section 1. Moreover, 
 we denote  by   $\textbf{X}(R(1))$  the set   of  isomorphism classes of Ringel-duals  of 1-quasi-hereditary algebras. Lemma ~\ref{33} implies  that  $\textbf{X}(R(1))\subseteq \textbf{X}$.  We denote  by   $\textbf{Y}(R(1))$   the image of $\Phi|_{\textbf{X}(R(1))}$. 
 Moreover, let 
$\mathcal{X}:=\textbf{X}(1)\cup \textbf{X}(R(1))$   and          $\mathcal{Y}:=\textbf{Y}(1)\cup \textbf{Y}(R(1))$  as well as          $\widehat{\mathcal{X}}:=\textbf{X}(1)\cap\textbf{X}(R(1))$    and          $\widehat{\mathcal{Y}}:=\textbf{Y}(1)\cap\textbf{Y}(R(1))$. 

The map $\mathcal{R}: \mathcal{X}\to \mathcal{X}$   with    $\mathcal{R}\left(\left[\A\right]\right)= \left[R(\A)\right]$  is obviously  bijective and $\mathcal{R}^2=\id_{\mathcal{X}}$.
The Morita-Tachikawa Theorem ~\ref{Morita} and  Theorem A  yield the  transfer  of  Ringel-duality  for  $\mathcal{X}$ on  $\mathcal{Y}$ (illustrated on the  picture below).

\newpage
\psset{xunit=0.75cm,yunit=0.75cm,runit=1cm}
\begin{pspicture}(-12.7,2.4)(0,0)
\rput(-0.7,-1.2){
\rput(-8,3.2){\rnode{cenr}{Let  $
\textbf{R}:=\Phi|_{\mathcal{X}}\circ \mathcal{R}|_{\mathcal{X}}\circ \Psi|_{\mathcal{Y}}
$, then }}
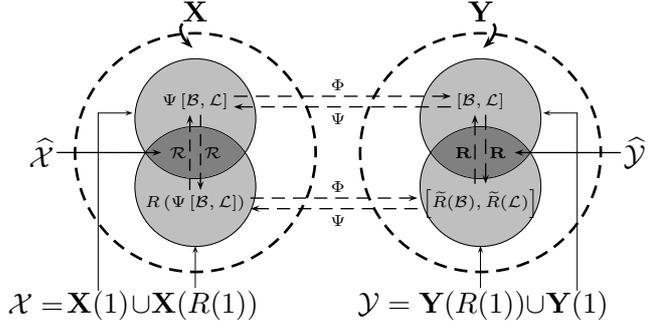
\rput(-8.5,2){\rnode{cen}{$
\textbf{R}: \mathcal{Y}  \longrightarrow  \mathcal{Y} 
$}}
\rput(-8.5,1){\rnode{cen0}{$
\textbf{R}( [\B,\LL]) : =  \left[\widetilde{R}(\B),\widetilde{R}(\LL)\right]
$}}
\rput(-8.2,0){\rnode{cen1}{and}}
\rput(-8.5,-1.2){\rnode{cen2}{$\textbf{R}\left(\textbf{R}[\B,\LL]\right)= [\B,\LL]$}}
}
\pscircle[linewidth=1pt,linestyle=dashed](0,0){1.6cm}
\pscircle[linewidth=1pt,linestyle=dashed](5,0){1.6cm}
\pscircle*[linewidth=1pt,linecolor=gray!50](0,0.6){0.8cm}
\pscircle*[linewidth=1pt,linecolor=gray!50](5,0.6){0.8cm}
\pscircle*[linewidth=1pt,linecolor=gray!50](0,-0.6){0.8cm}
\pscircle[linewidth=0.3pt](0,-0.6){0.8cm}
\pscircle*[linewidth=1pt,linecolor=gray!50](5,-0.6){0.8cm}
\pscircle[linewidth=0.3pt](5,-0.6){0.8cm}
\pscircle[linewidth=0.3pt](0,0.6){0.8cm}
\pscircle[linewidth=0.3pt](5,0.6){0.8cm}

\psclip{
\pscircle[linestyle=none,linecolor=black](0,0.6){0.8cm}
\pscircle[linestyle=none,linecolor=black](0,-0.6){0.8cm}
        }
        \psframe*[linecolor=gray](-1,-1)(5,5)
\endpsclip

\psclip{
\pscircle[linestyle=none,linecolor=black](5,0.6){0.8cm}
\pscircle[linestyle=none,linecolor=black](5,-0.6){0.8cm}
        }
        \psframe*[linecolor=gray](4,-1)(7,5)
\endpsclip

\begin{tiny} 
\rput(2.5,1.2){$\Phi$}
\rput(2.5,0.6){$\Psi$}
\rput(2.5,-1.2){$\Psi$}
\rput(2.5,-0.6){$\Phi$}
\rput(0.3,0){$\mathcal{R}$}
\rput(-0.3,0){$\mathcal{R}$}
\rput(5.3,0){$\textbf{R}$}
\rput(4.7,0){$\textbf{R}$}
\rput(5,0.9){\rnode{BM}{$\left[\B,\LL\right]$}}
\rput(0,0.9){\rnode{PsiBM}{$\Psi \left[\B,\LL\right]$}}
\rput(0,-0.9){\rnode{RPsiBM}{$R\left(\Psi \left[\B,\LL\right]\right)$}}
\rput(5,-0.9){\rnode{BMv}{$\left[\widetilde{R}(\B),\widetilde{R}(\LL)\right]$}}
\end{tiny}

\begin{small}
\rput(-2.7,0){\rnode{vc}{$\widehat{\mathcal{X}}$}}
\rput(-0.55,0){\rnode{vc1}{}}
\rput(7.7,0){\rnode{wc}{$\widehat{\mathcal{Y}}$}}
\rput(5.55,0){\rnode{wc1}{}}
\rput(0,2.5){\rnode{x}{$\textbf{X}$}}
\rput(0,1.8){\rnode{y}{}}
\rput(5,2.5){\rnode{xx}{$\textbf{Y}$}}
\rput(5,1.8){\rnode{yy}{}}
\rput(-2.8,-2.68){\rnode{a10}{$\mathcal{X}=$}}
\rput(-1.7,-2.7){\rnode{a1}{$\textbf{X}(1)$}}
\rput(-1.1,0.7){\rnode{a2}{}}
\rput(6.7,-2.7){\rnode{b1}{$\textbf{Y}(1)$}}
\rput(6.1,0.7){\rnode{b2}{}}
\rput(0,-2.7){\rnode{a1r}{$\cup \textbf{X}(R(1))$}}
\rput(0,-1.6){\rnode{a2r}{}}
\rput(5,-2.7){\rnode{b1r}{$\mathcal{Y}=\textbf{Y}(R(1))\cup \hspace{7.5mm}$}}
\rput(5,-1.6){\rnode{b2r}{}}
\end{small}

\ncangle[angleB=90,linewidth=0.27pt]{<-}{a2}{a1}
\ncangle[angleB=90,linewidth=0.27pt]{<-}{b2}{b1}

\psset{nodesep=1pt}
\ncline[linewidth=1pt]{->}{T1}{T2}
\ncline[linewidth=1pt]{->}{RT1}{RT2}
\ncarc[linewidth=1pt,arcangle=-45]{->}{x}{y}
\ncarc[linewidth=1pt,arcangle=45]{->}{xx}{yy}
\ncline[linewidth=0.5pt]{->}{vc}{vc1}
\ncline[linewidth=0.5pt]{->}{wc}{wc1}
\ncline[linewidth=0.27pt]{->}{a1r}{a2r}
\ncline[linewidth=0.27pt]{<-}{b2r}{b1r}
\psset{nodesep=2pt,offset=2pt}
\ncline[linestyle=dashed,linewidth=0.5pt]{->}{BM}{PsiBM}
\ncline[linestyle=dashed,linewidth=0.5pt]{->}{PsiBM}{BM}
\ncline[linestyle=dashed,linewidth=0.5pt]{->}{BMv}{RPsiBM}
\ncline[linestyle=dashed,linewidth=0.5pt]{->}{RPsiBM}{BMv}
\ncline[linestyle=dashed,linewidth=0.5pt]{->}{BM}{BMv}
\ncline[linestyle=dashed,linewidth=0.5pt]{->}{BMv}{BM}
\ncline[linestyle=dashed,linewidth=0.5pt]{->}{PsiBM}{RPsiBM}
\ncline[linestyle=dashed,linewidth=0.5pt]{->}{RPsiBM}{PsiBM}
\end{pspicture}
\\[2.7cm]
Obviously,  for any pair $[\B,\LL]\in \mathcal{Y} $ 
 the algebra $\A=\Psi\left([\B,\LL]\right)$    is 1-quasi-hereditary or $\A$ is the  Ringel dual of a 1-quasi-hereditary algebra  $R(\A)=\Psi\left([\widetilde{R}[\B],\widetilde{R}(\LL)]\right)$. The minimal faithful $\A$-module  $\LL$
is the  projective indecomposable $P(\textsl{1}_{\A})$  which   corresponds  to the minimal vertex $\textsl{1}_{\A}$ 
(see ~\ref{mini} and ~\ref{33}).  In particular, $\textsl{1}_{R(\A)}$ is the maximal vertex in $\Lambda$  with respect to the partial order corresponding to  $\A$.  The direct summands of  the $\B$-module $\LL$ are 
 isomorphic to  $\Hom_{\A}\left(P(i),P(\textsl{1}_{\A})\right)$ for any  $i\in \Lambda$.

\begin{subnum}\begin{normalfont}\textbf{Lemma.}\end{normalfont}  Let  $[\B,\LL]\in  \mathcal{Y}$   and let  $\A$ be a quasi-hereditary algebra  with $[\A]=\Psi\left([\B,\LL]\right)$. 
Then for $
\textbf{R}( [\B,\LL]) : =  \left[\widetilde{R}(\B),\widetilde{R}(\LL)\right]
$ the following hold: 
\begin{center}
$\B\cong \widetilde{R}(\B)$ \  and \  $ \displaystyle \widetilde{R}(\LL) \cong \bigoplus_{i\in\Lambda}\Hom_{\A}(T(i),T(\textsl{1}_{R(\A)})).$
\end{center}
In particular,  there exists a 1-quasi-hereditary algebra $A$   with $\B\cong \End_A(P(1_A))^{op}$  and    
$\displaystyle\LL\cong \bigoplus_{i\in \Lambda}\Hom_{A}(P(i),P(1_A))$      or      $\displaystyle \LL\cong \bigoplus_{i\in \Lambda}\Hom_{A}(T(i),T(n_A)).$
\label{Ringel}
\end{subnum}

\textit{Proof.} 
 It follows from  Theorem ~\ref{Morita} (with the preceding notation) that $\B\cong \End_{\A}(P(\textsl{1}_{\A}))^{op}$ and $\widetilde{R}(\B)\cong \End_{R(\A)}(P(\textsl{1}_{R(\A)}))^{op}$ as well as $\LL \cong \bigoplus_{i\in\Lambda}\Hom_{\A}\left(P(i),P(\textsl{1}_{\A})\right)$ and $  \widetilde{R}(\LL)\cong \bigoplus_{i\in \Lambda}\Hom_{R(\A)}(P_{(R)}(i),P_{(R)}(\textsl{1}_{R(\A)})) $.
      Since $R(\A)\cong \End_{\A}(T)^{op}$ and $P_{(R)}(i)\cong \Hom_{\A}(T,T(i))$, the   functor  $\Hom_{\A}(T,-):\modd \A\rightarrow \modd R(\A)$   yields  an isomorphism 
\begin{center}
$ \Hom_{R(\A)}(P_{(R)}(i),P_{(R)}(\textsl{1}_{R(\A)}))\cong \Hom_{\A}(T(i),T(\textsl{1}_{R(\A)}))$ 
\end{center}
 for all $i\in \Lambda$ (see \cite[2.1]{ARS}). In particular,    $\B\cong \End_{\A}(P(\textsl{1}_{\A}))^{op}\cong \End_{R(\A)}(P(\textsl{1}_{(R)}))^{op}\cong \widetilde{R}(\B)$  because  $P(\textsl{1}_{\A}) \cong T(\textsl{1}_{R(\A)})$ (see  \cite[Remark 5.3 and Lemma 6.2 (a)]{P}).
 
  If the algebra $\A$  is 1-quasi-hereditary, then $\LL \cong \bigoplus_{i\in\Lambda}\Hom_{\A}\left(P(i),P(\textsl{1}_{\A})\right)$.  If not,   then   $A=R(\A)$ is 1-quasi-hereditary with $1_{A}\ma i\ma n_{A}$ for all $i\in Q_0(A) $. For  $\Phi([A])=[B,L] $ we obtain   $[\B,\LL]= \textbf{R}([B,L] )=[B,\widetilde{R}(L)] $. Therefore  $\B\cong \End_{A}\left(P(1_{A}\right)^{op}$ and $\LL\cong   \bigoplus_{i\in \Lambda}\Hom_{A}(T(i),T(n_{A}))$.
  \hfill $\Box$ 
\\

Lemma ~\ref{Ringel} implies that for any $[\B,\LL]\in \mathcal{Y}$ we have $\textbf{R}([\B,\LL])=[\B,\widetilde{R}(\LL)]$. In particular, $\widetilde{R}(\LL)$ is a multiplicity-free generator-cogenerator of $\B$ (see ~\ref{Morita}).     The  Ringel-duality on $\mathcal{Y}$ conforms  with the duality on a subclass of multiplicity-free  generator-cogenerators  of local self-injective algebras, which arises  from  1-quasi-hereditary algebras (via $\B\cong \End_A(P(1_A))^{op}$).
\\[-3mm]

 According to  Theorem  A  the first component $B$ of a pair in  $\mathcal{Y}$  is a local self-injective  algebra   having a module  satisfying the properties in Definition~\ref{property_ma} for some partial order. We denote by  $\ger{L}(B)$ the set of all $B$-modules  $M$   such that     $(B,M)$ or $(B,\widetilde{R}(M))$ satisfies the condition  \fbox{$\ma$} for  some partial order $\ma$ on    $\Lambda=\left\{1,\ldots , \dimm_KB\right\}$.  Obviously, there exist  finitely many  partial orders  on $\Lambda$.   However,  there are examples of a  partial order  $\ma$ on $\Lambda$  with infinitely  many  pairwise non-isomorphic  $B$-modules $M$   such that   $(B,M)$ satisfies the property  \fbox{$\ma$}  (in  the next   example~\ref{non_Ringel_dual} the $B$-module $L_2$ depends  on the choice of $\mu\in K$,  there we have $\left\langle X^2+ \mu Y\right\rangle \not\cong \left\langle X^2+ \mu' Y\right\rangle$  if $\mu\neq \mu'$).  Moreover, if $(B,M)$  and $(B,M')$  satisfy the property  \fbox{$\ma$}, then $[B,M]\in \widehat{\mathcal{Y}}$ does not imply $[B,M']\in \widehat{\mathcal{Y}}$.

Let  $\textbf{L}$ be the set of isomorphism classes of  algebras $B$ with $\ger{L}(B)$$ \neq \emptyset$  (the function  $\textbf{X}(1) \to \textbf{L}$ with  $[A]\mapsto \left[\End_A(P(1))\right]$  is     surjective,  non-injective and the set $\textbf{L}$ is not finite). 
  For any $[B]\in \textbf{L}$    we define
 $\mathcal{Y}(B):=\left\{[B,M]\mid M\in \ger{L}(B)\right\}$  \     and \     $\mathcal{X}(B):=\left\{\left[\End_B(M)^{op}\right]\mid M\in \ger{L}(B)\right\}$.
  It is easy to see that for all  $\left[B\right],\left[B'\right]\in \textbf{L}$ with $\left[B\right]\neq  \left[B'\right]$ we have $\mathcal{Y}(B) \cap \mathcal{Y}(B')=\emptyset$ and a pair $\left[B,M\right]\in  \mathcal{Y}$ belongs to $\mathcal{Y}(B)$.    
  This implies  $\displaystyle \mathcal{Y}=\overset{\cdot}{\bigcup_{[B]\in \textbf{L}}}\mathcal{Y}(B)$ \   and \  
$\displaystyle \mathcal{X}=\overset{\cdot}{\bigcup_{[B]\in \textbf{L}}}\mathcal{X}(B)$.
    \\
  In the  picture   the sets  $\mathcal{X}$  and $\mathcal{Y}$ are  presented as the  disjoint union of $\mathcal{X}(B)$ and $\mathcal{Y}(B)$  (symbolized by the circles, they are closed under the  Ringel-duality $\textbf{R}$ and  $\mathcal{R}$)  respectively. The dark circles inside the circle corresponding to  $\mathcal{X}(B)$  symbolizes   $\mathcal{X}(B)\cap\widehat{\mathcal{X}}$.  Similarly,  a pair $[B,M]$ in the   dark circle of $\mathcal{Y}(B)$  has the property \fbox{$\ma$}  and $\textbf{R}([B,M])$ has the property \fbox{$\geqslant$}.    
  In particular, $\widehat{\mathcal{X}}$ and $\widehat{\mathcal{Y}}$ are  the disjoint 	  unions  of the \textit{dark circles}. They  are also closed under $\mathcal{R}$ resp.  $\textbf{R}$. The observation of Ringel-duality on $\mathcal{X}$ and  $\mathcal{Y}$ can by restricted to  $\mathcal{X}(B)$ and  $\mathcal{Y}(B)$ respectively.  The dark circles inside  the circle
  \\
\parbox{10cm}{ corresponding to  $\mathcal{X}(B)$  symbolizes   $\mathcal{X}(B)\cap\widehat{\mathcal{X}}$.  Similarly,  a pair $[B,M]$ in the   dark circle of $\mathcal{Y}(B)$  has the property \fbox{$\ma$}  and $\textbf{R}([B,M])$ has the property \fbox{$\geqslant$}, i.e., $[B,M]\in \widehat{\mathcal{Y}}$.   In particular, $\widehat{\mathcal{X}}$ and $\widehat{\mathcal{Y}}$ are  the disjoint 	  unions  of the \textit{dark circles}. They  are also closed under $\mathcal{R}$ resp.  $\textbf{R}$. The observation of Ringel-duality on $\mathcal{X}$ and  $\mathcal{Y}$ can by restricted to  $\mathcal{X}(B)$ and  $\mathcal{Y}(B)$ respectively.}\\[-2mm]

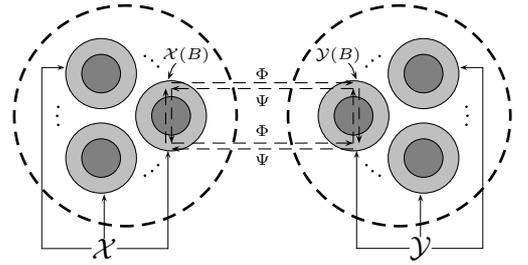
\psset{xunit=0.8cm,yunit=0.8cm,runit=1cm}
\begin{pspicture}(-13,0)(0,0)
\rput(1,3.5){\pscircle[linewidth=1pt,linestyle=dashed](0,0){1.48cm}
\rput(0.3,1){\rnode{G}{}}
\rput(0.55,0.8){\rnode{H}{}}
\rput(0.3,-1){\rnode{Gc}{}}
\rput(0.55,-0.8){\rnode{Hc}{}}
\rput(-1.1,0.17){\rnode{Gl}{}}
\rput(-1.1,-0.17){\rnode{Hl}{}}
\ncarc[linewidth=1pt,linestyle=dotted, arcangle=15]{G}{H}
\ncarc[linewidth=1pt,linestyle=dotted, arcangle=-15]{Gc}{Hc}
\ncarc[linewidth=1pt,linestyle=dotted, arcangle=-15]{Gl}{Hl}
\pscircle*[linewidth=1pt,linecolor=gray!50](0.75,0){0.47cm}
\pscircle*[linewidth=1pt,linecolor=gray!50](-0.4,0.7){0.47cm}
\pscircle*[linewidth=1pt,linecolor=gray!50](-0.4,-0.7){0.47cm}
\pscircle*[linewidth=1pt,linecolor=gray](0.75,0){0.26cm}
\pscircle*[linewidth=1pt,linecolor=gray](-0.4,0.7){0.26cm}
\pscircle*[linewidth=1pt,linecolor=gray](-0.4,-0.7){0.26cm}
\pscircle[linewidth=0.3pt](0.75,0){0.26cm}
\pscircle[linewidth=0.3pt](-0.4,0.7){0.26cm}
\pscircle[linewidth=0.3pt](-0.4,-0.7){0.26cm}
\pscircle[linewidth=0.3pt](0.75,0){0.47cm}
\pscircle[linewidth=0.3pt](-0.4,0.7){0.47cm}
\pscircle[linewidth=0.3pt](-0.4,-0.7){0.47cm}
\rput(4.5,0){
\pscircle[linewidth=1pt,linestyle=dashed](0,0){1.48cm}
\rput(-0.3,-1){\rnode{G}{}}
\rput(-0.55,-0.8){\rnode{H}{}}
\rput(-0.3,1){\rnode{Gc}{}}
\rput(-0.55,0.8){\rnode{Hc}{}}
\rput(1.1,-0.17){\rnode{Gl}{}}
\rput(1.1,0.17){\rnode{Hl}{}}
\ncarc[linewidth=1pt,linestyle=dotted, arcangle=15]{G}{H}
\ncarc[linewidth=1pt,linestyle=dotted, arcangle=-15]{Gc}{Hc}
\ncarc[linewidth=1pt,linestyle=dotted, arcangle=-15]{Gl}{Hl}
\pscircle*[linewidth=1pt,linecolor=gray!50](-0.75,0){0.47cm}
\pscircle*[linewidth=1pt,linecolor=gray!50](0.4,0.7){0.47cm}
\pscircle*[linewidth=1pt,linecolor=gray!50](0.4,-0.7){0.47cm}
\pscircle*[linewidth=1pt,linecolor=gray](-0.75,0){0.26cm}
\pscircle*[linewidth=1pt,linecolor=gray](0.4,0.7){0.26cm}
\pscircle*[linewidth=1pt,linecolor=gray](0.4,-0.7){0.26cm}
\pscircle[linewidth=0.3pt](-0.75,0){0.26cm}
\pscircle[linewidth=0.3pt](0.4,0.7){0.26cm}
\pscircle[linewidth=0.3pt](0.4,-0.7){0.26cm}
\pscircle[linewidth=0.3pt](-0.75,0){0.47cm}
\pscircle[linewidth=0.3pt](0.4,0.7){0.47cm}
\pscircle[linewidth=0.3pt](0.4,-0.7){0.47cm}
}
\begin{tiny} \rput(2.25,0.7){$\Phi$}
\rput(2.25,0.25){$\Psi$}
\rput(2.25,-0.25){$\Phi$}
\rput(2.25,-0.73){$\Psi$}
\rput(1,1){\rnode{BM}{$\mathcal{X}(B)$}}
\rput(3.5,1){\rnode{cM}{$\mathcal{Y}(B)$}}
\rput(0.7,0.6){\rnode{BMf}{}}
\rput(3.8,0.6){\rnode{cMf}{}}
\ncarc[linewidth=0.4pt,arcangle=-15]{->}{BM}{BMf}
\ncarc[linewidth=0.4pt,arcangle=15]{->}{cM}{cMf}
\end{tiny}
\begin{small}
\rput(0,1.6){\rnode{y}{}}
\rput(4.5,1.6){\rnode{yy}{}}
\rput(-0.35,-2.2){\rnode{A}{$\mathcal{X}$}}
\rput(-1,0.8){\rnode{B}{}}
\rput(-0.35,-1.3){\rnode{C}{}}
\rput(0.7,-0.58){\rnode{D}{}}
\ncangles[angleA=180,angleB=-180,arm=0.3,linewidth=0.27pt]{<-}{B}{A}
\ncangle[angleA=90,linewidth=0.27pt]{<-}{D}{A}
\ncline[linewidth=0.27pt]{<-}{C}{A}
\rput(4.5,0){
\rput(0.35,-2.2){\rnode{A}{$\mathcal{Y}$}}
\rput(1,0.8){\rnode{B}{}}
\rput(0.35,-1.3){\rnode{Z}{}}
\ncangles[angleA=0,angleB=0,arm=0.3,linewidth=0.27pt]{<-}{B}{A}
\rput(-0.7,-0.58){\rnode{B}{}}
\ncangle[angleA=-180,arm=-0,linewidth=0.27pt]{->}{A}{B}
\ncline[linewidth=0.27pt]{<-}{Z}{A}
}
\end{small}

\rput(0.71,0.5){\rnode{p1}{}}
\rput(0.71,-0.5){\rnode{p2}{}}
\put(3.6,0){
\rput(-0.71,0.5){\rnode{q1}{}}
\rput(-0.71,-0.5){\rnode{q2}{}}}
\psset{nodesep=1pt, offset=1.1pt}
\ncline[linewidth=0.4pt,linestyle=dashed]{->}{q1}{p1}
\ncline[linewidth=0.4pt,linestyle=dashed]{->}{p1}{q1}
\ncline[linewidth=0.4pt,linestyle=dashed]{->}{q2}{p2}
\ncline[linewidth=0.4pt,linestyle=dashed]{->}{p2}{q2}
\ncline[linewidth=0.4pt,linestyle=dashed]{->}{p1}{p2}
\ncline[linewidth=0.4pt,linestyle=dashed]{->}{p2}{p1}

\ncline[linewidth=0.4pt,linestyle=dashed]{->}{q1}{q2}
\ncline[linewidth=0.4pt,linestyle=dashed]{->}{q2}{q1}

\psset{nodesep=1pt}
\ncline[linewidth=1pt]{->}{T1}{T2}
\ncline[linewidth=1pt]{->}{RT1}{RT2}
\ncarc[linewidth=1pt,arcangle=-45]{->}{x}{y}
\ncarc[linewidth=1pt,arcangle=45]{->}{xx}{yy}
\ncline[linewidth=0.6pt]{<-}{l2}{l1}
\ncline[linewidth=0.6pt]{->}{s1}{s2}
\ncline[linewidth=0.6pt]{<-}{d2}{d1}
\ncline[linewidth=1pt]{<-}{b2}{b1}
\ncline[linewidth=0.27pt]{->}{a1r}{a2r}
\ncline[linewidth=0.27pt]{<-}{b2r}{b1r}
}
\end{pspicture}
\\[-7mm]

\begin{subnum}\begin{normalfont}\textbf{Example.} The algebra  $B=K[x,y]/\left\langle xy,\ x^4-y^2\right\rangle$ is local and self-injective with  $\dimm_KB=6$. The pair $(B,L_i)$ satisfies the property \fbox{$\ma_{(i)}$},  where   $L_i=\bigoplus_{k=1}^6L_i(k)$  and $\ma_{(i)}$ for $i=1,2$   are  presented in the following diagrams  in  the same way  as in  Example ~\ref{exampleB}. Both pairs belong  to $\mathcal{Y}(B)$, however,   $[B,L_1]\in \widehat{\mathcal{Y}}$ and  $[B,L_2]\not\in \widehat{\mathcal{Y}}$. 
\end{normalfont}
\label{non_Ringel_dual}
\end{subnum}

 \hspace*{8mm}  \psset{xunit=0.75mm,yunit=0.9mm,runit=1.2mm}
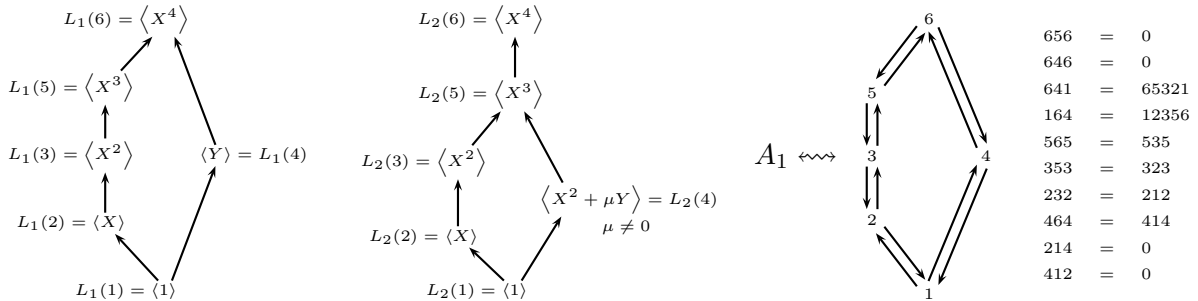
\begin{pspicture}(-10,-17)(0,2)
\begin{tiny}
              \rput(1,0){\rnode{1}{$L_1(6)=\left\langle X^4\right\rangle \hspace{11mm}$}}
              \rput(-10,-10){\rnode{2}{$L_1(5)=\left\langle X^3\right\rangle \hspace{9mm}$}}
              \rput(10,-20){\rnode{3}{$\hspace{9mm}\left\langle Y\right\rangle = L_1(4)$}}
              \rput(-10,-20){\rnode{4}{$L_1(3)= \left\langle X^2\right\rangle \hspace{9mm}$}}
              \rput(-10,-30){\rnode{5}{$ L_1(2)= \left\langle  X\right\rangle \hspace{9mm}$}}
              \rput(1,-40){\rnode{6}{$L_1(1)=\left\langle 1\right\rangle \hspace{11mm} $}} \end{tiny}
              \psset{nodesep=1pt}
              \ncline{<-}{1}{2}
              \ncline{<-}{1}{3}
              \ncline{<-}{3}{6}
              \ncline{<-}{2}{4}
              \ncline{<-}{4}{5}
              \ncline{<-}{5}{6}
\end{pspicture} 
\psset{xunit=0.75mm,yunit=0.9mm,runit=1.2mm}
\begin{pspicture}(-60,-17)(0,0)
\begin{tiny}
              \rput(0,0){\rnode{1}{$L_2(6)=\left\langle X^4\right\rangle \hspace{9mm}$}}
              \rput(0,-11){\rnode{2}{$L_2(5)=\left\langle X^3\right\rangle \hspace{9mm}$}}
              \rput(10,-26.5){\rnode{3}{$\hspace{15mm}\left\langle X^2 + \mu  Y\right\rangle = L_2(4)$}}
              \rput(19.5,-30.5){\rnode{3mu}{$\mu\neq 0$}}
              \rput(-10,-21){\rnode{4}{$L_2(3)= \left\langle X^2\right\rangle \hspace{9mm}$}}
              \rput(-10,-32){\rnode{5}{$ L_2(2)= \left\langle  X\right\rangle \hspace{9mm}$}}
              \rput(0,-40){\rnode{6}{$L_2(1)=\left\langle 1\right\rangle \hspace{10mm} $}} \end{tiny}
              \psset{nodesep=1pt}
              \ncline{<-}{1}{2}
              \ncline{<-}{2}{3}
              \ncline{<-}{3}{6}
              \ncline{<-}{2}{4}
              \ncline{<-}{4}{5}
              \ncline{<-}{5}{6}
\end{pspicture} 
\psset{xunit=2.5mm,yunit=2.8mm,runit=5mm}
\begin{pspicture}(-5,0)(0,0)
\rput(16.5,5){
\begin{small}\rput(-7,-6){$A_1  \leftrightsquigarrow$} \end{small}
\begin{tiny}
              \rput(0,0.5){\rnode{0}{6}}
              \rput(-3,-3){\rnode{1}{5}}
              \rput(3,-6){\rnode{2}{4}}
              \rput(-3,-6){\rnode{3}{3}}
              \rput(-3,-9){\rnode{4}{2}}
              \rput(0,-12.5){\rnode{5}{1}}   \end{tiny}        }
\rput[l](21,-1){
\begin{tiny}$
\begin{array}{ccl} 
 656& = &0\\[3pt]
 646& = &0\\[3pt]
 641& = & 65321\\[3pt]
 164&=& 12356\\[3pt]
565& = &535\\[3pt]
353& = &323\\[3pt]
232& = &212\\[3pt]
464& = &414\\[3pt]
214& = &0\\[3pt]
412& = &0\\
\end{array}
$\end{tiny}}  
              \psset{nodesep=2pt,offset=2pt,arrows=<-}
              \ncline{0}{1}
              \ncline{1}{0}
              \ncline{1}{3}
              \ncline{3}{1}
              \ncline{3}{4}
              \ncline{4}{3}
              \psset{nodesep=2pt,offset=2pt,arrows=->}
              \ncline{2}{5}
              \ncline{5}{2}
              \ncline{0}{2}
              \ncline{2}{0}
              \ncline{4}{5}
              \ncline{5}{4}
              
\end{pspicture}
\\[2.4cm]
On the right-hand side we present the quiver and relations of the 1-quasi-hereditary algebra $A_1=\End_B(L_1)^{op}$. The algebras $A_1$ and $R(A_1)$ are isomorphic as quasi-hereditary algebras  (i.e.,  $A_1$  is Ringel self-dual), because $[B,L_1]=[B,\widetilde{R}(L_1)]$ (see Theorem ~\ref{ringel-dual_Y}).  The quiver and relations of $A_2=\End_B(L_2)^{op}$ can be found in \cite[Example 4]{P1} (there $q=1+\mu^2$). The Ringel-dual  $\left[B,\widetilde{R}(L_2)\right]$   of $[B,L_2]$ is given by $\widetilde{R}(L_2) =\bigoplus_{i=1}^6\widetilde{R}(L_2(i))$  with 
\begin{center}
$\begin{array}{lll}
	\widetilde{R}(L_2(1))=\left\langle X^4\right\rangle, \    &  \widetilde{R}(L_2(2))=\left\langle X^3 \right\rangle, &      \widetilde{R}(L_2(3))=\left\langle \mu X^2-Y\right\rangle, \\ \widetilde{R}(L_2(4))=\left\langle Y \right\rangle, & \widetilde{R}(L_2(5))=\left\langle (X,Y) \right\rangle +\left\langle (Y,0)\right\rangle \subset B\oplus B, \   & \widetilde{R}(L_2(6))=B. 
\end{array}
$
\end{center}

\subsection{Ringel-duality on $\widehat{\mathcal{Y}}$}
\begin{small}In \cite[Sec. 5 and 6]{P} 1-quasi-hereditary  algebras whose  isomorphism classes belong to   $\widehat{\mathcal{X}}$ have been considered.  These  results imply a precise  description of the Ringel-duality $\textbf{R}$ on $\widehat{\mathcal{Y}}$. \end{small}
\\

Let $[B,L]$ be in $\textbf{Y}(1)$,  then  $B\cong \End_A(P(1))^{op}$ and $L\cong \bigoplus_{i\in \Lambda}B\circ \ger{f}_{(i)}$, where $A=\End_B(L)^{op}$ is a 1-quasi-hereditary algebra  with $(\Lambda,\ma)$ (here $1\ma i\ma n$ for all $i\in \Lambda$) and $\ger{f}_{(i)}:=\ger{f}_{(1,i,1)}$ is the endomorphism of $P(1)$ corresponding to the path $p(1,i,1)$  of $A$ (see ~\ref{End(P(1))}).  Let $R(A)=\End_A(T)^{op}$ with $(\Lambda,\ma_{(R)})$  also  be  1-quasi-hereditary.  According to  Theorem 6.1 in \cite{P},  the  direct summand   $T(i)$ of the characteristic tilting $A$-module $T$   is a submodule      and a factor module of $P(1)\cong I(1)$  for any  $i\in \Lambda$ or  more precisely $T(i)\cong P(1)/\left(\sum_{l\in \Lambda\backslash \Lambda_{(i)}}P(l)\right)\cong \bigcap_{l\in \Lambda\backslash \Lambda_{(i)}}\ker\left(P(1)\sur I(l)\right)$.
   Consequently,    the subspace of $P(1)$ corresponding to the vertex $1$ contains an element   $\ger{t}(i)$ which generates    $T(i)$.
 For any $i\in \Lambda $    we denote  by $\mathcal{T}_{(i)}$ the following endomorphism of $P(1)$:
\begin{center}
$\mathcal{T}_{(i)}=(\iota(i)\circ \pi(i)):\left(P(1)\stackrel{\pi(i)}{\sur} T(i)\stackrel{\iota(i)}{\hookrightarrow} P(1)\right)$ \  with  \  $ e_1\stackrel{\pi(i)}{\mapsto}  \ger{t}(i)$  \  and \  $ \ger{t}(i)\stackrel{\iota(i)}{\mapsto}  \ger{t}(i) $
\end{center}
 Using  Lemma 3.2 \cite{P1} it is easy to show, that the pair $(B,\LL)$ with $\LL\cong \bigoplus_{i\in \Lambda}B\circ \mathcal{T}_{(i)}$ satisfies  the
    condition   \fbox{$\ma_{(R)}$}. We recall that  using   the notations of   Section  2,    for all $i,j\in \Lambda$  with   $i\not\ma j$  (and therefore  $j\not\ma_{(R)} i$)   we have $j\in \Lambda\backslash \Lambda^{(i)}$ resp. $i\in \Lambda\backslash \Lambda_{(j)}$.

     Ringel-duality $\textbf{R}$ on $\widehat{\mathcal{Y}}$ provides a relationship between  the endomorphisms $\mathcal{T}_{(1)}, \ldots , \mathcal{T}_{(n)}$  and  $\ger{f}_{(1)},\ldots , \ger{f}_{(n)}$  of $P(1)$.  (The following statement yields  Theorem C from the introduction.)

  \begin{subnum}
\begin{normalfont}\textbf{Theorem.}\end{normalfont} 
\textit{Let $(A,\ma)$ be  a 1-quasi-hereditary algebra  and let  $[(B,L),\ma]$ be the corresponding pair in  $ \textbf{Y}(1)$  with  $B=\End_A(P(1))^{op}$  and  $\displaystyle L\cong \bigoplus_{i\in \Lambda}L(i)$   where   $L(i)= B\circ \ger{f}_{(i)}$  for any $i\in \Lambda$.    Let  $\textbf{R}([B,L])=[B,\widetilde{R}(L)]$  with  $\displaystyle \widetilde{R}(L)\cong  \bigoplus_{i\in \Lambda} \widetilde{R}(L(i))   $.  Then the  following statements are equivalent:
\begin{itemize}
	\item[(i)]  $[B,L]\in \widehat{\mathcal{Y}}$. 
\item[(ii)]    $\displaystyle  \widetilde{R}(L(i))\cong B \circ \mathcal{T}_{(i)} $, where $\mathcal{T}_{(i)} \in B$   satisfies $\im\mathcal{T}_{(i)}=T(i)$ for every  $i\in \Lambda$. 
\item[(iii)]    $\displaystyle  \widetilde{R}(L(i))\cong B/\left(\sum_{j\in \Lambda\backslash \Lambda_{(i)}}L(j)\right) \cong \bigcap_{j\in \Lambda\backslash \Lambda_{(i)}} \ker\left(B\sur L(j)\right)$ \    for every  $i\in \Lambda$.
\end{itemize}}
\label{ringel-dual_Y}
\end{subnum}

\textit{Proof.} Let $A$ be a 1-quasi-hereditary algebra with $[A]=\Psi([B,L])$.  According to ~\ref{Ringel} for any $i\in \Lambda$  we have  $\widetilde{R}(L(i))\cong \Hom_{A}(T(i),T(n))$.

 $(i)\Rightarrow (ii)$ The assumption  $[B,L]\in \widehat{\mathcal{Y}}$ implies $[A]\in \widehat{\mathcal{X}}$.   As already explained,  we can define endomorphisms $\mathcal{T}_{(i)}=(\iota(i)\circ \pi(i))$ for any $i\in \Lambda$.  Since  $\Hom_A(-,P(1))$ is  exact,   
  the inclusion $T(i)\stackrel{\iota(i)}{\hookrightarrow} P(1)$  yields   a surjective  $B$-map   $\Hom_{A}(P(1),P(1))\sur \Hom_A(T(i),P(1))$ with  $F\mapsto F\circ \iota(i)$  and therefore   $\Hom_A(T(i),P(1))= B\circ \iota_{(i)}$ for all $i\in \Lambda$. The map  $B\circ \iota(i)  \rightarrow B\circ  \mathcal{T}_{(i)}$ given by  $F\circ \iota(i) \mapsto  F\circ  \iota(i) \circ \pi(i) = F\circ \mathcal{T}_{(i)}$ for all $F\in B$ is obviously  a $B$-module isomorphism.   We obtain  $\widetilde{R}(L(i))\cong B\circ \iota(i)  \cong  B\circ  \mathcal{T}_{(i)}$ for all $i\in \Lambda$. 
  
  $(i)\Leftarrow (ii)$ Since  $\im\mathcal{T}_{(i)}=T(i)$, the top of any direct summand of the characteristic tilting module of $A$ is simple.  According to Theorem 5.1 and   6.1 in \cite{P} the algebra  $(R(A),\ma_{(R)})$ is 1-quasi-hereditary. This implies  $\textbf{R}[B,L]\in \textbf{Y}(1)$.

 $(i) \Rightarrow (iii)$   Let   $i\in  \Lambda$ and $\Upsilon_{(i)}:B\sur L(i)$ be  the  surjective  $B$-map  given by   $\Upsilon_{(i)} (F)= F\circ \ger{f}_{(i)}$ for all $F\in B$.  Since 
   $\left(\mathcal{T}_{(j)}\circ \ger{f}_{(i)}\right):\left(P(1)\stackrel{f_{(i,i,1)}}{\rightarrow}P(i)\stackrel{f_{(1,i,i)}}{\hookrightarrow} P(1)\stackrel{\pi(j)}{\sur} T(j)\stackrel{\iota(j)}{\hookrightarrow} P(1)\right)$ (in the notations of Subsection  2.1)  and since  for any $j\in \Lambda\backslash\Lambda^{(i)}$  holds      $\dimm_K\Hom_A(P(i),T(j))=[T(j):S(i)]=0$ (because
of the  properties of $T(i)$), we    obtain 
    $\pi(j)\circ f_{(1,i,i)} =0$ and hence $\mathcal{T}_{(j)}\circ \ger{f}_{(i)}=0$. For all $i,j\in \Lambda$ with $i\not\ma j$ we have $\widetilde{R}(L(j)) = B\circ \mathcal{T}_{(j)} \subseteq \ker (\Upsilon_{(i)})$, thus 
\begin{center}
$\displaystyle \sum_{j\in \Lambda\backslash \Lambda^{(i)}}\widetilde{R}(L(j)) \subseteq \ker (\Upsilon_{(i)})$ \  \  and \  \   $\displaystyle \widetilde{R}(L(j)) \subseteq \bigcap_{i\in \Lambda\backslash \Lambda_{(j)}}\ker (\Upsilon_{(i)}) $
\end{center}
 By  assumption  the pairs $(B,L)$ and  $(B,\widetilde{R}(L))$ satisfy  the properties  \fbox{$\ma$} and   \fbox{$\ma_{(R)}$}  respectively. Therefore   $\dimm_KL(i)\stackrel{~\ref{L5}\textit{(1)}}{=}\left|\Lambda^{(i)}\right|$  implies  $\dimm_K\ker\left(\Upsilon_{(i)}\right)= \left|\Lambda\backslash\Lambda^{(i)}\right|$, moreover,    $  \dimm_K\left(\sum_{j\in \Lambda\backslash\Lambda^{(i)}}\widetilde{R}(L(j))\right) \stackrel{~\ref{L5}\textit{(1)}}{=} \left|\bigcup_{j\in \Lambda \atop j\not\ma_{(R)}i}\left\{k\in \Lambda \mid j\ma_{(R)}k\right\} \right|=  \left|\Lambda\backslash \Lambda^{(i)}\right|$.  This implies      $\sum_{j\in \Lambda\backslash \Lambda^{(i)}}\widetilde{R}(L(j))= \ker (\Upsilon_{(i)}) $  and consequently  $\bigcap_{i\in \Lambda\backslash \Lambda_{(j)}}\ker (\Upsilon_{(i)})= \sum_{j\ma_{(R)}k}\widetilde{R}(L(k))=\widetilde{R}(L(j))$, because $\widetilde{R}(L(k))\subseteq \widetilde{R}(L(j))$ for all $k\in \Lambda$ with $j\ma_{(R)}k$.  Moreover we obtain  $L(i) \cong B/\ker (\Upsilon_{(i)})  \cong B/\left(\sum_{j\in \Lambda\backslash \Lambda^{(i)}}\widetilde{R}(L(j)) \right) = B/\left(\sum_{j\not\ma_{(R)}i}\widetilde{R}(L(j)) \right)$. Using the dual argumentation we obtain $\widetilde{R}(L(i))\cong B/\left(\sum_{j\in \Lambda\backslash \Lambda_{(i)}}L(j)\right)$ for any $i\in \Lambda$.

$(iii)\Rightarrow (i)$ We  have to show   that  $(B,\widetilde{R}(L)=\bigoplus_{i\in \Lambda}\widetilde{R}(L(i)))$ satisfies the property \fbox{$\ma_{(R)}$}.\\
  Since $L(i)=B\circ \ger{f}_{(i)} \stackrel{~\ref{L5}\textit{(1)}}{=} \ger{f}_{(i)} \circ B$,  for  all   $f,g\in B$  there exists $f'\in B$  with   $g\circ f \circ \ger{f}_{(i)} = g\circ \ger{f}_{(i)}\circ f' $.   Therefore, if   $g\in \ker\left(B \stackrel{\circ \ger{f}_{(i)}}{\sur} L(i)\right)$, then  $g\circ f\in \ker\left(B \stackrel{\circ \ger{f}_{(i)}}{\sur} L(i)\right)$ for  all  $f\in B$. The assumption   $\widetilde{R}(L(j))\cong \bigcap_{i\in \Lambda\backslash\Lambda_{(i)}}\ker\left(B\sur L(i)\right) $  implies  $\widetilde{R}(L(j))\circ B\subseteq \widetilde{R}(L(j))$  and consequently    $\widetilde{R}(L(j))$ is two-sided local  ideal of $B$, since $\widetilde{R}(L(j))$ is a factor module and a submodule  of $B$.  
  
   The property \fbox{$\ma$} of $L(i)$  implies  $\widetilde{R}(L(k))\sur \widetilde{R}(L(j))$ resp. $\widetilde{R}(L(j)) \subseteq    \widetilde{R}(L(k))$  if and only  if $k\ma_{(R)} j$,   because  $\Lambda\backslash\Lambda_{(k)}\subseteq \Lambda\backslash\Lambda_{(j)}$ if and only if $j\ma k$.     We have  $\rad \widetilde{R}(L(k)) = \sum_{k<_{(R)}j}\widetilde{R}(L(j))$ for all $k\in \Lambda$ since $\widetilde{R}(L(k))/\left(\sum_{k<_{(R)}j}\widetilde{R}(L(j))\right)$ is simple.   \hfill $\Box$

\begin{subnum}
\begin{normalfont}\textbf{Remark.} If  $[B,L]\in \widehat{\mathcal{Y}}$, then the socle of $B/\left(\sum_{j\in \Lambda\backslash \Lambda_{(i)}}L(j)\right)$ is simple for all  $i\in \Lambda$, because these factor modules of $B$ are also  submodules of $B$.

 In  Example ~\ref{non_Ringel_dual} the pair $[B,L_2]\in \textbf{Y}(1)$  is not in $\widehat{\mathcal{Y}}$, because for $i=5$  we obtain  \begin{small}  $\soc \left(B/\sum_{j\in \Lambda\backslash \Lambda_{(5)}}L_2(j)\right)=\soc \left( B/L_2(6)\right) \cong (\soc^2B)/\soc B \cong \left(\left\langle X^3\right\rangle +\left\langle Y\right\rangle\right)/\left\langle X^4\right\rangle \cong K \oplus K$.\end{small}
\end{normalfont}
\end{subnum}

\begin{subnum}
\begin{normalfont}\textbf{Example.} Let $B:=B_n(C)$ be the algebra given in Example~\ref{Example_C}  and $L=B\oplus \bigoplus_{i=2}^{n-1}\left\langle X_i\right\rangle \oplus (\soc B)$.  The Ringel dual of  the corresponding 1-quasi-hereditary algebra $A_n(C)$  is also 1-quasi-hereditary (see \cite[Lemma 1.2]{P1}). Thus  $[B,L]\in  \widehat{\mathcal{Y}}$,   
   in particular,  $\widetilde{R}(L(1))\cong \soc B$, $\widetilde{R}(L(n))\cong B$ and $\widetilde{R}(L(i))\cong B/\left(\sum_{j=2\atop i\neq j}^{n-1}L(j)\right)\cong \left\langle \sum_{j=2}^{n-1}d_{ji}X_j\right\rangle$ for  $2\leq j\leq n-1$, where $C^{-1}$$=$$(d_{ij})_{2\leq i,j\leq n-1}$.
\end{normalfont}
\label{Example_matrix}
\end{subnum}

\begin{subnum}
\begin{normalfont}\textbf{Remark.} If a 1-quasi-hereditary algebra $A$  is Ringel self-dual, then   $(A,\ma)\cong (R(A),\ma_{(R)})$ implies   $[B,L]=\Phi([A])=\Phi([R(A)])=[B,\widetilde{R}(L)]$. In other words: The pair $[B,L]$ satisfies the property \fbox{$\ma$} and the property \fbox{$\ma_{(R)}$}.  Thus, there exists a permutation $\sigma\in \text{Sym}(\dimm_KB)$ with $L(\sigma(i))\cong \widetilde{R}(L(i)) $.

The algebras associated with blocks of the category $\OO(\ger{g})$ are Ringel self-dual.
In the Example~\ref{exampleB}  the 1-quasi-hereditary algebra  $(A_3,\ma_{(3)})$  corresponding to the pair $(B,L_3)$  is related to a regular block of $\OO(\ger{sl}_3)$.  It
is easy to check that    the   permutation $\sigma=(\sigma(1),\ldots , \sigma(6))= (6,5,4,3,2,1) \in \text{Sym}(6)$ yields     $L_3(\sigma(i))\cong \widetilde{R}(L_3(i))$. 

In the same example the algebra $(A_1,\ma_{(1)})$ is also Ringel self-dual. The permutation $\tau \in \text{Sym}(6)$ with $L_1(\tau (i))\cong \widetilde{R}(L_1(i))$ is given by $\tau=(6,4,5,2,3,1)$.
	The algebra $R(A_2)$ is not 1-quasi-hereditary, because $\soc (B/\left(\sum_{j\in \Lambda\backslash \Lambda_{(4)}}L_2(j)\right))\cong \soc\left(B/L_2(5)\right)$ is not simple.
 
In general, for some  $[B,L]\in \widehat{\mathcal{Y}}$ the equation $[B,L]= [B,\widetilde{R}(L)]$ is not satisfied (see Example~\ref{Example_matrix}).
Consequently  a 1-quasi-hereditary algebra $A$ with $[A]\in \widehat{\mathcal{X}}$ is not Ringel self-dual, in general. 
\end{normalfont}
\end{subnum}

\textbf{Acknowledgments.}  
I would like to thank   Jean-Marie 	Bois,  Rolf 	Farnsteiner,    Julian Külshammer  and Julia 	Worch    for their helpful and constructive comments.

\begin{small}

\end{small}

\end{document}